
\magnification=1200
\hsize=13.50cm
\vsize=18cm
\hoffset=-2mm
\voffset=.8cm
\parindent=12pt   \parskip=0pt
\hfuzz=1pt

\pretolerance=500 \tolerance=1000  \brokenpenalty=5000

\catcode`\@=11

\font\eightrm=cmr8
\font\eighti=cmmi8
\font\eightsy=cmsy8
\font\eightbf=cmbx8
\font\eighttt=cmtt8
\font\eightit=cmti8
\font\eightsl=cmsl8
\font\sevenrm=cmr7
\font\seveni=cmmi7
\font\sevensy=cmsy7
\font\sevenbf=cmbx7

\font\sixrm=cmr6
\font\sixi=cmmi6
\font\sixsy=cmsy6
\font\sixbf=cmbx6

\font\douzebf=cmbx10 at 12pt

\font\twelvebf=cmbx10 at 12pt

\font\tencal=eusm10

\font\sevencal=eusm7

\font\fivecal=eusm5
\newfam\calfam
\textfont\calfam=\tencal
\scriptfont\calfam=\sevencal
\scriptscriptfont\calfam=\fivecal
\def\cal#1{{\fam\calfam\relax#1}}

\skewchar\eighti='177 \skewchar\sixi='177
\skewchar\eightsy='60 \skewchar\sixsy='60

\def\tenpoint{%
   \textfont0=\tenrm \scriptfont0=\sevenrm
   \scriptscriptfont0=\fiverm
   \def\rm{\fam\z@\tenrm}%
   \textfont1=\teni  \scriptfont1=\seveni
   \scriptscriptfont1=\fivei
   \def\oldstyle{\fam\@ne\teni}\let\old=\oldstyle
   \textfont2=\tensy \scriptfont2=\sevensy
   \scriptscriptfont2=\fivesy
   \textfont\itfam=\tenit
   \def\it{\fam\itfam\tenit}%
   \textfont\slfam=\tensl
   \def\sl{\fam\slfam\tensl}%
   \textfont\bffam=\tenbf
   \scriptfont\bffam=\sevenbf
   \scriptscriptfont\bffam=\fivebf
   \def\bf{\fam\bffam\tenbf}%
   \textfont\ttfam=\tentt
   \def\tt{\fam\ttfam\tentt}%
   \abovedisplayskip=12pt plus 3pt minus 9pt
   \belowdisplayskip=\abovedisplayskip
   \abovedisplayshortskip=0pt plus 3pt
   \belowdisplayshortskip=4pt plus 3pt
   \smallskipamount=3pt plus 1pt minus 1pt
   \medskipamount=6pt plus 2pt minus 2pt
   \bigskipamount=12pt plus 4pt minus 4pt
   \normalbaselineskip=12pt
   \setbox\strutbox=\hbox{\vrule height8.5pt depth3.5pt width0pt}%
   \let\bigf@nt=\tenrm
   \let\smallf@nt=\sevenrm
   \normalbaselines\rm}

\def\eightpoint{%
   \textfont0=\eightrm \scriptfont0=\sixrm
   \scriptscriptfont0=\fiverm
   \def\rm{\fam\z@\eightrm}%
   \textfont1=\eighti  \scriptfont1=\sixi
   \scriptscriptfont1=\fivei
   \def\oldstyle{\fam\@ne\eighti}\let\old=\oldstyle
   \textfont2=\eightsy \scriptfont2=\sixsy
   \scriptscriptfont2=\fivesy
   \textfont\itfam=\eightit
   \def\it{\fam\itfam\eightit}%
   \textfont\slfam=\eightsl
   \def\sl{\fam\slfam\eightsl}%
   \textfont\bffam=\eightbf
   \scriptfont\bffam=\sixbf
   \scriptscriptfont\bffam=\fivebf
   \def\bf{\fam\bffam\eightbf}%
   \textfont\ttfam=\eighttt
   \def\tt{\fam\ttfam\eighttt}%
   \abovedisplayskip=9pt plus 3pt minus 9pt
   \belowdisplayskip=\abovedisplayskip
   \abovedisplayshortskip=0pt plus 3pt
   \belowdisplayshortskip=3pt plus 3pt
   \smallskipamount=2pt plus 1pt minus 1pt
   \medskipamount=4pt plus 2pt minus 1pt
   \bigskipamount=9pt plus 3pt minus 3pt
   \normalbaselineskip=9pt
   \setbox\strutbox=\hbox{\vrule height7pt depth2pt width0pt}%
   \let\bigf@nt=\eightrm
   \let\smallf@nt=\sixrm
   \normalbaselines\rm}
\tenpoint


\def\pc#1{\bigf@nt#1\smallf@nt}

\catcode`\;=\active
\def;{\relax\ifhmode\ifdim\lastskip>\z@\unskip\fi \kern\fontdimen2
 -1.2 \fontdimen3 \string;}

\catcode`\:=\active
\def:{\relax\ifhmode\ifdim\lastskip>\z@\unskip\fi\penalty\@M\
\fi\string:}

\catcode`\!=\active
\def!{\relax\ifhmode\ifdim\lastskip>\z@ \unskip\fi\kern\fontdimen2
 -1.1 \fontdimen3 \string!}

\catcode`\?=\active
\def?{\relax\ifhmode\ifdim\lastskip>\z@ \unskip\fi\kern\fontdimen2
 -1.1 \fontdimen3 \string?}

\frenchspacing

\newtoks\auteurcourant
\auteurcourant={\hfil}

\newtoks\titrecourant
\titrecourant={\hfil}

\newtoks\hautpagetitre
\hautpagetitre={\hfil}

\newtoks\baspagetitre
\baspagetitre={\hfil}

\newtoks\hautpagegauche
\hautpagegauche={\eightpoint\rlap{\folio}\hfil\the
\auteurcourant\hfil}

\newtoks\hautpagedroite
\hautpagedroite={\eightpoint\hfil\the
\titrecourant\hfil\llap{\folio}}

\newtoks\baspagegauche
\baspagegauche={\hfil}

\newtoks\baspagedroite
\baspagedroite={\hfil}

\newif\ifpagetitre
\pagetitretrue


\headline={\ifpagetitre\the\hautpagetitre
\else\ifodd\pageno\the\hautpagedroite\else\the
\hautpagegauche\fi\fi}

\footline={\ifpagetitre\the\baspagetitre\else
\ifodd\pageno\the\baspagedroite\else\the
\baspagegauche\fi\fi
\global\pagetitrefalse}

\def\raggedbottom{\topskip 10pt plus 36pt\r@ggedbottomtrue}

\def\point{\raise.2ex\hbox{\douzebf .}}
\def\pointir{\unskip . --- \ignorespaces}


\def\Medbreak{\vskip-\lastskip\medbreak}

\def\rem#1\endrem{%
\Medbreak {\it#1\unskip} : }


\long\def\th#1 #2\enonce#3\endth{%
\Medbreak {\pc#1} {#2\unskip}\pointir{\it #3}\medskip}


\long\def\thm#1 #2\enonce#3\endthm{%
\Medbreak {\pc#1} {#2\unskip}\pointir{\it #3}\medskip}

\def\decale#1{\smallbreak\hskip 28pt\llap{#1}\kern 5pt}

\def\decaledecale#1{\smallbreak\hskip 34pt\llap{#1}\kern 5pt}

\let\@ldmessage=\message

\def\message#1{{\def\pc{\string\pc\space}%
\def\'{\string'}\def\`{\string`}%
\def\^{\string^}\def\"{\string"}%
\@ldmessage{#1}}}

\def\({{\rm (}}

\def\){{\rm )}}


\def\up#1{\raise 1ex\hbox{\smallf@nt#1}}

\def\diagram#1{\def\normalbaselines{\baselineskip=0pt
\lineskip=5pt}\matrix{#1}}


\def\longmaprightover#1#2{\smash{\mathop{\hbox to#2
{\rightarrowfill}}\limits^{\scriptstyle#1}}}

\def\longmapleftover#1#2{\smash{\mathop{\hbox to#2
{\leftarrowfill}}\limits^{\scriptstyle#1}}}

\def\longmaprightunder#1#2{\smash{\mathop{\hbox to#2
{\rightarrowfill}}\limits_{\scriptstyle#1}}}

\def\longmapleftunder#1#2{\smash{\mathop{\hbox to#2
{\leftarrowfill}}\limits_{\scriptstyle#1}}}

\def\longhookrightarrowover#1#2{\smash{\mathop{\lhook\joinrel
\mathrel{\hbox to #2{\rightarrowfill}}}\limits^{\scriptstyle#1}}}

\def\longhookrightarrowunder#1#2{\smash{\mathop{\lhook\joinrel
\mathrel{\hbox to #2{\rightarrowfill}}}\limits_{\scriptstyle#1}}}

\def\longhookleftarrowover#1#2{\smash{\mathop{{\hbox to #2
{\leftarrowfill}}\joinrel\kern -0.9mm\mathrel\rhook}
\limits^{\scriptstyle#1}}}

\def\longhookleftarrowunder#1#2{\smash{\mathop{{\hbox to #2
{\leftarrowfill}}\joinrel\kern -0.9mm\mathrel\rhook}
\limits_{\scriptstyle#1}}}

\def\longtwoheadrightarrowover#1#2{\smash{\mathop{{\hbox to #2
{\rightarrowfill}}\kern -3.25mm\joinrel\mathrel\rightarrow}
\limits^{\scriptstyle#1}}}

\def\longtwoheadrightarrowunder#1#2{\smash{\mathop{{\hbox to #2
{\rightarrowfill}}\kern -3.25mm\joinrel\mathrel\rightarrow}
\limits_{\scriptstyle#1}}}

\def\longtwoheadleftarrowover#1#2{\smash{\mathop{\joinrel\mathrel
\leftarrow\kern -3.8mm{\hbox to #2{\leftarrowfill}}}
\limits^{\scriptstyle#1}}}

\def\longtwoheadleftarrowunder#1#2{\smash{\mathop{\joinrel\mathrel
\leftarrow\kern -3.8mm{\hbox to #2{\leftarrowfill}}}
\limits_{\scriptstyle#1}}}


\def\longmapsto#1{\mapstochar\mathrel{\joinrel \kern-0.2mm\hbox to
#1mm{\rightarrowfill}}}

\def\og{\leavevmode\raise.3ex\hbox{$\scriptscriptstyle
\langle\!\langle\,$}}
\def\fg{\leavevmode\raise.3ex\hbox{$\scriptscriptstyle
\,\rangle\!\rangle$}}

\def\section#1#2{\vskip 5mm {\bf {#1}. {#2}}\vskip 5mm}
\def\subsection#1#2{\vskip 3mm {\it #2}\vskip 3mm}

\catcode`\@=12

\showboxbreadth=-1  \showboxdepth=-1


\message{`lline' & `vector' macros from LaTeX}

\def\Grille{\setbox13=\vbox to 5\unitlength{\hrule width 109mm \vfill}
\setbox13=\vbox to 65mm
{\offinterlineskip\leaders\copy13\vfill\kern-1pt\hrule}
\setbox14=\hbox to 5\unitlength{\vrule height 65mm\hfill}
\setbox14=\hbox to 109mm{\leaders\copy14\hfill\kern-2mm \vrule height
65mm}
\ht14=0pt\dp14=0pt\wd14=0pt \setbox13=\vbox to 0pt
{\vss\box13\offinterlineskip\box14} \wd13=0pt\box13}

\def\rule(#1,#2)\dir(#3,#4)\long#5{%
\noalign{\leftput(#1,#2){\lline(#3,#4){#5}}}}

\def\arrow(#1,#2)\dir(#3,#4)\length#5{%
\noalign{\leftput(#1,#2){\vector(#3,#4){#5}}}}

\def\put(#1,#2)#3{\noalign{\setbox1=\hbox{%
\kern #1\unitlength \raise #2\unitlength\hbox{$#3$}}%
\ht1=0pt \wd1=0pt \dp1=0pt\box1}}

\catcode`@=11

\def\{{\relax\ifmmode\lbrace\else$\lbrace$\fi}

\def\}{\relax\ifmmode\rbrace\else$\rbrace$\fi}

\def\newcount{\alloc@0\count\countdef\insc@unt}

\def\newdimen{\alloc@1\dimen\dimendef\insc@unt}

\def\newwrite{\alloc@7\write\chardef\sixt@@n}

\newwrite\@unused
\newcount\@tempcnta
\newcount\@tempcntb
\newdimen\@tempdima
\newdimen\@tempdimb
\newbox\@tempboxa

\def\@spaces{\space\space\space\space}

\def\@whilenoop#1{}

\def\@whiledim#1\do #2{\ifdim #1\relax#2\@iwhiledim{#1\relax#2}\fi}

\def\@iwhiledim#1{\ifdim #1\let\@nextwhile=\@iwhiledim
\else\let\@nextwhile=\@whilenoop\fi\@nextwhile{#1}}

\def\@badlinearg{\@latexerr{Bad \string\line\space or \string\vector
\space argument}}

\def\@latexerr#1#2{\begingroup
\edef\@tempc{#2}\expandafter\errhelp\expandafter{\@tempc}%

\def\@eha{Your command was ignored.^^JType \space I <command> <return>
\space to replace it with another command,^^Jor \space <return> \space
to continue without it.}

\def\@ehb{You've lost some text. \space \@ehc}

\def\@ehc{Try typing \space <return> \space to proceed.^^JIf that
doesn't work, type \space X <return> \space to quit.}

\def\@ehd{You're in trouble here.  \space\@ehc}

\typeout{LaTeX error.  \space See LaTeX manual for explanation.^^J
\space\@spaces\@spaces\@spaces Type \space H <return> \space for
immediate help.}\errmessage{#1}\endgroup}

\def\typeout#1{{\let\protect\string\immediate\write\@unused{#1}}}

\font\tenln = line10
\font\tenlnw = linew10

\newdimen\@wholewidth
\newdimen\@halfwidth
\newdimen\unitlength

\unitlength =1pt

\def\thinlines{\let\@linefnt\tenln \let\@circlefnt\tencirc
\@wholewidth\fontdimen8\tenln \@halfwidth .5\@wholewidth}

\def\thicklines{\let\@linefnt\tenlnw \let\@circlefnt\tencircw
\@wholewidth\fontdimen8\tenlnw \@halfwidth .5\@wholewidth}

\def\linethickness#1{\@wholewidth #1\relax \@halfwidth .5
\@wholewidth}

\newif\if@negarg

\def\lline(#1,#2)#3{\@xarg #1\relax \@yarg #2\relax
\@linelen=#3\unitlength \ifnum\@xarg =0 \@vline \else \ifnum\@yarg =0
\@hline \else \@sline\fi \fi}

\def\@sline{\ifnum\@xarg< 0 \@negargtrue \@xarg -\@xarg \@yyarg
-\@yarg \else \@negargfalse \@yyarg \@yarg \fi
\ifnum \@yyarg >0 \@tempcnta\@yyarg \else \@tempcnta - \@yyarg \fi
\ifnum\@tempcnta>6 \@badlinearg\@tempcnta0 \fi
\setbox\@linechar\hbox{\@linefnt\@getlinechar(\@xarg,\@yyarg)}%
\ifnum \@yarg >0 \let\@upordown\raise \@clnht\z@
\else\let\@upordown\lower \@clnht \ht\@linechar\fi
\@clnwd=\wd\@linechar
\if@negarg \hskip -\wd\@linechar \def\@tempa{\hskip -2\wd \@linechar}
\else \let\@tempa\relax \fi
\@whiledim \@clnwd <\@linelen \do {\@upordown\@clnht\copy\@linechar
\@tempa \advance\@clnht \ht\@linechar \advance\@clnwd \wd\@linechar}%
\advance\@clnht -\ht\@linechar \advance\@clnwd -\wd\@linechar
\@tempdima\@linelen\advance\@tempdima -\@clnwd
\@tempdimb\@tempdima\advance\@tempdimb -\wd\@linechar
\if@negarg \hskip -\@tempdimb \else \hskip \@tempdimb \fi
\multiply\@tempdima \@m\@tempcnta \@tempdima \@tempdima \wd\@linechar
\divide\@tempcnta \@tempdima \@tempdima \ht\@linechar
\multiply\@tempdima \@tempcnta \divide\@tempdima \@m \advance\@clnht
\@tempdima
\ifdim \@linelen <\wd\@linechar \hskip \wd\@linechar
\else\@upordown\@clnht\copy\@linechar\fi}

\def\@hline{\ifnum \@xarg <0 \hskip -\@linelen \fi
\vrule height \@halfwidth depth \@halfwidth width \@linelen
\ifnum \@xarg <0 \hskip -\@linelen \fi}

\def\@getlinechar(#1,#2){\@tempcnta#1\relax
\multiply\@tempcnta 8\advance\@tempcnta -9
\ifnum #2>0 \advance\@tempcnta #2\relax
  \else\advance\@tempcnta -#2\relax\advance\@tempcnta 64 \fi
\char\@tempcnta}

\def\vector(#1,#2)#3{\@xarg #1\relax \@yarg #2\relax
\@linelen=#3\unitlength
\ifnum\@xarg =0 \@vvector \else \ifnum\@yarg =0 \@hvector \else
\@svector\fi \fi}

\def\@hvector{\@hline\hbox to 0pt{\@linefnt \ifnum \@xarg <0
\@getlarrow(1,0)\hss\else \hss\@getrarrow(1,0)\fi}}

\def\@vvector{\ifnum \@yarg <0 \@downvector \else \@upvector \fi}

\def\@svector{\@sline\@tempcnta\@yarg \ifnum\@tempcnta <0
\@tempcnta=-\@tempcnta\fi \ifnum\@tempcnta <5 \hskip -\wd\@linechar
\@upordown\@clnht \hbox{\@linefnt \if@negarg
\@getlarrow(\@xarg,\@yyarg) \else \@getrarrow(\@xarg,\@yyarg)
\fi}\else\@badlinearg\fi}

\def\@getlarrow(#1,#2){\ifnum #2 =\z@ \@tempcnta='33\else
\@tempcnta=#1\relax\multiply\@tempcnta \sixt@@n \advance\@tempcnta -9
\@tempcntb=#2\relax \multiply\@tempcntb \tw@ \ifnum \@tempcntb >0
\advance\@tempcnta \@tempcntb\relax \else\advance\@tempcnta
-\@tempcntb\advance\@tempcnta 64 \fi\fi \char\@tempcnta}

\def\@getrarrow(#1,#2){\@tempcntb=#2\relax \ifnum\@tempcntb < 0
\@tempcntb=-\@tempcntb\relax\fi \ifcase \@tempcntb\relax
\@tempcnta='55 \or \ifnum #1<3 \@tempcnta=#1\relax\multiply\@tempcnta
24 \advance\@tempcnta -6 \else \ifnum #1=3 \@tempcnta=49
\else\@tempcnta=58 \fi\fi\or \ifnum #1<3
\@tempcnta=#1\relax\multiply\@tempcnta 24 \advance\@tempcnta -3 \else
\@tempcnta=51\fi\or \@tempcnta=#1\relax\multiply\@tempcnta \sixt@@n
\advance\@tempcnta -\tw@ \else \@tempcnta=#1\relax\multiply\@tempcnta
\sixt@@n \advance\@tempcnta 7 \fi \ifnum #2<0 \advance\@tempcnta 64
\fi \char\@tempcnta}

\def\@vline{\ifnum \@yarg <0 \@downline \else \@upline\fi}

\def\@upline{\hbox to \z@{\hskip -\@halfwidth \vrule width
\@wholewidth height \@linelen depth \z@\hss}}

\def\@downline{\hbox to \z@{\hskip -\@halfwidth \vrule width
\@wholewidth height \z@ depth \@linelen \hss}}

\def\@upvector{\@upline\setbox\@tempboxa
\hbox{\@linefnt\char'66}\raise \@linelen \hbox to\z@{\lower
\ht\@tempboxa \box\@tempboxa\hss}}

\def\@downvector{\@downline\lower \@linelen \hbox to
\z@{\@linefnt\char'77\hss}}

\thinlines

\newcount\@xarg
\newcount\@yarg
\newcount\@yyarg
\newcount\@multicnt
\newdimen\@xdim
\newdimen\@ydim
\newbox\@linechar
\newdimen\@linelen
\newdimen\@clnwd
\newdimen\@clnht
\newdimen\@dashdim
\newbox\@dashbox
\newcount\@dashcnt
\catcode`@=12

\newbox\tbox
\newbox\tboxa

\def\leftzer#1{\setbox\tbox=\hbox to 0pt{#1\hss}%
\ht\tbox=0pt \dp\tbox=0pt \box\tbox}

\def\rightzer#1{\setbox\tbox=\hbox to 0pt{\hss #1}%
\ht\tbox=0pt \dp\tbox=0pt \box\tbox}

\def\centerzer#1{\setbox\tbox=\hbox to 0pt{\hss #1\hss}%
\ht\tbox=0pt \dp\tbox=0pt \box\tbox}

\def\leftput(#1,#2)#3{\setbox\tboxa=\hbox{%
\kern #1\unitlength \raise #2\unitlength\hbox{\leftzer{#3}}}%
\ht\tboxa=0pt \wd\tboxa=0pt \dp\tboxa=0pt\box\tboxa}

\def\rightput(#1,#2)#3{\setbox\tboxa=\hbox{%
\kern #1\unitlength \raise #2\unitlength\hbox{\rightzer{#3}}}%
\ht\tboxa=0pt \wd\tboxa=0pt \dp\tboxa=0pt\box\tboxa}

\def\centerput(#1,#2)#3{\setbox\tboxa=\hbox{%
\kern #1\unitlength \raise #2\unitlength\hbox{\centerzer{#3}}}%
\ht\tboxa=0pt \wd\tboxa=0pt \dp\tboxa=0pt\box\tboxa}

\unitlength=1mm

\expandafter\ifx\csname amssym.def\endcsname\relax \else
\endinput\fi
%
\expandafter\edef\csname amssym.def\endcsname{%
        \catcode`\noexpand\@=\the\catcode`\@\space}
\catcode`\@=11
%

\def\undefine#1{\let#1\undefined}
\def\newsymbol#1#2#3#4#5{\let\next@\relax
  \ifnum#2=\@ne\let\next@\msafam@\else
  \ifnum#2=\tw@\let\next@\msbfam@\fi\fi
  \mathchardef#1="#3\next@#4#5}
\def\mathhexbox@#1#2#3{\relax
  \ifmmode\mathpalette{}{\m@th\mathchar"#1#2#3}%
  \else\leavevmode\hbox{$\m@th\mathchar"#1#2#3$}\fi}
\def\hexnumber@#1{\ifcase#1 0\or 1\or 2\or 3
\or 4\or 5\or 6\or 7\or 8\or
  9\or A\or B\or C\or D\or E\or F\fi}

\font\tenmsa=msam10
\font\sevenmsa=msam7
\font\fivemsa=msam5
\newfam\msafam
\textfont\msafam=\tenmsa
\scriptfont\msafam=\sevenmsa
\scriptscriptfont\msafam=\fivemsa
\edef\msafam@{\hexnumber@\msafam}
\mathchardef\dabar@"0\msafam@39
\def\dashrightarrow{\mathrel{\dabar@\dabar@\mathchar"0
\msafam@4B}}
\def\dashleftarrow{\mathrel{\mathchar"0\msafam@4C
\dabar@\dabar@}}

\def\ulcorner{\delimiter"4\msafam@70\msafam@70 }
\def\urcorner{\delimiter"5\msafam@71\msafam@71 }
\def\llcorner{\delimiter"4\msafam@78\msafam@78 }
\def\lrcorner{\delimiter"5\msafam@79\msafam@79 }
\def\yen{{\mathhexbox@\msafam@55}}
\def\checkmark{{\mathhexbox@\msafam@58}}
\def\circledR{{\mathhexbox@\msafam@72}}
\def\maltese{{\mathhexbox@\msafam@7A}}

\font\tenmsb=msbm10
\font\sevenmsb=msbm7
\font\fivemsb=msbm5
\newfam\msbfam
\textfont\msbfam=\tenmsb
\scriptfont\msbfam=\sevenmsb
\scriptscriptfont\msbfam=\fivemsb
\edef\msbfam@{\hexnumber@\msbfam}
\def\Bbb#1{{\fam\msbfam\relax#1}}
\def\widehat#1{\setbox\z@\hbox{$\m@th#1$}%
  \ifdim\wd\z@>\tw@ em\mathaccent"0\msbfam@5B{#1}%
  \else\mathaccent"0362{#1}\fi}
\def\widetilde#1{\setbox\z@\hbox{$\m@th#1$}%
  \ifdim\wd\z@>\tw@ em\mathaccent"0\msbfam@5D{#1}%
  \else\mathaccent"0365{#1}\fi}
\font\teneufm=eufm10
\font\seveneufm=eufm7
\font\fiveeufm=eufm5
\newfam\eufmfam
\textfont\eufmfam=\teneufm
\scriptfont\eufmfam=\seveneufm
\scriptscriptfont\eufmfam=\fiveeufm
\def\frak#1{{\fam\eufmfam\relax#1}}

\csname amssym.def\endcsname

\expandafter\ifx\csname pre amssym.tex at\endcsname\relax \else
\endinput\fi
\expandafter\chardef\csname pre amssym.tex at\endcsname=\the
\catcode`\@
\catcode`\@=11
\begingroup\ifx\undefined\newsymbol \else\def\input#1
{\endgroup}\fi
\input amssym.def \relax
\newsymbol\boxdot 1200
\newsymbol\boxplus 1201
\newsymbol\boxtimes 1202
\newsymbol\square 1003
\newsymbol\blacksquare 1004
\newsymbol\centerdot 1205
\newsymbol\lozenge 1006
\newsymbol\blacklozenge 1007
\newsymbol\circlearrowright 1308
\newsymbol\circlearrowleft 1309
\undefine\rightleftharpoons
\newsymbol\rightleftharpoons 130A
\newsymbol\leftrightharpoons 130B
\newsymbol\boxminus 120C
\newsymbol\Vdash 130D
\newsymbol\Vvdash 130E
\newsymbol\vDash 130F
\newsymbol\twoheadrightarrow 1310
\newsymbol\twoheadleftarrow 1311
\newsymbol\leftleftarrows 1312
\newsymbol\rightrightarrows 1313
\newsymbol\upuparrows 1314
\newsymbol\downdownarrows 1315
\newsymbol\upharpoonright 1316
  
\newsymbol\downharpoonright 1317
\newsymbol\upharpoonleft 1318
\newsymbol\downharpoonleft 1319
\newsymbol\rightarrowtail 131A
\newsymbol\leftarrowtail 131B
\newsymbol\leftrightarrows 131C
\newsymbol\rightleftarrows 131D
\newsymbol\Lsh 131E
\newsymbol\Rsh 131F
\newsymbol\rightsquigarrow 1320
\newsymbol\leftrightsquigarrow 1321
\newsymbol\looparrowleft 1322
\newsymbol\looparrowright 1323
\newsymbol\circeq 1324
\newsymbol\succsim 1325
\newsymbol\gtrsim 1326
\newsymbol\gtrapprox 1327
\newsymbol\multimap 1328
\newsymbol\therefore 1329
\newsymbol\because 132A
\newsymbol\doteqdot 132B
  
\newsymbol\triangleq 132C
\newsymbol\precsim 132D
\newsymbol\lesssim 132E
\newsymbol\lessapprox 132F
\newsymbol\eqslantless 1330
\newsymbol\eqslantgtr 1331
\newsymbol\curlyeqprec 1332
\newsymbol\curlyeqsucc 1333
\newsymbol\preccurlyeq 1334
\newsymbol\leqq 1335
\newsymbol\leqslant 1336
\newsymbol\lessgtr 1337
\newsymbol\backprime 1038
\newsymbol\risingdotseq 133A
\newsymbol\fallingdotseq 133B
\newsymbol\succcurlyeq 133C
\newsymbol\geqq 133D
\newsymbol\geqslant 133E
\newsymbol\gtrless 133F
\newsymbol\sqsubset 1340
\newsymbol\sqsupset 1341
\newsymbol\vartriangleright 1342
\newsymbol\vartriangleleft 1343
\newsymbol\trianglerighteq 1344
\newsymbol\trianglelefteq 1345
\newsymbol\bigstar 1046
\newsymbol\between 1347
\newsymbol\blacktriangledown 1048
\newsymbol\blacktriangleright 1349
\newsymbol\blacktriangleleft 134A
\newsymbol\vartriangle 134D
\newsymbol\blacktriangle 104E
\newsymbol\triangledown 104F
\newsymbol\eqcirc 1350
\newsymbol\lesseqgtr 1351
\newsymbol\gtreqless 1352
\newsymbol\lesseqqgtr 1353
\newsymbol\gtreqqless 1354
\newsymbol\Rrightarrow 1356
\newsymbol\Lleftarrow 1357
\newsymbol\veebar 1259
\newsymbol\barwedge 125A
\newsymbol\doublebarwedge 125B
\undefine\angle
\newsymbol\angle 105C
\newsymbol\measuredangle 105D
\newsymbol\sphericalangle 105E
\newsymbol\varpropto 135F
\newsymbol\smallsmile 1360
\newsymbol\smallfrown 1361
\newsymbol\Subset 1362
\newsymbol\Supset 1363
\newsymbol\Cup 1264
  
\newsymbol\Cap 1265
  
\newsymbol\curlywedge 1266
\newsymbol\curlyvee 1267
\newsymbol\leftthreetimes 1268
\newsymbol\rightthreetimes 1269
\newsymbol\subseteqq 136A
\newsymbol\supseteqq 136B
\newsymbol\bumpeq 136C
\newsymbol\Bumpeq 136D
\newsymbol\lll 136E
  
\newsymbol\ggg 136F
  
\newsymbol\circledS 1073
\newsymbol\pitchfork 1374
\newsymbol\dotplus 1275
\newsymbol\backsim 1376
\newsymbol\backsimeq 1377
\newsymbol\complement 107B
\newsymbol\intercal 127C
\newsymbol\circledcirc 127D
\newsymbol\circledast 127E
\newsymbol\circleddash 127F
\newsymbol\lvertneqq 2300
\newsymbol\gvertneqq 2301
\newsymbol\nleq 2302
\newsymbol\ngeq 2303
\newsymbol\nless 2304
\newsymbol\ngtr 2305
\newsymbol\nprec 2306
\newsymbol\nsucc 2307
\newsymbol\lneqq 2308
\newsymbol\gneqq 2309
\newsymbol\nleqslant 230A
\newsymbol\ngeqslant 230B
\newsymbol\lneq 230C
\newsymbol\gneq 230D
\newsymbol\npreceq 230E
\newsymbol\nsucceq 230F
\newsymbol\precnsim 2310
\newsymbol\succnsim 2311
\newsymbol\lnsim 2312
\newsymbol\gnsim 2313
\newsymbol\nleqq 2314
\newsymbol\ngeqq 2315
\newsymbol\precneqq 2316
\newsymbol\succneqq 2317
\newsymbol\precnapprox 2318
\newsymbol\succnapprox 2319
\newsymbol\lnapprox 231A
\newsymbol\gnapprox 231B
\newsymbol\nsim 231C
\newsymbol\ncong 231D
\newsymbol\diagup 201E
\newsymbol\diagdown 201F
\newsymbol\varsubsetneq 2320
\newsymbol\varsupsetneq 2321
\newsymbol\nsubseteqq 2322
\newsymbol\nsupseteqq 2323
\newsymbol\subsetneqq 2324
\newsymbol\supsetneqq 2325
\newsymbol\varsubsetneqq 2326
\newsymbol\varsupsetneqq 2327
\newsymbol\subsetneq 2328
\newsymbol\supsetneq 2329
\newsymbol\nsubseteq 232A
\newsymbol\nsupseteq 232B
\newsymbol\nparallel 232C
\newsymbol\nmid 232D
\newsymbol\nshortmid 232E
\newsymbol\nshortparallel 232F
\newsymbol\nvdash 2330
\newsymbol\nVdash 2331
\newsymbol\nvDash 2332
\newsymbol\nVDash 2333
\newsymbol\ntrianglerighteq 2334
\newsymbol\ntrianglelefteq 2335
\newsymbol\ntriangleleft 2336
\newsymbol\ntriangleright 2337
\newsymbol\nleftarrow 2338
\newsymbol\nrightarrow 2339
\newsymbol\nLeftarrow 233A
\newsymbol\nRightarrow 233B
\newsymbol\nLeftrightarrow 233C
\newsymbol\nleftrightarrow 233D
\newsymbol\divideontimes 223E
\newsymbol\varnothing 203F
\newsymbol\nexists 2040
\newsymbol\Finv 2060
\newsymbol\Game 2061
\newsymbol\mho 2066
\newsymbol\eth 2067
\newsymbol\eqsim 2368
\newsymbol\beth 2069
\newsymbol\gimel 206A
\newsymbol\daleth 206B
\newsymbol\lessdot 236C
\newsymbol\gtrdot 236D
\newsymbol\ltimes 226E
\newsymbol\rtimes 226F
\newsymbol\shortmid 2370
\newsymbol\shortparallel 2371
\newsymbol\smallsetminus 2272
\newsymbol\thicksim 2373
\newsymbol\thickapprox 2374
\newsymbol\approxeq 2375
\newsymbol\succapprox 2376
\newsymbol\precapprox 2377
\newsymbol\curvearrowleft 2378
\newsymbol\curvearrowright 2379
\newsymbol\digamma 207A
\newsymbol\varkappa 207B
\newsymbol\Bbbk 207C
\newsymbol\hslash 207D
\undefine\hbar
\newsymbol\hbar 207E
\newsymbol\backepsilon 237F
\catcode`\@=\csname pre amssym.tex at\endcsname


\centerline{\twelvebf Fibres de Springer et jacobiennes
compactifi\'{e}es}
\vskip 10mm
\centerline{G\'{e}rard Laumon\footnote{${}^{\ast}$}{\sevenrm CNRS et
Universit\'{e} Paris-Sud, UMR 8628, Math\'{e}matique, B\^{a}t. 425,
F-91405 Orsay Cedex, France, Gerard.Laumon@math.u-psud.fr}}
\vskip 25mm

L'objet de ce travail est d'identifier, \`{a} hom\'{e}omorphisme
pr\`{e}s, les fibres de Springer affines pour $\mathop{\rm GL}(n)$ sur
un corps local d'\'{e}gales caract\'{e}ristiques \`{a} des
rev\^{e}tements de jacobiennes compactifi\'{e}es de courbes
projectives singuli\`{e}res.  Ce lien permet de d\'{e}montrer
certaines propri\'{e}t\'{e}s g\'{e}om\'{e}triques de ces fibres de
Springer, dont une propri\'{e}t\'{e} d'irr\'{e}ductibilit\'{e}, et
aussi d'en construire des {\og}{d\'{e}formations}{\fg}.
\vskip 10mm

\centerline{0. INTRODUCTION}
\vskip 5mm

Soient $F$ un corps local non archim\'{e}dien d'\'{e}gales
caract\'{e}ristiques, ${\cal O}_{F}$ son anneau des entiers, $k$ son
corps r\'{e}siduel et $E$ un $F$-espace vectoriel de dimension finie.

La grassmannienne affine pour le $F$-sch\'{e}ma en groupes
$\mathop{\rm Aut}\nolimits_{F}(E)$ des automorphismes de $E$ est le
ind-$k$-sch\'{e}ma des ${\cal O}_{F}$-r\'{e}seaux $M$ dans $E$.  Pour
tout endomorphisme r\'{e}gulier semi-simple et topologiquement
nilpotent $\gamma$ de $E$, on peut consid\'{e}rer le ferm\'{e}
r\'{e}duit $X_{\gamma}$ de la grassmannienne affine form\'{e} des
r\'{e}seaux $M\subset E$ tels que $\gamma (M)\subset M$.  Ce ferm\'{e}
est appel\'{e} la {\it fibre de Springer affine en} $\gamma$ par
analogie avec les fibres de Springer classiques dans les
vari\'{e}t\'{e}s de drapeaux.  D'apr\`{e}s Kazhdan et Lusztig [K-L],
$X_{\gamma}$ est un vrai sch\'{e}ma, localement de type fini sur $k$,
qui est muni d'une action libre naturelle d'un groupe ab\'{e}lien
libre de type fini $\Lambda_{\gamma}$, et le quotient $Z_{\gamma}=
X_{\gamma}/\Lambda_{\gamma}$ est lui un $k$-sch\'{e}ma projectif.

Dans ce travail, nous attachons \`{a} $\gamma$ une courbe int\`{e}gre
et projective $C_{\gamma}$ sur $k$, qui n'a au plus qu'un seul point
singulier, point en lequel l'anneau local compl\'{e}t\'{e} de la
courbe n'est autre que ${\cal O}_{F}[\gamma]\subset F[\gamma] \subset
\mathop{\rm Au}\nolimits_{F}(E)$.  Puis, nous relions la fibre de
Springer affine $X_{\gamma}$ et son quotient $Z_{\gamma}$ \`{a} la
jacobienne compactifi\'{e}e de $C_{\gamma}$.  Nous d\'{e}duisons alors
des r\'{e}sultats d'Altman et Kleiman sur les jacobiennes
compactifi\'{e}es, un \'{e}nonc\'{e} d'irr\'{e}ductibilit\'{e} pour
$Z_{\gamma}$ (Corollaire 2.2.2) et la possibilit\'{e} de d\'{e}former
\`{a} hom\'{e}omorphisme pr\`{e}s $Z_{\gamma}$ (et aussi dans une
certaine mesure $X_{\gamma}$) en faisant varier $\gamma$ (cf.
Chapitre 4).

Cette \'{e}tude des fibres de Springer affines est motiv\'{e}e par une
conjecture de Langlands et Shelstad, connue sous le nom de {\og}{Lemme
fondamental}{\fg}, dans le cas particulier des groupes unitaires sur
les corps locaux non archim\'{e}diens.  Nous explicitons cette
conjecture dans la section 1.3 et nous en donnons, dans la section
1.4, une version g\'{e}om\'{e}trique, sous une forme due \`{a}
Kottwitz.

Pour simplifier l'exposition, nous nous somme partout limit\'{e} au cas
o\`{u} la $F$-alg\`{e}bre commutative, semi-simple et de dimension
finie, $F[\gamma]$, est un produit d'extensions {\it totalement
ramifi\'{e}es} de $F$.
\vskip 2mm

{\sevenrm Je remercie J.-B. Bost, E. Esteves, L. G\"{o}ttsche, L. Illusie, S.
Kleiman, Y. Laszlo, F. Loeser et M. Raynaud pour l'aide qu'ils m'ont
apport\'{e}e durant la pr\'{e}paration de ce travail.}
\vskip 5mm

\centerline{1. FIBRES DE SPRINGER}
\vskip 2mm

\section{1.1}{Les donn\'{e}es}

On fixe un corps $k$ parfait.  Dans cet travail, on appelle simplement
{\it corps local} tout corps $K$ contenant $k$ et muni d'une valuation
discr\`{e}te $v_{K}:K\rightarrow {\Bbb Z}\cup\{\infty\}$ pour laquelle
$K$ est complet et de corps r\'{e}siduel $k$.  Pour un tel corps
local, on note ${\cal O}_{K}\subset K$ l'anneau des entiers de la
valuation discr\`{e}te et ${\frak p}_{K}$ l'id\'{e}al maximal de
${\cal O}_{K}$.  Le choix d'une uniformisante $\varpi_{K}$ de $K$
identifie ${\frak p}_{K}\subset {\cal O}_{K}\subset K$ \`{a}
$\varpi_{K}k[[\varpi_{K}]]\subset k[[\varpi_{K}]]\subset
k((\varpi_{K}))$.
\vskip 2mm

On fixe un corps local $F$ et une famille finie $(E_{i})_{i\in I}$ non
vide d'extensions finies, s\'{e}parables et totalement ramifi\'{e}es
de $F$.  On note $n_{i}$ de degr\'{e} de $E_{i}$ sur $F$.  Pour chaque
$i\in I$, on se donne un \'{e}l\'{e}ment $\gamma_{i}$ de ${\frak
p}_{E_{i}}\subset {\cal O}_{E_{i}}\subset E_{i}$ qui engendre $E_{i}$
sur $F$, de sorte que $E_{i}\cong F[T]/(P_{i}(T))$ o\`{u} $P_{i}(T)\in
{\cal O}_{F}[T]$ est le polyn\^{o}me minimal de $\gamma_{i}$ sur $F$.
\vskip 1mm

{\it On suppose que les polyn\^{o}mes $P_{i}(T)$ unitaires et
irr\'{e}ductibles dans $F[T]$ sont deux \`{a} deux distincts}.
\vskip 2mm

On note $A_{i}$ la $k$-alg\`{e}bre int\`{e}gre
$$
A_{i}={\cal O}_{F}[\gamma_{i}]\subset {\cal O}_{E_{i}}.
$$
Elle est locale d'id\'{e}al maximal ${\frak m}_{i}={\frak
p}_{E_{i}}\cap A_{i}$, son corps des fractions est $E_{i}$, et sa
normalisation $\widetilde{A}_{i}\subset E_{i}$ n'est autre que ${\cal
O}_{E_{i}}$.
\vskip 2mm

On note $E_{I}=\prod_{i\in I}E_{i}$, $n_{I}=\sum_{i\in I}n_{i}$ la
dimension de cet espace vectoriel sur $F$, $\gamma_{I}=
(\gamma_{i})_{i\in I}\in\prod_{i\in I}A_{i}$ et
$$
A_{I}={\cal O}_{F}[\gamma_{I}]\subset\prod_{i\in I}A_{i}.
$$
La $k$-alg\`{e}bre $A_{I}$ est locale d'id\'{e}al maximal ${\frak
m}_{I}=\left(\prod_{i\in I}{\frak m}_{i}\right)\cap A_{I}$, son anneau
total des fractions est $E_{I}$ et sa normalisation est \'{e}gale
\`{a} $\widetilde{A}_{I}:=\prod_{i\in I}\widetilde{A}_{i}=\prod_{i\in
I}{\cal O}_{E_{i}}=:{\cal O}_{E_{I}}\subset E_{I}$.

Comme
$$
A_{I}\cong {\cal O}_{F}[T]/(P_{I}(T))=k[[\varpi_{F}]][T]/(P_{I}(T))
\cong k[[\varpi_{F},T]]/(P_{I}(T))
$$
o\`{u} $P_{I}(T)=\prod_{i\in I}P_{i}(T)$, la $k$-alg\`{e}bre locale
(int\`{e}gre et de dimension $1$) $A_{I}$ est de Gorenstein et son
module dualisant $\omega_{A_{I}}$ est \'{e}gal \`{a}
$$
\omega_{A_{I}}=\{y\in \Omega_{E_{I}/k}^{1}\mid \mathop{\rm
Res}\nolimits_{I}(xy)=0,~\forall x\in A_{I}\}\supset
\Omega_{\widetilde{A}_{I}/k}^{1},
$$
o\`{u} l'application $k$-lin\'{e}aire $\mathop{\rm Res}\nolimits_{I}:
\Omega_{E_{I}/k}^{1} \rightarrow k$ est la somme des applications
r\'{e}sidus $\mathop{\rm Res}\nolimits_{i}:\Omega_{E_{i}/k}^{1}
\rightarrow k$.

\thm PROPOSITION 1.1.1 (Rosenlicht, cf. [A-K 1] VIII, Proposition
(1.16))
\enonce
L'accou\-ple\-ment $(\widetilde{A}_{I}/A_{I})\times (\omega_{A_{I}}/
\Omega_{\widetilde{A}_{I}/k}^{1})\rightarrow k$ qui envoie
$(x+A_{I},y+\Omega_{\widetilde{A}_{I}/k}^{1})$ sur $\mathop{\rm
Res}\nolimits_{I}(xy)$ est un accouplement parfait.
\hfill\hfill$\square$
\endthm

On pose
$$
\delta_{I}=\mathop{\rm dim}\nolimits_{k}(\widetilde{A}_{I}/A_{I})
$$
et on note ${\frak a}_{I}$ le conducteur de $\widetilde{A}_{I}$ dans
$A_{I}$, c'est-\`{a}-dire l'id\'{e}al de $\widetilde{A}_{I}$ form\'{e}
des $x\in\widetilde{A}_{I}$ tels que $x\widetilde{A}_{I}\subset
A_{I}$. Cet id\'{e}al est contenu dans $A_{I}$ et il r\'{e}sulte de
la proposition ci-dessus que
$$
\mathop{\rm dim}\nolimits_{k}(A_{I}/{\frak a}_{I})=\delta_{I}.
$$
Pour chaque $i\in I$, on pose
$$
\delta_{i}=\mathop{\rm dim}\nolimits_{k}(\widetilde{A}_{i}/A_{i})
$$
et, pour chaque $i\not=j$ dans $I$, on note $r_{ij}=r_{ji}$ la
valuation du r\'{e}sultant dans ${\cal O}_{F}$ des polyn\^{o}mes
$P_{i}(T)$ et $P_{j}(T)$, c'est-\`{a}-dire la valuation de
$P_{j}(\gamma_{i})$ dans $E_{i}$ ou, ce qui revien au m\^{e}me, de
$P_{i}(\gamma_{j})$ dans $E_{j}$.  On v\'{e}rifie que
$$
\delta_{I}=\sum_{i\in I}\delta_{i}+{1\over 2}\sum_{{\scriptstyle
i,j\in I\atop\scriptstyle i\not=j}}r_{ij}
$$
et que
$$
{\frak a}_{I}=\prod_{i\in I}{\frak p}_{E_{i}}^{2\delta_{i}+\sum_{j\in
I\setminus\{i\}}r_{ij}}\subset\prod_{i\in I}{\cal O}_{E_{i}}=
{\cal O}_{E_{I}}.
$$

\section{1.2}{Fibres de Springer}

Rappelons qu'un ${\cal O}_{F}$-r\'{e}seau dans un $F$-espace
vectoriel de dimension finie $E$ est un sous-${\cal O}_{F}$-Module
de $E$ de rang \'{e}gal \`{a} la dimension de $E$. Si $M$ et $N$
sont deux tels ${\cal O}_{F}$-r\'{e}seaux, l'indice de $M$
relativement \`{a} $N$ est l'entier relatif
$$
[M:N]=\mathop{\rm dim}\nolimits_{k}(M/P)-\mathop{\rm
dim}\nolimits_{k}(N/P)
$$
o\`{u} $P$ est n'importe quel ${\cal O}_{F}$-r\'{e}seau de $E$ contenu
\`{a} la fois dans $M$ et dans $N$.  Par exemple, $A_{I}$ et
$\widetilde{A}_{I}$ sont des ${\cal O}_{F}$-r\'{e}seaux dans $E_{I}$
et $[\widetilde{A}_{I}:A_{I}]=\delta_{I}$.
\vskip 2mm

Soient $N\geq 0$ et $d$ des entiers et soit
$$
M\subset\varpi_{F}^{-N}A_{I}\subset E_{I}
$$
un ${\cal O}_{F}$-r\'{e}seau qui est d'indice $d$ relativement au
${\cal O}_{F}$-r\'{e}seau particulier $A_{I}$ (bien entendu, pour
qu'il existe un tel $M$, il faut que $d\leq n_{I}N$). La
multiplication par $\varpi_{F}$ induit un endomorphisme nilpotent du
$k$-espace vectoriel $\varpi_{F}^{-N}A_{I}/M$ de dimension $n_{I}N-d$.
Par suite, $M$ contient automatiquement $\varpi_{F}^{(n_{I}-1)N-d}
A_{I}$ et la donn\'{e}e de $M$ \'{e}quivaut \`{a} la donn\'{e}e du
sous-espace vectoriel
$$
M/\varpi_{F}^{(n_{I}-1)N-d}A_{I}\subset\varpi_{F}^{-N}A_{I}/
\varpi_{F}^{(n_{I}-1)N-d}A_{I}
$$
stable par l'endomorphisme nilpotent induit par la multiplication par
$\varpi_{F}$.

Pour toute extension $k'$ de $k$, on note
$F'=k'\widehat{\otimes}_{k}F=k'((\varpi_{F}))$,
$E_{i}'=k'\widehat{\otimes}_{k}E_{i}=k'((\varpi_{E_{i}}))$, etc ...
les compl\'{e}t\'{e}s $\varpi_{F}$-adiques de $k'\otimes_{k}F$,
$k'\otimes_{k}E_{i}$, etc ...  Ce que l'on vient de dire vaut encore
apr\`{e}s que l'on ait remplac\'{e} $k$ par $k'$, $F$ par $F'$, etc
...

Pour $k'$ variable, les ${\cal O}_{F'}$-r\'{e}seaux dans $E_{I}'$ qui
sont d'indice $d$ relativement au ${\cal O}_{F'}$-r\'{e}seau
particulier $A_{I}'$ et qui sont contenus dans $\varpi_{F}^{-N}A_{I}'$
sont donc naturellement les $k'$-points d'un $k$-sch\'{e}ma projectif
r\'{e}duit $R_{I,N}^{d}$, \`{a} savoir le ferm\'{e} (r\'{e}duit) de la
grassmannienne des $(n_{I}-1)(n_{I}N-d)$-plans dans le $k$-espace
vectoriel $\varpi_{F}^{-N}A_{I}/\varpi_{F}^{(n_{I}-1)N-d}A_{I}$ de
dimension $n_{I}(n_{I}N-d)$, form\'{e} des plans qui sont stables par
l'endomorphisme nilpotent induit par la multiplication par
$\varpi_{F}$.

Pour $d$ fix\'{e}, les $k$-sch\'{e}mas projectifs $R_{I,N}^{d}$
s'organisent en un syst\`{e}me inductif d'immersions ferm\'{e}es
$$
\cdots\hookrightarrow R_{I,N}^{d}\hookrightarrow R_{I,N+1}^{d}
\hookrightarrow\cdots
$$
et on note $R_{I}^{d}$ le ind-$k$-sch\'{e}ma {\og}{limite}{\fg}. Le
ind-$k$-sch\'{e}ma des ${\cal O}_{F}$-r\'{e}seaux de $E_{I}$ est
par d\'{e}finition la somme disjointe
$$
R_{I}=\coprod_{d\in {\Bbb Z}}R_{I}^{d}.
$$
\vskip 2mm

Toujours pour $k'$ variable, les ${\cal O}_{F'}$-r\'{e}seaux
$(M\subset E_{I}')\in R_{I,N}^{d}(k')$ tels que
$$
\gamma_{I}M\subset M
$$
sont les $k'$-points d'un sous-$k$-sch\'{e}ma ferm\'{e} r\'{e}duit
$X_{I,N}^{d}$ de $R_{I,N}^{d}$. Pour $d$ fix\'{e}, les
$X_{I,N}^{d}\subset R_{I,N}^{d}$ s'organisent en un syst\`{e}me
inductif et on note $X_{I}^{d}\subset R_{I}^{d}$ le
sous-ind-$k$-sch\'{e}ma ferm\'{e} r\'{e}duit {\og}{limite}{\fg}.

\thm D\'{E}FINITION 1.2.1
\enonce
La {\rm fibre de Springer} en $\gamma_{I}$ est le
sous-ind-$k$-sch\'{e}ma ferm\'{e} r\'{e}duit
$$
X_{I}=\coprod_{d\in {\Bbb Z}}X_{I}^{d}\subset\coprod_{d\in {\Bbb
Z}}R_{I}^{d}=R_{I}
$$
des ${\cal O}_{F}$-r\'{e}seaux de $E_{I}$ stabilis\'{e} par
$\gamma_{I}$.
\endthm

Bien entendu, chaque $X_{I}^{d}=X_{I}\cap R_{I}^{d}$ est une partie ouverte
et ferm\'{e}e de $X_{I}$.
\vskip 2mm

La remarque \'{e}vidente suivante est essentielle pour la suite:
\vskip 2mm

{\it On peut identifier les $k'$-points de $X_{I}$ aux
sous-$A_{I}'$-modules $M\subset E_{I}'$ qui sont de rang $1$ en chaque
point g\'{e}n\'{e}rique de $\mathop{\rm Spec}(A_{I}')$}.
\vskip 2mm

Pour $k'$ variable, le groupe $E_{I}'^{\times}/A_{I}'^{\times}$ est de
mani\`{e}re naturelle le groupe des $k'$-points d'un $k$-sch\'{e}ma en
groupes commutatifs $G_{I}$ qui est lisse et de dimension $\delta_{I}$
sur $k$.  Plus pr\'{e}cis\'{e}ment, le groupe des composantes connexes
de $G_{I}$ est le quotient $E_{I}^{\times}/{\cal O}_{E_{I}}^{\times}
=\prod_{i\in I}(E_{i}^{\times}/{\cal O}_{E_{i}}^{\times})$ qui est
canoniquement isomorphe \`{a} $\Lambda_{I}:={\Bbb Z}^{I}$; la
composante neutre $G_{I}^{0}$ de $G_{I}$, qui admet ${\cal
O}_{E_{I}'}^{\times}/ A_{I}'^{\times}$ pour groupe des $k'$-points, est
une extension d'un tore $T_{I}$ (le quotient de ${\Bbb G}_{{\rm
m},k}^{I}$ par ${\Bbb G}_{{\rm m},k}$ plong\'{e} diagonalement) par un
sch\'{e}ma en groupes unipotents $U_{I}$ de type fini dont le groupe
des $k'$-points est
$$
U_{I}(k')=\left(\textstyle\prod_{i\in I}(1+{\frak p}_{E_{I}'})\right)/
(1+{\frak m}_{I}')
$$
o\`{u} ${\frak m}_{I}'$ est l'id\'{e}al maximal de $A_{I}'$.

L'action par homoth\'{e}ties de $E_{I}'^{\times}/A_{I}'^{\times}$
sur les r\'{e}seaux $M\in X_{I}(k')$ provient d'une action
alg\'{e}brique naturelle de $G_{I}$ sur $X_{I}$. Cette action permute
les composantes $X_{I}^{d}$ de $X_{I}$ suivant la r\`{e}gle
$$
g\cdot X_{I}^{d}=X_{I}^{d+|\lambda |}
$$
o\`{u} $\lambda$ est l'image de $g\in G_{I}$ dans $\Lambda_{I}$ et
o\`{u} on a pos\'{e} $|\lambda |=\sum_{\iota\in I}\lambda_{i}$ pour
chaque $\lambda\in\Lambda_{I}$.

Le $k$-sch\'{e}ma $X_{I}$ contient le $k$-point particulier $M=A_{I}$.
Le fixateur dans $G_{I}$ de ce point particulier est r\'{e}duit \`{a}
l'\'{e}l\'{e}ment neutre et son orbite $X_{I}^{\circ}=G_{I}\cdot
A_{I}$ est donc une partie de $X_{I}$ isomorphe \`{a} $G_{I}$.

\thm LEMME 1.2.2
\enonce
La $G_{I}$-orbite $X_{I}^{\circ}$ est l'ouvert de $X_{I}$ dont les
$k'$-points sont les $M'$ qui sont libres de rang $1$ en tant que
$A_{I}'$-modules.
\hfill\hfill$\square$
\endthm

Tout scindage $\sigma :\Lambda_{I}^{0}\hookrightarrow G_{I}$ de
l'extension
$$
1\rightarrow G_{I}^{0}\rightarrow G_{I}\rightarrow
\Lambda_{I}\rightarrow 0
$$
au-dessus du sous-groupe
$$
\Lambda_{I}^{0}:=\{\lambda\in\Lambda_{I}\mid |\lambda|=0\}
\subset\Lambda_{I}
$$
d\'{e}finit une action libre de $\Lambda_{I}^{0}$ sur $X_{I}$ qui
pr\'{e}serve les composantes $X_{I}^{d}$. On notera $Z_{I}=X_{I}/
\sigma (\Lambda_{I}^{0})$ le $k$-espace quotient correspondant et
$Z_{I}^{d}=X_{I}^{d}/\sigma (\Lambda_{I}^{0})$ ses composantes.

\thm TH\'{E}OR\`{E}ME 1.2.3 (Kazhdan-Lusztig, cf. [K-L] \S 3)
\enonce
La fibre de Springer $X_{I}$ est en fait un $k$-sch\'{e}ma localement
de type fini et de dimension finie, dont les composantes connexes sont
exactement les $X_{I}^{d}$, $d\in {\Bbb Z}$.

Pour tout scindage $\sigma :\Lambda_{I}^{0}\hookrightarrow G_{I}$
comme ci-dessus, les $k$-espaces quotients $Z_{I}^{d}$ correspondants
sont des $k$-sch\'{e}mas projectifs, $Z_{I}$ est le $k$-sch\'{e}ma
somme disjointe des $Z_{I}^{d}$ et l'application quotient
$X_{I}\rightarrow Z_{I}$ est un rev\^{e}tement \'{e}tale galoisien de
groupe de Galois $\Lambda_{I}^{0}$.
\hfill\hfill$\square$
\endthm

\section{1.3}{Lemme fondamental arithm\'{e}tique (d'apr\`{e}s
Kottwitz)}

Soit $q$ une puissance d'un nombre premier $p$ et ${\Bbb F}_{q}\subset
{\Bbb F}_{q^{2}}$ les corps finis \`{a} $q$ et $q^{2}$
\'{e}l\'{e}ments respectivement, qui sont contenus dans une
cl\^{o}ture alg\'{e}brique fix\'{e}e $\overline{{\Bbb F}}_{p}$ du
corps premier ${\Bbb F}_{p}$.  On note $x\mapsto x^{\ast}$ les
\'{e}l\'{e}ments de Frobenius g\'{e}om\'{e}triques des groupes de
Galois $\mathop{\rm Gal}(\overline{{\Bbb F}}_{p}/{\Bbb F}_{q})$ et
$\mathop{\rm Gal}({\Bbb F}_{q^{2}}/{\Bbb F}_{q})$.  On fixe aussi un
\'{e}l\'{e}ment $\varepsilon\in {\Bbb F}_{q^{2}}^{\times}$ tel que
$\varepsilon^{\ast}=-\varepsilon$.  On prend pour $k$ le corps ${\Bbb
F}_{q^{2}}$.

On se donne un corps local $F_{0}$ de corps r\'{e}siduel ${\Bbb
F}_{q}$ et une famille $(E_{i,0})_{i\in I}$ d'extensions finies,
s\'{e}parables et totalement ramifi\'{e}es de $F_{0}$, et on prend
pour $F$ le corps local $F={\Bbb F}_{q^{2}}\otimes_{{\Bbb
F}_{q}}F_{0}$ et pour $E_{i}$, $i\in I$, les corps locaux $E_{i}={\Bbb
F}_{q^{2}}\otimes_{{\Bbb F}_{q}}E_{i,0}$.

Pour toute extension finie $k'$ de ${\Bbb F}_{q}$ contenue dans
$\overline{{\Bbb F}}_{p}$, on note $F'=k'\otimes_{{\Bbb F}_{q}}F_{0}$,
$E_{i}'=k'\otimes_{{\Bbb F}_{q}}E_{i,0}$, etc ...  et on note encore
$x\mapsto x^{\ast}$ les rel\`{e}vements naturels de l'\'{e}l\'{e}ment
de Frobenius g\'{e}om\'{e}trique de $\mathop{\rm Gal}(k'/{\Bbb F}_{q})$ aux
groupes de Galois $\mathop{\rm Gal}(F'/F_{0})$, $\mathop{\rm
Gal}(E_{i}'/E_{i,0})$, etc ...  Bien s\^{u}r, si $k'$ contient $k$, on
a $F'=k'\otimes_{k}F$, $E_{i}'=k'\otimes_{k}E_{i}$, etc ...

On se donne des $\gamma_{i}\in {\frak p}_{E_{i}}$ comme dans la
section 1.1 et on fait l'hypoth\`{e}se suppl\'{e}mentaire suivante:
{\it pour chaque $i\in I$, on a}
$$
\mathop{\rm Tr}\nolimits_{E_{i}/E_{i,0}}(\gamma_{i})
(=\gamma_{i}^{\ast}+\gamma_{i})=0.\leqno{(1.3.1)}
$$

Il r\'{e}sulte de cette hypoth\`{e}se que l'\'{e}l\'{e}ment
$$
\alpha_{I}=\varepsilon^{n_{I}-1}{{\rm d}P_{I}\over {\rm
d}T}(\gamma_{I})\in E_{I}^{\times}
$$
est en fait dans $E_{I,0}^{\times}\subset E_{I}^{\times}$.  On munit le
$F$-espace vectoriel $E_{I}$ (de dimension finie $n_{I}$) de la forme
hermitienne
$$
\langle x,y\rangle_{I}=\mathop{\rm Tr}
\nolimits_{E_{I}/F}(\alpha_{I}^{-1}x^{\ast}y).
$$
On a choisi $\alpha_{I}$ pour que $A_{I}\subset E_{I}$ soit son propre
orthogonal pour cette forme hermitienne.

La forme sesquilin\'{e}aire ci-dessus induit sur la fibre de Springer
$X_{I}$ une ${\Bbb F}_{q}$-structure rationnelle ${\cal X}_{I}$ dont
l'endomorphisme de Frobenius $\mathop{\rm Frob}\nolimits_{{\cal
X}_{I}/{\Bbb F}_{q}}$ envoie, quel que soit l'extension finie $k'$ de
$k$ contenue dans $\overline{{\Bbb F}}_{p}$, un $k'$-point $M\subset
E_{I}'$ de $X_{I}$ sur le $k'$-point
$$
M^{\perp_{I}}=\{y\in E_{I}'\mid \langle M,y\rangle_{I}' \subset
O_{F'}\}\subset E_{I}'
$$
de $X_{I}$, o\`{u} on a not\'{e}
$$
\langle x,y\rangle_{I}'=\mathop{\rm Tr}
\nolimits_{E_{I}'/F'}(\alpha_{I}^{-1}x^{\ast}y).
$$
On a $\gamma_{I}M^{\perp_{I}}\subset M^{\perp_{I}}$ puisque
$\gamma_{I}^{\ast}=-\gamma_{I}$ par hypoth\`{e}se, et le fl\`{e}che
canonique
$$
(M^{\ast})^{\ast}\hookrightarrow (M^{\perp_{I}})^{\perp_{I}}\subset E_{I}'
$$
est bijective, de sorte que
$$
\mathop{\rm Frob}\nolimits_{{\cal X}_{I}/{\Bbb F}_{q}}\circ
\mathop{\rm Frob}\nolimits_{{\cal X}_{I}/{\Bbb F}_{q}}=\mathop{\rm
Frob}\nolimits_{X_{I}/{\Bbb F}_{q^{2}}}
$$
est l'endomorphisme de Frobenius de $X_{I}$ sur ${\Bbb F}_{q^{2}}$.

Le ${\Bbb F}_{q^{2}}$-point $A_{I}\subset E_{I}$ de $X_{I}$ est fix\'{e} par
l'endomorphisme de Frobenius $\mathop{\rm Frob}\nolimits_{{\cal
X}_{I}/{\Bbb F}_{q}}$ et, quel que soit le point $M\subset E_{I}$ de
$X_{I}$, on a
$$
[\mathop{\rm Frob}\nolimits_{{\cal X}_{I}/{\Bbb
F}_{q}}(M):A_{I}]=-[M:A_{I}].
$$
Par suite, $\mathop{\rm Frob}\nolimits_{{\cal X}_{I}/{\Bbb F}_{q}}$
induit une ${\Bbb F}_{q}$-structure rationnelle ${\cal X}_{I}^{0}$ sur
le ${\Bbb F}_{q^{2}}$-sch\'{e}ma $X_{I}^{0}$ et $A_{I}$ est un ${\Bbb
F}_{q}$-point de ${\cal X}_{I}^{0}$.

On a aussi une ${\Bbb F}_{q}$-structure rationnelle ${\cal G}_{I}$ sur
le ${\Bbb F}_{q^{2}}$-sch\'{e}ma en groupes $G_{I}$ dont
l'endomorphisme de Frobenius $\mathop{\rm Frob}\nolimits_{{\cal G}_{I}
/{\Bbb F}_{q}}$ envoie $g\in G_{I}(k')=E_{I}'^{\times}/
A_{I}'^{\times}$ sur $(g^{\ast})^{-1}$.  Cet endomorphisme induit sur
le quotient $\Lambda_{I}=G_{I}/G_{I}^{0}$ l'automorphisme
$\lambda\mapsto -\lambda$.

Les ${\Bbb F}_{q}$-structures rationnelles ${\cal X}_{I}$ et ${\cal
G}_{I}$ sont compatibles au sens o\`{u}
$$
\mathop{\rm Frob}\nolimits_{{\cal X}_{I}/{\Bbb F}_{q}}(g\cdot M)=
\mathop{\rm Frob}\nolimits_{{\cal G}_{I}/{\Bbb F}_{q}}(g)\cdot
\mathop{\rm Frob}\nolimits_{{\cal X}_{I}/{\Bbb F}_{q}}(M)
$$
pour tout $g\in G_{I}$ et tout $M\in X_{I}$.  Fixons un scindage
$\sigma :\Lambda_{I}^{0}\hookrightarrow G_{I}$ de telle sorte que
$$
\mathop{\rm Frob}\nolimits_{{\cal G}_{I}/{\Bbb F}_{q}}(\sigma
(\lambda ))=\sigma (-\lambda ),~\forall \lambda\in\Lambda_{I}^{0},
$$
ce qui est toujours possible (il suffit de choisir pour chaque $i\in
I$, une uniformisante $\varpi_{E_{i}}$ de $E_{i}$ qui est d\'{e}j\`{a}
dans $E_{i,0}$).  Alors la ${\Bbb F}_{q}$-structure rationnelle ${\cal
X}_{I}$ en induit une ${\cal Z}_{I}$ sur le quotient
$Z_{I}=X_{I}/\sigma (\Lambda_{I}^{0})$.  Comme on a
$$
\mathop{\rm Frob}\nolimits_{{\cal Z}_{I}/{\Bbb
F}_{q}}(Z_{I}^{i})=Z_{I}^{-i},~\forall i \in {\Bbb Z},
$$
l'endomorphisme de Frobenius $\mathop{\rm Frob}\nolimits_{{\cal
Z}_{I}/{\Bbb F}_{q}}$ induit une ${\Bbb F}_{q}$-structure rationnelle
${\cal Z}_{I}^{0}$ sur $Z_{I}^{0}=X_{I}/\sigma (\Lambda_{I}^{0})$ et
les ${\Bbb F}_{q}$-points de ${\cal Z}_{I}$, c'est-\`{a}-dire les
points fixes de l'endomorphisme de Frobenius $\mathop{\rm
Frob}\nolimits_{{\cal Z}_{I}/{\Bbb F}_{q}}$, sont
n\'{e}cessairement contenus dans ${\cal Z}_{I}^{0}$.

\thm LEMME 1.3.2
\enonce
Les ${\Bbb F}_{q}$-points de ${\cal Z}_{I}$, ou ce qui revient au
m\^{e}me de ${\cal Z}_{I}^{0}$, sont les $\Lambda_{I}^{0}$-orbites $z$
des ${\Bbb F}_{q^{2}}$-points $x$ de $X_{I}^{0}$ pour lesquels il
existe $\lambda_{x}\in\Lambda_{I}^{0}$ tel que
$$
\mathop{\rm Frob}\nolimits_{{\cal X}_{I}/{\Bbb F}_{q}}(x)=\sigma
(\lambda_{x})\cdot x.
$$

De plus, pour tout $z\in {\cal Z}_{I}^{0}({\Bbb F}_{q})$, la classe
$\overline{\lambda}_{x}$ dans $\Lambda_{I}^{0}/2\Lambda_{I}^{0}$ de
l'\'{e}l\'{e}ment $\lambda_{x}$ ci-dessus est uniquement
d\'{e}termin\'{e}e par $z$ et on a donc une application
$$
{\cal Z}_{I}^{0}({\Bbb F}_{q})
\rightarrow\Lambda_{I}^{0}/2\Lambda_{I}^{0}.
$$
\hfill\hfill$\square$
\endthm

Pour chaque $\overline{\lambda}\in\Lambda_{I}^{0}/2\Lambda_{I}^{0}$,
on pose
$$
\mathop{\rm O}\nolimits_{\overline{\lambda}}^{I}=|\{z\in {\cal
Z}_{I}^{0} ({\Bbb F}_{q})\mid\overline{\lambda}_{z}=
\overline{\lambda}\}|.
$$
Pour chaque caract\`{e}re d'ordre $2$, $\kappa :\Lambda_{I}^{0}
\rightarrow \{\pm 1\}$, ou ce qui revient au m\^{e}me, pour chaque
caract\`{e}re $\kappa :\Lambda_{I}^{0}/2\Lambda_{I}^{0} \rightarrow
\{\pm 1\}$, on pose
$$
\mathop{\rm O}\nolimits^{I,\kappa}=\sum_{\overline{\lambda}\in
\Lambda_{I}^{0}/2\Lambda_{I}^{0}}\kappa (\overline{\lambda})
\mathop{\rm O}\nolimits_{\overline{\lambda}}^{I}.
$$
En particulier, pour $\kappa =1$ le caract\`{e}re trivial, on pose
$$
\mathop{\rm SO}\nolimits^{I}=\mathop{\rm O}\nolimits^{I,1}=
\sum_{\overline{\lambda}\in\Lambda_{I}^{0}/2\Lambda_{I}^{0}}
\mathop{\rm O}\nolimits_{\overline{\lambda}}^{I} =|{\cal Z}_{I}^{0}
({\Bbb F}_{q})|.
$$

\rem Remarque {\rm 1.3.3}
\endrem
Soit $U_{I}$ le groupe unitaire sur $F_{0}$ d\'{e}fini par le
$F$-espace vectoriel hermitien $(E_{I}, \langle~,~\rangle_{I})$ et
${\frak u}_{I}$ son alg\`{e}bre de Lie.  Le groupe r\'{e}ductif
$U_{I}$ est naturellement muni d'un tore maximal elliptique $T_{I}$
tel que $T_{I}(F_{0})$ soit form\'{e} des \'{e}l\'{e}ments $t$ de
$E_{I}^{\times}$ tels que $t^{\ast}t=1$ et $\gamma_{I}$ est un
\'{e}l\'{e}ment r\'{e}gulier de l'alg\`{e}bre de Lie ${\frak t}_{I}$
de $T_{I}$.

Les classes de $U_{I}(F_{0})$-conjugaison dans la classe de
conjugaison stable de $\gamma_{I}$ peuvent \^{e}tre
param\'{e}tr\'{e}es par $\Lambda_{I}^{0}/2\Lambda_{I}^{0}$ et chaque
$\mathop{\rm O}\nolimits_{\overline{\lambda}}^{I}$ est l'{\it
int\'{e}grale orbitale} sur la classe de conjugaison correspondante
\`{a} $\overline{\lambda}$ de la fonction caract\'{e}ristique du
r\'{e}seau dans ${\frak u}_{I}(F_{0})$ qui fixe $A_{I}\subset E_{I}$.

Par suite, $\mathop{\rm O}\nolimits^{I,\kappa}$ est un
$\kappa$-int\'{e}grale orbitale et $\mathop{\rm SO}\nolimits^{I}$
une int\'{e}grale orbitale stable, au sens de Labesse et Langlands
\hfill\hfill$\square$
\vskip 3mm

Bien entendu, dans tout ce qui pr\'{e}c\`{e}de on peut remplacer $I$
par une partie $J$ de $I$, $E_{I}$ par $E_{J}=\prod_{i\in J}E_{i}$,
$\gamma_{J}=(\gamma_{i})_{i\in J}$, etc ...  {\it On fera cependant
attention au fait que la forme hermitienne $\langle~,~\rangle_{J}$
n'est pas la forme hermitienne sur $E_{J}\subset E_{I}$ qui est
induite par la forme hermitienne $\langle~,~\rangle_{I}$ {\rm (}voir
$1.4.8${\rm )}.}

La donn\'{e}e d'un caract\`{e}re $\kappa :\Lambda_{I}^{0}/
2\Lambda_{I}^{0}\rightarrow \{\pm 1\}$ \'{e}quivaut \`{a} la
donn\'{e}e d'une partition $I=I_{1}\amalg I_{2}$ de $I$, le
caract\`{e}re $\kappa_{I_{1},I_{2}}$ associ\'{e} \`{a} une telle
partition \'{e}tant donn\'{e} par la formule
$$
\kappa (\lambda )=(-1)^{\sum_{i\in I_{1}}\lambda_{i}}
=(-1)^{\sum_{i\in I_{2}}\lambda_{i}},~\forall
\lambda\in\Lambda_{I}^{0}.
$$

Langlands et Shelstad (cf.  [La-Sh]) ont conjectur\'{e} toute une
famille de relations entre $\kappa$-int\'{e}grales orbitales sur un
groupe r\'{e}ductif $G$ sur un corps local non archim\'{e}dien et
int\'{e}grales orbitales stables sur le groupe endoscopique $H$ de $G$
correspondant \`{a} $\kappa$.  Cette famille de relations
conjecturales est commun\'{e}ment appel\'{e}e le {\og}{Lemme
fondamental}{\fg}.  En particulier, pour $G=U_{I}$ et
$\kappa=\kappa_{I_{1},I_{2}}$, le groupe endoscopique $H$ n'est autre
que $U_{I_{1}}\times U_{I_{2}}$ et, compte tenu des calculs des
facteurs de transferts pour les groupes classiques, effectu\'{e}s par
Waldspurger (cf.  [Wa], Chapitre X), la conjecture de Langlands et
Shelstad s'\'{e}nonce simplement:

\thm CONJECTURE 1.3.4 (Lemme fondamental arithm\'{e}tique)
\enonce
Pour toute partition $I=I_{1}\amalg I_{2}$ de $I$, on a la relation
$$
\mathop{\rm O}\nolimits^{I,\kappa_{I_{1},I_{2}}}=q^{r_{I_{1},I_{2}}}\times
\mathop{\rm SO}\nolimits^{I_{1}}\times \mathop{\rm SO}\nolimits^{I_{2}}.
$$
o\`{u} $ r_{I_{1},I_{2}}=[A_{I_{1}}\times A_{I_{2}}:A_{I}]$.
\endthm

\rem Remarque $1.3.5$
\endrem
On a $r_{I_{1},I_{2}}=\delta_{I}-(\delta_{I_{1}}+\delta_{I_{2}})=\sum_{i_{1}\in
I_{1},i_{2}\in I_{2}}r_{i_{1},i_{2}}$.

\hfill\hfill$\square$
\vskip 3mm

Fixons un nombre premier $\ell\not=p$.  D'apr\`{e}s Grothendieck [Gr
1], le nombre des ${\Bbb F}_{q}$-points d'un ${\Bbb F}_{q}$-sch\'{e}ma
projectif est la trace de l'action de l'endomorphisme de Frobenius sur
la cohomologie $\ell$-adique de ce sch\'{e}ma.  On a donc
$$
\mathop{\rm SO}\nolimits^{J}=\mathop{\rm tr}(\mathop{\rm Frob}
\nolimits_{{\cal Z}_{J}^{0}},R\Gamma (\overline{{\Bbb F}}_{p}
\otimes_{{\Bbb F}_{q}}{\cal Z}_{J}^{0},{\Bbb Q}_{\ell}))\leqno{(1.3.6)}
$$
pour toute partie $J$ de $I$.

Pour toute partition $I=I_{1}\amalg I_{2}$, le caract\`{e}re
$\kappa_{I_{1},I_{2}}$ peut \^{e}tre vu \`{a} valeurs dans $\{\pm 1\}
{\Bbb Q}_{\ell}$ et on peut {\og}{pousser}{\fg} le rev\^{e}tement
\'{e}tale galoisien $X_{I}^{0}\rightarrow X_{I}^{0}/\sigma
(\Lambda_{I}^{0})=Z_{I}^{0}$ de groupe de Galois $\Lambda_{I}^{0}$ par
$\kappa_{I_{1},I_{2}}$ pour obtenir un syst\`{e}me local $\ell$-adique
${\cal L}_{I_{1},I_{2}}$ de rang $1$ sur $Z_{I}^{0}$.  L'endomorphisme
de Frobenius $\mathop{\rm Frob}\nolimits_{{\cal Z}_{J}^{0}}$ se
rel\`{e}ve \`{a} ${\cal L}_{I_{1},I_{2}}$ vu notre choix de la section
$\sigma$ et, toujours d'apr\`{e}s Grothendieck, on d\'{e}duit du lemme
1.3.2 que
$$
\mathop{\rm O}\nolimits^{I,\kappa_{I_{1},I_{2}}}=\mathop{\rm tr}
(\mathop{\rm Frob}\nolimits_{{\cal Z}_{I}},R\Gamma (\overline{{\Bbb
F}}_{p}\otimes_{{\Bbb F}_{q}}{\cal Z}_{I}^{0},{\cal
L}_{I_{1},I_{2}})).\leqno{(1.3.7)}
$$

\section{1.4}{Lemme fondamental g\'{e}om\'{e}trique}

Revenons \`{a} la situation de la section 1.1 o\`{u} $k$ est un corps
parfait arbitraire de caract\'{e}ristique $\not=\ell$ et fixons une
cl\^{o}ture alg\'{e}brique $\overline{k}$ de $k$.

Le lemme fondamental arithm\'{e}tique nous am\`{e}ne tout
naturellement \`{a} esp\'{e}rer que, pour chaque partition
$I=I_{1}\amalg I_{2}$, il existe un isomorphisme canonique
$$
R\Gamma (\overline{k}\otimes_{k}(Z_{I_{1}}^{0}\times
Z_{I_{2}}^{0}),{\Bbb Q}_{\ell})[-2r_{I_{1},I_{2}}](-r_{I_{1},I_{2}})
\buildrel\sim\over\longrightarrow R\Gamma (\overline{k}
\otimes_{k}Z_{I}^{0},{\cal L}_{I_{1},I_{2}})
$$
o\`{u} le ${\Bbb Q}_{\ell}$-syst\`{e}me local ${\cal L}_{I_{1},I_{2}}$
de rang $1$ sur $Z_{I}^{0}$ est d\'{e}fini, comme dans la section
pr\'{e}c\'{e}dente, en poussant le rev\^{e}tement \'{e}tale galoisien
$X_{I}^{0}\rightarrow X_{I}^{0}/\sigma (\Lambda_{I}^{0})=Z_{I}^{0}$ de
groupe de Galois $\Lambda_{I}^{0}$ par $\kappa_{I_{1},I_{2}}$.

Bien s\^{u}r, dans la situation arithm\'{e}tique de la section 1.3, un
tel isomorphisme canonique, qui s'\'{e}crit encore
$$
R\Gamma (\overline{{\Bbb F}}_{p}\otimes_{{\Bbb F}_{q}}({\cal
Z}_{I_{1}}^{0}\times {\cal Z}_{I_{2}}^{0}),{\Bbb Q}_{\ell})
[-2r_{I_{1},I_{2}}](-r_{I_{1},I_{2}})
\buildrel\sim\over\longrightarrow R\Gamma (\overline{{\Bbb
F}}_{p}\otimes_{{\Bbb F}_{q}}{\cal Z}_{I},{\cal L}_{I_{1},I_{2}}),
$$
serait automatiquement compatible aux actions des endomorphismes de
Frobenius relatifs \`{a} ${\Bbb F}_{q}$ et son existence impliquerait
la conjecture 1.3.4.

Nous allons voir qu'un tel isomorphisme existe
{\og}{virtuellement}{\fg}.  Plus pr\'{e}cis\'{e}ment, nous allons
d\'{e}montrer dans cette section:

\thm TH\'{E}OR\`{E}ME 1.4.1
\enonce
Dans le groupe de Grothendieck des ${\Bbb Q}_{\ell}[\mathop{\rm
Gal}(\overline{k}/k)]$-modules, on a
$$
\sum_{i}(-1)^{i}[H^{i}(\overline{k}\otimes_{k}Z_{I}^{0},{\cal
L}_{I_{1},I_{2}})]=\sum_{i}(-1)^{i}[H^{i-2r_{I_{1},I_{2}}}(\overline{k}
\otimes_{k}(Z_{I_{1}}^{0}\times Z_{I_{2}}^{0}),{\Bbb Q}_{\ell})
(-r_{I_{1},I_{2}})].
$$
\endthm

Pour all\'{e}ger la r\'{e}daction, nous noterons simplement
$$
r=r_{I_{1},I_{2}}=[A_{I_{1}}\times A_{I_{2}}:A_{I}]
$$
et nous ne consid\'{e}rerons dans ce qui suit que des points
rationnels sur $k$.  Bien s\^{u}r, les m\^{e}mes arguments valent pour
des points rationnels sur une extension finie arbitraire de $k$
apr\`{e}s avoir remplacer $F$ par $F'=k'\otimes_{k}F$, $E_{i}$ par
$E_{i}'=k'\otimes_{k}E_{i}$, etc ...

Pour chaque $(M\subset E_{I})\in X_{I}^{0}$, soient $M_{1}'$ et
$M_{2}'$ ses intersections avec $E_{I_{1}}=E_{I_{1}}\times \{0\}$ et
$E_{I_{1}}=\{0\}\times E_{I_{2}}$ dans $E_{I_{1}}\times
E_{I_{2}}=E_{I}$, et $M_{1}''$ et $M_{2}''$ ses projections sur les
facteurs $E_{I_{1}}$ et $E_{I_{2}}$ de $E_{I_{1}}\times
E_{I_{2}}=E_{I}$. On a des suites exactes
$$
0\rightarrow M_{1}'\rightarrow M\rightarrow M_{2}''\rightarrow 0
$$
et
$$
0\rightarrow M_{2}'\rightarrow M\rightarrow M_{1}''\rightarrow 0,
$$
et les inclusions
$$
M_{1}'\subset M_{1}''\subset E_{I_{1}}\hbox{ et }
M_{2}'\subset M_{2}''\subset E_{I_{2}}.
$$
Par suite, si on pose
$$
\mathop{\rm ind}\nolimits_{\alpha}^{\prime}(M)=
[M_{\alpha}':A_{I_{\alpha}}]\hbox{ et }\mathop{\rm
ind}\nolimits_{\alpha}^{\prime\prime}(M)=
[M_{\alpha}'':A_{I_{\alpha}}]
$$
pour $\alpha =1,2$, on a les relations
$$
\mathop{\rm ind}\nolimits_{1}^{\prime}(M)+
\mathop{\rm ind}\nolimits_{2}^{\prime\prime}(M)=
\mathop{\rm ind}\nolimits_{2}^{\prime}(M)+
\mathop{\rm ind}\nolimits_{1}^{\prime\prime}(M)=[M:A_{I_{1}}\times
A_{I_{2}}]=-r.
$$

\thm LEMME 1.4.2
\enonce
Pour chaque $(M\subset E_{I})\in X_{I}^{0}$, on a
$$
0\leq\mathop{\rm ind}\nolimits_{1}^{\prime\prime}(M)-
\mathop{\rm ind}\nolimits_{1}^{\prime}(M)=
\mathop{\rm ind}\nolimits_{2}^{\prime\prime}(M)-
\mathop{\rm ind}\nolimits_{2}^{\prime}(M)\leq r
$$
\endthm

\rem Preuve
\endrem
Pour $\{\alpha ,\beta\}=\{1,2\}$, on a
$$
P_{I_{\beta}}(\gamma_{I_{\alpha}})M_{\alpha}''\subset M_{\alpha}'\subset
M_{\alpha}''
$$
puisque $P_{I_{\beta}}(\gamma_{I})M\subset M$ et que
$P_{I_{\beta}}(\gamma_{I_{\beta}})=0$.  Par suite, on a
$$
[M_{\alpha}'':M_{\alpha}']\leq [M_{\alpha}'':P_{I_{\beta}}(\gamma_{I_{\alpha}})
M_{\alpha}'']
$$
o\`{u}
$$\displaylines{
\qquad [M_{\alpha}'':P_{I_{\beta}}(\gamma_{I_{\alpha}})
M_{\alpha}'']
\hfill\cr\hfill
\eqalign{&=[M_{\alpha}'':{\cal O}_{E_{I_{\alpha}}}]+
[{\cal O}_{E_{I_{\alpha}}}:P_{I_{\beta}}(\gamma_{I_{\alpha}})
{\cal O}_{E_{I_{\alpha}}}]
-[P_{I_{\beta}}(\gamma_{I_{\alpha}})M_{\alpha}''
:P_{I_{\beta}}(\gamma_{I_{\alpha}}){\cal O}_{E_{I_{\alpha}}}]\cr
&= [{\cal O}_{E_{I_{\alpha}}}:P_{I_{\beta}}(\gamma_{I_{\alpha}})
{\cal O}_{E_{I_{\alpha}}}]=\sum_{i_{\alpha}\in I_{\alpha}}
[{\cal O}_{E_{i_{\alpha}}}:P_{I_{\beta}}(\gamma_{i_{\alpha}})
{\cal O}_{E_{i_{\alpha}}}]\cr
&=\sum_{i_{\alpha}\in I_{\alpha}}v_{E_{I_{\alpha}}}
(P_{I_{\beta}}(\gamma_{i_{\alpha}}))=\sum_{{i_{\alpha}\in
I_{\alpha}\atop i_{\beta}\in
I_{\beta}}}r_{i_{\alpha},i_{\beta}}=r.\cr}\qquad}
$$
\hfill\hfill$\square$
\vskip 3mm

Pour chaque entier $j$, les parties
$$
\{M\in X_{I}^{0}\mid \mathop{\rm ind}\nolimits_{1}^{\prime}(M)\geq j\}
=\{M\in X_{I}^{0}\mid \mathop{\rm ind}\nolimits_{2}^{\prime\prime}(M)
\leq -r-j\}
$$
et
$$
\{M\in X_{I}^{0}\mid \mathop{\rm ind}\nolimits_{1}^{\prime}(M)\leq j\}
=\{M\in X_{I}^{0}\mid \mathop{\rm ind}\nolimits_{2}^{\prime\prime}(M)
\geq -r-j\}
$$
sont respectivement ferm\'{e}e et ouverte dans $X_{I}^{0}$.  Par
suite, pour chaque $\rho=0,1,\ldots ,r$, la partie
$X_{I_{1},I_{2};\rho}^{0}$ de $X_{I}^{0}$ form\'{e}e des $M$ tels que
$$
\mathop{\rm ind}\nolimits_{1}^{\prime\prime}(M)-
\mathop{\rm ind}\nolimits_{1}^{\prime}(M)=
\mathop{\rm ind}\nolimits_{2}^{\prime\prime}(M)-
\mathop{\rm ind}\nolimits_{2}^{\prime}(M)=\rho
$$
est localement ferm\'{e}e dans $X_{I}^{0}$ et la r\'{e}union disjointe
$$
X_{I_{1},I_{2};\leq\rho}^{0}=X_{I_{1},I_{2};0}^{0}\cup
X_{I_{1},I_{2};1}^{0}\cup\cdots\cup X_{I_{1},I_{2};\rho}^{0}
$$
est ferm\'{e}e dans $X_{I}^{0}$.  D'apr\`{e}s le lemme 1.4.2, on a en
fait $X_{I_{1},I_{2};\leq r}^{0}=X_{I}^{0}$.

Pour chaque $\rho =0,1,\ldots ,r$, l'action par translations de
$\Lambda_{I}^{0}$ sur $X_{I}^{0}$ (via le scindage $\sigma$) stabilise
$X_{I_{1},I_{2};\rho}^{0}$ et il r\'{e}sulte de ce qui pr\'{e}c\`{e}de
que le quotient $Z_{I_{1},I_{2};\rho}^{0}=X_{I_{1},I_{2};\rho}^{0}/
\sigma (\Lambda_{I}^{0})$ est une partie localement ferm\'{e}e de
$Z_{I}^{0}$ et que la r\'{e}union disjointe
$$
Z_{I_{1},I_{2};\leq\rho}^{0}=Z_{I_{1},I_{2};0}^{0}\cup
Z_{I_{1},I_{2};1}^{0}\cup\cdots\cup Z_{I_{1},I_{2};\rho}^{0}
$$
est un ferm\'{e} de $Z_{I}^{0}$, qui est en fait $Z_{I}^{0}$ tout
entier pour $\rho =r$.

On v\'{e}rifie que:
\vskip 1mm

\itemitem{-} pour chaque $\rho =0,1,\ldots ,r$, le $k$-sch\'{e}ma
r\'{e}duit $X_{I_{1},I_{2};\rho}^{0}$ est r\'{e}union disjointe de ses
parties ouvertes et ferm\'{e}es
$$\eqalign{
X_{I_{1},I_{2};\rho ,i}^{0} &=\{M\in X_{I}^{0}\mid \mathop{\rm
ind}\nolimits_{1}^{\prime}(M)= i\hbox{ et }\mathop{\rm
ind}\nolimits_{1}^{\prime\prime}(M)=\rho +i\}\cr
&=\{M\in X_{I}^{0}\mid \mathop{\rm ind}\nolimits_{2}^{\prime}(M)=
-r-\rho -i\hbox{ et }\mathop{\rm ind}\nolimits_{2}^{\prime\prime}(M)
=-r-i\},\cr}
$$
pour $i\in {\Bbb Z}$,
\vskip 1mm

\itemitem{-} pour chaque entier $i$, l'action de
$\Lambda_{I_{1}}^{0}\times\Lambda_{I_{2}}^{0} \subset\Lambda_{I}^{0}$
via $\sigma$ sur $X_{I}^{0}$ stabilise $X_{I_{1},I_{2};\rho ,i}^{0}$,
\vskip 1mm

\itemitem{-} le morphisme quotient $X_{I}^{0}\rightarrow
X_{I}^{0}/\sigma (\Lambda_{I}^{0})=Z_{I}^{0}$ induit un isomorphisme
de $X_{I_{1},I_{2};\rho ,i}^{0}/\sigma (\Lambda_{I_{1}}^{0}\times
\Lambda_{I_{2}}^{0})$ sur $Z_{I_{1},I_{2};\rho}^{0}$.
\vskip 1mm

En particulier, le rev\^{e}tement
$$
X_{I_{1},I_{2};\rho}^{0}/\sigma (\Lambda_{I_{1}}^{0}\times
\Lambda_{I_{2}}^{0})\rightarrow Z_{I_{1},I_{2};\rho}^{0}
$$
induit par le morphisme quotient $X_{I}^{0}/\sigma
(\Lambda_{I_{1}}^{0}\times \Lambda_{I_{2}}^{0})\rightarrow X_{I}^{0}/
\sigma(\Lambda_{I}^{0})=Z_{I}^{0}$ est trivial.

\thm LEMME 1.4.3
\enonce
On a
$$
\sum_{i}(-1)^{i}[H^{i}(\overline{k}\otimes_{k}Z_{I}^{0},{\cal
L}_{I_{1},I_{2}})]=\sum_{\rho =1}^{r}\sum_{i}(-1)^{i}[H^{i}_{{\rm c}}
(\overline{k}\otimes_{k}Z_{I_{1},I_{2};\rho}^{0},{\Bbb Q}_{\ell})].
$$
\endthm

\rem Preuve
\endrem
Comme $(Z_{I_{1},I_{2};\rho}^{0})_{\rho =0,1,\ldots ,r}$ est une
stratification du $k$-sch\'{e}ma projectif $Z_{I}^{0}$ en parties
localement ferm\'{e}es, on a une suite spectrale
$$
E_{\bullet}({\cal F}):~
E_{1}^{\rho ,\sigma}({\cal F})=H_{{\rm c}}^{\rho +\sigma}
(\overline{k}\otimes_{k}Z_{I_{1},I_{2};\rho}^{0},{\cal F})\Rightarrow
H^{\rho +\sigma}(\overline{k}\otimes_{k}Z_{I}^{0},{\cal F}),
$$
et donc la relation
$$
\sum_{i}(-1)^{i}[H^{i}(\overline{k}\otimes_{k}Z_{I}^{0},{\cal F})]=
\sum_{\rho =1}^{r}\sum_{i}(-1)^{i}[H_{{\rm c}}^{i}(\overline{k}
\otimes_{k}Z_{I_{1},I_{2};\rho}^{0},{\cal F})],
$$
quel que soit le ${\Bbb Q}_{\ell}$-faisceau ${\cal F}$ sur $Z_{I}^{0}$

La restriction du ${\Bbb Q}_{\ell}$-syst\`{e}me local ${\cal
L}_{I_{1},I_{2}}$ de rang $1$ \`{a} chaque partie localement
ferm\'{e}e $Z_{I_{1},I_{2};\rho}^{0}$ est triviale puisque le
rev\^{e}tement $X_{I_{1},I_{2};\rho}^{0}/\sigma
(\Lambda_{I_{1}}^{0}\times \Lambda_{I_{2}}^{0})\rightarrow
Z_{I_{1},I_{2};\rho}^{0}$ l'est, d'o\`{u} le lemme.
\hfill\hfill$\square$
\vskip 3mm

\rem Remarque $1.4.4$
\endrem
Si les suites spectrales $E_{\bullet}({\Bbb Q}_{\ell})$ et
$E_{\bullet}({\cal L}_{I_{1},I_{2}})$ ci-dessus ont le m\^{e}me terme
$E_{1}$, elles n'ont cependant pas les m\^{e}mes diff\'{e}rentielles
$d_{1}$, et donc elles diff\`{e}rent d\'{e}j\`{a} en leurs termes $E_{2}$.
\hfill\hfill$\square$
\vskip 3mm

La partie
$$
U_{I_{1},I_{2}}^{0}=\{M\in X_{I}^{0}\mid \mathop{\rm ind}
\nolimits_{1}^{\prime}(M)=0\}=\{M\in X_{I}^{0}\mid \mathop{\rm ind}
\nolimits_{2}^{\prime\prime}(M)=-r\}
$$
de $X_{I}^{0}$ est localement ferm\'{e}e et est munie d'un
$k$-morphisme
$$
\pi_{I_{1},I_{2}}^{0}:U_{I_{1},I_{2}}^{0}\rightarrow X_{I_{1}}^{0}
\times X_{I_{2}}^{-r}
$$
qui envoie $M$ sur $(M_{1}',M_{2}'')$.  Ce morphisme est
\'{e}videmment \'{e}quivariant pour les actions du sous-groupe
$\Lambda_{I_{1}}^{0}\times \Lambda_{I_{2}}^{0}\subset \Lambda_{I}^{0}$
via le scindage $\sigma$, et passe donc au quotient en un
$k$-morphisme
$$
V_{I_{1},I_{2}}^{0}=U_{I_{1},I_{2}}^{0}/(\Lambda_{I_{1}}^{0}\times
\Lambda_{I_{2}}^{0}) \rightarrow Z_{I_{1}}^{0}\times Z_{I_{2}}^{-r}
$$

\thm PROPOSITION 1.4.5 (voir [L-R] Theorem 6.1, et aussi [K-L] \S 5)
\enonce
Le $k$-morphisme $\pi_{I_{1},I_{2}}^{0} :U_{I_{1},I_{2}}^{0}
\rightarrow X_{I_{1}}^{0}\times X_{I_{2}}^{-r}$ ci-dessus est un
fibr\'{e} vectoriel $(\Lambda_{I_{1}}^{0}\times
\Lambda_{I_{2}}^{0})$-\'{e}quivariant de rang $r$ et le $k$-morphisme
induit, $V_{I_{1},I_{2}}^{0}\rightarrow Z_{I_{1}}^{0}\times
Z_{I_{2}}^{-r}$, est un fibr\'{e} vectoriel de rang $r$.
\endthm

\rem Preuve
\endrem
La fibre de $\pi_{I_{1},I_{2}}^{0}$ en un point $(M_{1},M_{2})$ de
$X_{I_{1}}^{0}\times X_{I_{2}}^{-r}$ est le sch\'{e}ma
des ${\cal O}_{F}$-r\'{e}seaux $M\subset E_{I}$, stables par
$\gamma_{I}$, qui s'ins\`{e}rent dans le diagramme commutatif
$$
\matrix{E_{I_{1}} & \hookrightarrow & E_{I} & \twoheadrightarrow &
E_{I_{2}}\cr
\cup &  & \cup &  & \cup\cr
M_{1} & \hookrightarrow & M & \twoheadrightarrow & M_{2}\cr}
$$
de ${\cal O}_{F}$-modules.  Pour tout $M$ comme ci-dessus, on a
n\'{e}cessairement $M_{1}\oplus (0)\subset M\subset E_{I_{1}}\oplus
M_{2}$, et $M/(M_{1}\oplus (0))\subset (E_{I_{1}}/M_{1})\oplus M_{2}$
est le graphe d'un homomorphisme de ${\cal O}_{F}$-modules
$$
M_{2}\rightarrow E_{I_{1}}/M_{1}
$$
qui \'{e}change la multiplication par $\gamma_{I_{2}}$ sur $M_{2}$
avec celle par $\gamma_{I_{1}}$ sur $E_{I_{1}}/M_{1}$.  La fibre de
$\pi_{I_{1},I_{2}}^{0}$ en $(M_{1},M_{2})$ peut donc \^{e}tre
identifi\'{e}e \`{a} l'espace vectoriel
$$
\mathop{\rm Hom}\nolimits_{{\cal O}_{F}[T]}(M_{2},E_{I_{1}}/M_{1}),
$$
o\`{u} $T$ agit sur $M_{2}$ par multiplication par $\gamma_{I_{2}}$ et
sur $E_{I_{1}}/M_{1}$ par multiplication par $\gamma_{I_{1}}$.  Une
variante de cet argument o\`{u} on remplace $k$ (ou $k'$) par une
$k$-alg\`{e}bre arbitraire, montre que $\pi_{I_{1},I_{2}}^{0}$ est un
fibr\'{e} vectoriel g\'{e}n\'{e}ralis\'{e} au sens de Grothendieck.
Comme $\gamma_{I}$ est topologiquement nilpotent, on a encore
$$
\mathop{\rm Hom}\nolimits_{{\cal O}_{F}[T]}(M_{2},E_{I_{1}}/M_{1})
=\mathop{\rm Hom}\nolimits_{{\cal O}_{F}[[T]]}(M_{2},E_{I_{1}}/M_{1}).
$$

Pour conclure la preuve de la proposition, il suffit donc de
d\'{e}montrer que le rang de l'espace vectoriel
$$
\mathop{\rm Hom}\nolimits_{{\cal O}_{F}[[T]]}(M_{2},E_{I_{1}}/M_{1})
$$
ne d\'{e}pend pas de $(M_{1},M_{2})$ et est \'{e}gal \`{a}
$r$.  Pour cela on va montrer que le complexe de
$k$-espaces vectoriels
$$
R\mathop{\rm Hom}\nolimits_{{\cal O}_{F}[[T]]}(M_{2},E_{I_{1}}/M_{1}),
$$
est concentr\'{e} en degr\'{e} $0$ et a une caract\'{e}ristique
d'Euler-Poincar\'{e} \'{e}gale \`{a} $r$.

On a un triangle distingu\'{e}
$$
R\mathop{\rm Hom}\nolimits_{{\cal O}_{F}[[T]]}(M_{2},M_{1})\rightarrow
R\mathop{\rm Hom}\nolimits_{{\cal O}_{F}[[T]]}(M_{2},E_{I_{1}})\rightarrow
R\mathop{\rm Hom}\nolimits_{{\cal O}_{F}[[T]]}(M_{2},E_{I_{1}}/M_{1})
\rightarrow
$$
o\`{u} $R\mathop{\rm Hom}\nolimits_{{\cal O}_{F}[[T]]}
(M_{2},E_{I_{1}}) =(0)$ puisque la multiplication par $P_{I_{2}}(T)\in
{\cal O}_{F}[[T]]$ est l'endomorphisme nul de $M_{2}$ et un
automorphisme de $E_{I_{1}}$.  Il ne nous reste plus qu'\`{a}
d\'{e}montrer que le complexe $R\mathop{\rm Hom}\nolimits_{{\cal
O}_{F}[[T]]}(M_{2},M_{1})$ est concentr\'{e} en degr\'{e} $1$ et a une
caract\'{e}ristique d'Euler-Poincar\'{e} \'{e}gale \`{a} $-r$.

Adoptons un point de vue plus g\'{e}om\'{e}trique en introduisant le
germe formel de surface ${\cal S}=\mathop{\rm Spf}({\cal O}_{F}[[T]])$
sur $k$ et, pour $\alpha =1,2$, la courbe
$$
\iota_{\alpha}:{\cal C}_{I_{\alpha}}=\mathop{\rm
Spf}(A_{I_{\alpha}})\hookrightarrow {\cal S}
$$
d'\'{e}quation $P_{I_{\alpha}}(T)=0$, trac\'{e}e sur cette surface, et
le ${\cal O}_{{\cal C}_{I_{\alpha}}}$-Module coh\'{e}rent ${\cal
M}_{\alpha}$, sans torsion et partout de rang g\'{e}n\'{e}rique $1$,
d\'{e}fini par $M_{I_{\alpha}}$. On a bien entendu
$$
R\mathop{\rm Hom}\nolimits_{{\cal O}_{F}[[T]]}(M_{2},M_{1})
=R\mathop{\rm Hom}\nolimits_{{\cal O}_{{\cal S}}}(\iota_{2,\ast}{\cal
M}_{2},\iota_{1,\ast}{\cal M}_{1})
$$
et $R\mathop{\rm Hom}\nolimits_{{\cal O}_{{\cal S}}}
(\iota_{2,\ast}{\cal M}_{2},\iota_{1,\ast}{\cal M}_{1})$ est aussi la
fibre de $R\mathop{{\cal H}{\it om}}\nolimits_{{\cal O}_{{\cal
S}}}(\iota_{2,\ast}{\cal M}_{2},\iota_{1,\ast}{\cal M}_{1})$ en
l'origine $s=(\varpi_{F}=0,T=0)$ de ${\cal S}$ puisque l'intersection
des deux courbes ${\cal C}_{I_{1}}$ et ${\cal C}_{I_{2}}$ est par
hypoth\`{e}se support\'{e}e par $s$.

Comme ${\cal S}$ est r\'{e}guli\`{e}re, on a $\mathop{\rm Ext}
\nolimits_{{\cal S}}^{i}(\iota_{2,\ast}{\cal M}_{2},
\iota_{1,\ast}{\cal M}_{1})=(0)$ pour $i\notin\{0,1,2\}$.  On a aussi
$$
\mathop{\rm Hom}\nolimits_{{\cal O}_{{\cal S}}}(\iota_{2,\ast}{\cal
M}_{2}, \iota_{1,\ast}{\cal M}_{1})= \mathop{\rm Hom}\nolimits_{{\cal
O}_{{\cal C}_{I_{1}}}}(\iota_{1}^{\ast}\iota_{2,\ast}{\cal M}_{2},
{\cal M}_{1})=(0)
$$
puisque $\iota_{1}^{\ast}\iota_{2,\ast}{\cal M}_{2}$ est de torsion et
${\cal M}_{1}$ est sans torsion.

Soient $\rho_{2}:\widetilde{{\cal C}}_{I_{2}}=\mathop{\rm Spf}({\cal
O}_{E_{I_{1}}})\rightarrow {\cal C}_{I_{2}}$ la normalisation de la
courbe ${\cal C}_{I_{2}}$ et $\widetilde{{\cal M}}_{2}$ le ${\cal
O}_{\widetilde{{\cal C}}_{I_{2}}}$-Module libre de rang $1$ d\'{e}fini
par le ${\cal O}_{E_{I_{2}}}$-module $\widetilde{M}_{2}={\cal
O}_{E_{I_{2}}}M_{2}\subset E_{I_{2}}$.  On a un homomorphisme injectif
de ${\cal O}_{{\cal S}}$-Modules $\iota_{2,\ast}{\cal M}_{2}
\hookrightarrow \iota_{2,\ast}\rho_{2,\ast}\widetilde{{\cal M}}_{2}$,
et donc un \'{e}pimorphisme de $k$-espaces vectoriels
$$
\mathop{\rm Ext}\nolimits^{2}\nolimits_{{\cal O}_{{\cal S}}}
(\iota_{2,\ast} \rho_{2,\ast}\widetilde{{\cal M}}_{2},
\iota_{1,\ast}{\cal M}_{1}) \twoheadrightarrow \mathop{\rm
Ext}\nolimits^{2}\nolimits_{{\cal O}_{{\cal S}}} (\iota_{2,\ast}{\cal
M}_{2},\iota_{1,\ast}{\cal M}_{1}).
$$
Or, par dualit\'{e} de Grothendieck, on a
$$
\mathop{\rm Ext}\nolimits_{{\cal O}_{{\cal S}}}^{i}(\iota_{2,\ast}
\rho_{2,\ast}\widetilde{{\cal M}}_{2}, \iota_{1,\ast}{\cal M}_{1})
=\mathop{\rm Ext}\nolimits^{i}\nolimits_{{\cal O}_{\widetilde{{\cal
C}}_{I_{2}}}}(\widetilde{{\cal M}}_{2},L(\iota_{2}\circ\rho_{2})^{!}
\iota_{1,\ast}{\cal M}_{1})
$$
et
$$
L(\iota_{2}\circ\rho_{2})^{!}\iota_{1,\ast}{\cal M}_{1}=
\omega_{\widetilde{{\cal C}}_{I_{2}}}\otimes
(\iota_{2}\circ\rho_{2})^{\ast}\omega_{{\cal S}}^{\otimes -1}\otimes
L(\iota_{2}\circ\rho_{2})^{\ast}\iota_{1,\ast}{\cal M}_{1}[-1]
$$
puisque les complexes dualisants de ${\cal S}$ et $\widetilde{{\cal
C}}_{I_{2}}$ sont de la forme $\omega_{S}[2]$ et
$\omega_{\widetilde{{\cal C}}_{I_{2}}}[1]$ o\`{u} les Modules
$\omega_{S}$ et $\omega_{\widetilde{{\cal C}}_{I_{2}}}$ sont
inversibles.  Par suite, $\mathop{\rm
Ext}\nolimits^{2}\nolimits_{{\cal S}}(\iota_{2,\ast}
\rho_{2,\ast}\widetilde{{\cal M}}_{2},\iota_{1,\ast}{\cal M}_{1})=(0)$
et on a montr\'{e} que $R\mathop{\rm Hom}\nolimits_{{\cal O}_{{\cal
S}}} (\iota_{2,\ast}{\cal M}_{2},\iota_{1,\ast}{\cal M}_{1})$ est
concentr\'{e} en degr\'{e} $1$.

On termine la d\'{e}monstration de la proposition \`{a} l'aide
du corollaire 1.4.7 ci-dessous.

\hfill\hfill$\square$
\vskip 3mm

\thm PROPOSITION 1.4.6 (Deligne)
\enonce
Soient ${\cal X}=\mathop{\rm Spf}(\widehat{{\cal O}}_{X,x})$ le
compl\'{e}t\'{e} formel d'une vari\'{e}t\'{e} lisse $X$ sur un corps
$k$ en un point ferm\'{e} $x$, suppos\'{e} rationnel sur $k$ pour
simplifier.  Soient $\iota_{1}:{\cal Y}_{1} \hookrightarrow {\cal X}$
et $\iota_{2}:{\cal Y}_{2}\hookrightarrow {\cal X}$ deux
sous-sch\'{e}mas formels ferm\'{e}s de ${\cal X}$, de dimensions pures
$\mathop{\rm dim}({\cal Y}_{1})$ et $\mathop{\rm dim}({\cal Y}_{2})$
respectivement, tels que
$$
({\cal Y}_{1}\cap {\cal Y}_{2})_{{\rm red}}=\{x\}
$$
et
$$
\mathop{\rm dim}({\cal Y}_{1})+\mathop{\rm dim}({\cal Y}_{2})
=\mathop{\rm dim}({\cal X}).
$$

Soient $K_{1}$ et $K_{2}$ deux complexes born\'{e}s de ${\cal
O}_{{\cal X}}$-Modules quasi-coh\'{e}rents plats dont les faisceaux
de cohomologie sont support\'{e}s par ${\cal Y}_{1}$ et ${\cal Y}_{2}$
respectivement.

Alors, le complexe $K_{1}\otimes_{{\cal O}_{{\cal X}}}K_{2}$ est
support\'{e} par le point ferm\'{e} $x$ de ${\cal X}$ et sa fibre
{\rm (}na\"{\i}ve{\rm )} en
ce point est un complexe born\'{e} de $k$-espaces vectoriels dont la
caract\'{e}ristique d'Euler-Poincar\'{e} est \'{e}gale \`{a}
$$
m({\cal Y}_{1},{\cal Y}_{2})\mathop{\rm rang}(K_{1})\mathop{\rm
rang}(K_{2})
$$
o\`{u} $m({\cal Y}_{1},{\cal Y}_{2})$ est la multiplicit\'{e}
d'intersection de ${\cal Y}_{1}$ et ${\cal Y}_{2}$ et o\`{u} $\mathop{\rm
rang}(K_{\alpha})$ est le rang g\'{e}n\'{e}rique de $K_{\alpha}$
pour $\alpha =1,2$.
\endthm

\rem Preuve
\endrem
Voir [De 1] Th\'{e}or\`{e}me 2.3.8 (iii).  Les arguments de Deligne
peuvent \^{e}tre facilement modifi\'{e}s pour \'{e}viter l'utilisation
de la cohomologie $\ell$-adique.
\hfill\hfill$\square$
\vskip 3mm

\thm COROLLAIRE 1.4.7
\enonce
Supposons de plus que ${\cal Y}_{2}$ soit d'intersection compl\`{e}te
dans ${\cal X}$ et, pour $\alpha =1,2$, soit ${\cal M}_{\alpha}$ un ${\cal
O}_{{\cal Y}_{\alpha}}$-Module coh\'{e}rent de rang g\'{e}n\'{e}rique
$r_{\alpha}$. Alors le complexe de $k$-espaces vectoriels
$$
R\mathop{\rm Hom}\nolimits_{{\cal O}_{{\cal X}}}(\iota_{2,\ast}{\cal
M}_{2},\iota_{1,\ast}{\cal M}_{1})
$$
est born\'{e} et \`{a} cohomologie de dimension finie, et sa
caract\'{e}ristique d'Euler-Poincar\'{e} est \'{e}gale \`{a}
$$
(-1)^{\mathop{\rm dim}({\cal Y}_{2})}m({\cal Y}_{1},{\cal Y}_{2})
r_{1}r_{2}.
$$
\endthm

\rem Preuve
\endrem
Comme ${\cal X}$ est r\'{e}gulier, on peut trouver, pour $\alpha
=1,2$, une r\'{e}solution finie
$$
K_{\alpha}\rightarrow\iota_{\alpha ,\ast}{\cal M}_{\alpha}\rightarrow
0
$$
par des ${\cal O}_{{\cal X}}$-Modules libres de rang fini. Le complexe
d'espaces vectoriels
$$
R\mathop{\rm Hom}\nolimits_{{\cal O}_{{\cal
X}}} (\iota_{2,\ast}{\cal M}_{2},\iota_{1,\ast}{\cal M}_{1})=
\mathop{\rm Hom}\nolimits_{{\cal O}_{{\cal
X}}}(K_{2},K_{1})
$$
est alors la fibre en $s$ du complexe $K_{2}^{\vee}\otimes_{{\cal
O}_{{\cal X}}}K_{1}$, o\`{u}
$$
K_{2}^{\vee}=\mathop{{\cal H}{\it om}}\nolimits_{{\cal O}_{{\cal X}}}
(K_{2},{\cal O}_{{\cal X}})
$$
est aussi un complexe born\'{e} de ${\cal O}_{{\cal X}}$-Modules
libres de rang fini.

Par dualit\'{e} de Grothendieck, on a
$$\eqalign{
K_{2}^{\vee}=R\mathop{{\cal H}{\it om}}\nolimits_{{\cal O}_{{\cal
X}}}(\iota_{2,\ast}{\cal M}_{2},{\cal O}_{{\cal X}})&=\iota_{2,\ast}
R\mathop{{\cal H}{\it om}}\nolimits_{{\cal O}_{{\cal Y}_{2}}}
({\cal M}_{2},\iota_{2}^{!}{\cal O}_{{\cal X}})\cr
&\cong\iota_{2,\ast}
R\mathop{{\cal H}{\it om}}\nolimits_{{\cal O}_{{\cal Y}_{2}}}
({\cal M}_{2},{\cal O}_{{\cal Y}_{2}})[\mathop{\rm dim}({\cal
X})-\mathop{\rm dim}({\cal Y}_{2})]\cr}
$$
et, en particulier,
$$
\mathop{\rm rang}(K_{2}^{\vee})=(-1)^{\mathop{\rm dim}({\cal Y}_{1})}r_{2}.
$$

Il ne reste plus qu'\`{a} appliquer la proposition.
\hfill\hfill$\square$
\vskip 3mm

\rem Preuve du th\'{e}or\`{e}me $1.4.1$
\endrem
On stratifie $U_{I_{1},I_{2}}^{0}$ par les parties localement
ferm\'{e}es $X_{I_{1},I_{2};\rho ,0}^{0}$ pour $\rho =0,1,\ldots ,r$,
et donc $V_{I_{1},I_{2}}^{0}$ par les parties localement ferm\'{e}es
$X_{I_{1},I_{2};\rho ,0}^{0}/\sigma (\Lambda_{I_{1}}^{0}\times
\Lambda_{I_{2}}^{0})\cong Z_{I_{1},I_{2};\rho}^{0}$ pour $\rho
=0,1,\ldots ,r$.

On a alors une suite spectrale
$$
E_{1}^{\rho ,\sigma}=H_{{\rm c}}^{\rho +\sigma} (\overline{k}
\otimes_{k}Z_{I_{1},I_{2};\rho}^{0},{\Bbb Q}_{\ell})\Rightarrow
H_{{\rm c}}^{\rho +\sigma}(\overline{k} \otimes_{k}
V_{I_{1},I_{2}}^{0},{\Bbb Q}_{\ell}),
$$
et la relation
$$
\sum_{\rho =1}^{r}\sum_{i}(-1)^{i}[H^{i}_{{\rm c}}
(\overline{k}\otimes_{k}Z_{I_{1},I_{2};\rho}^{0},{\Bbb Q}_{\ell})]=
\sum_{i}(-1)^{i}[H_{{\rm c}}^{i}(\overline{k}\otimes_{k}
V_{I_{1},I_{2}}^{0},{\Bbb Q}_{\ell})].
$$
Mais on a
$$
H_{{\rm c}}^{i}(\overline{k}\otimes_{k}V_{I_{1},I_{2}}^{0},{\Bbb
Q}_{\ell})=H^{i-2r}(\overline{k}\otimes_{k} (Z_{I_{1}}^{0}\times
Z_{I_{2}}^{-r}),{\Bbb Q}_{\ell})(-r)
$$
d'apr\`{e}s la proposition 1.4.5, et $Z_{I_{2}}^{-r}$ est isomorphe
\`{a} $Z_{I_{2}}^{0}$, d'o\`{u} la conclusion d'apr\`{e}s le lemme 1.4.3.
\hfill\hfill$\square$

\rem Remarque $1.4.8$
\endrem
Dans la situation arithm\'{e}tique de la section 1.3, on v\'{e}rifie
que
$$
\mathop{\rm ind}\nolimits_{\alpha}^{\prime}(\mathop{\rm Frob}
\nolimits_{{\cal X}_{I}/{\Bbb F}_{q}}(M))=-r-\mathop{\rm ind}
\nolimits_{\alpha}^{\prime\prime}(M)\hbox{ et }\mathop{\rm ind}
\nolimits_{\alpha}^{\prime\prime}(\mathop{\rm Frob} \nolimits_{{\cal
X}_{I}/{\Bbb F}_{q}}(M))=-r-\mathop{\rm ind}
\nolimits_{\alpha}^{\prime}(M)
$$
pour $\alpha =1,2$ et tout point $M$ de ${\cal X}_{I}$.  En effet,
pour $\{\alpha ,\beta\}=\{1,2\}$, on a
$$
\varepsilon^{n_{I}-1}{{\rm d}P_{I}\over {\rm d}T}
(\gamma_{I_{\alpha}})=\varepsilon^{n_{I_{\beta}}}P_{I_{\beta}}
(\gamma_{I_{\alpha}})\times\varepsilon^{n_{I_{\alpha}}-1}{{\rm d}
P_{I_{\alpha}}\over {\rm d}T}(\gamma_{I_{\alpha}})
$$
avec $v_{E_{I_{\alpha}}}(\varepsilon^{n_{I_{\beta}}}P_{I_{\beta}}
(\gamma_{I_{\alpha}}))=r$, de sorte que
$$
(A_{I_{1}}\times A_{I_{2}})^{\perp_{I}}=
\varepsilon^{n_{2}}P_{I_{2}}(\gamma_{I_{1}})A_{I_{1}}\times
\varepsilon^{n_{1}}P_{I_{1}}(\gamma_{I_{2}})A_{I_{2}}\subset A_{I}\subset
A_{I_{1}}\times A_{I_{2}}.
$$
Par suite, chaque partie localement ferm\'{e}e
$Z_{I_{1},I_{2};\rho}^{0}$ de $Z_{I}^{0}$ est stable par
l'endomorphisme de Frobenius $\mathop{\rm Frob} \nolimits_{{\cal
Z}_{I}^{0}/{\Bbb F}_{q}}$ et est munie d'une structure rationnelle sur
${\Bbb F}_{q}$ que l'on note bien entendu ${\cal
Z}_{I_{1},I_{2};\rho}^{0}$.  De plus, le rel\`{e}vement de
$\mathop{\rm Frob}\nolimits_{{\cal Z}_{I_{1},I_{2};\rho}^{0}/{\Bbb
F}_{q}}$ \`{a} ${\cal L}_{I_{1},I_{2}}|{\cal
Z}_{I_{1},I_{2};\rho}^{0}={\Bbb Q}_{\ell}$ est la multiplication par
$(-1)^{r-\rho}$ puisque l'on a
$$
\mathop{\rm ind} \nolimits_{\alpha}^{\prime}(\sigma (\lambda )
\cdot M)-\mathop{\rm ind}\nolimits_{\alpha}^{\prime}(M)=
\mathop{\rm ind} \nolimits_{\alpha}^{\prime\prime}(\sigma (\lambda )
\cdot M)-\mathop{\rm ind}\nolimits_{\alpha}^{\prime\prime}(M)
=\sum_{i\in I_{\alpha}}\lambda_{i}
$$
pour $\alpha =1,2$, tout $M\in X_{I}^{0}$ et tout
$\lambda\in\Lambda_{I}^{0}$, et que l'\'{e}galit\'{e}
$$
\mathop{\rm Frob}\nolimits_{{\cal X}_{I}/{\Bbb F}_{q}}(M)=\sigma
(\lambda_{M})\cdot M
$$
pour $M\in X_{I_{1},I_{2};\rho}$ et $\lambda_{M}\in\Lambda_{I}^{0}$
force
$$
\sum_{i\in I_{\alpha}}\lambda_{M,i}=-r-\mathop{\rm ind}
\nolimits_{\alpha}^{\prime}(M)-\mathop{\rm ind}
\nolimits_{\alpha}^{\prime\prime}(M)
$$
\`{a} \^{e}tre de m\^{e}me parit\'{e} que $r-\rho =r-\mathop{\rm ind}
\nolimits_{\alpha}^{\prime}(M)+\mathop{\rm ind}
\nolimits_{\alpha}^{\prime\prime}(M)$.  On a donc
$$\eqalign{
\mathop{\rm O}\nolimits^{I,\kappa_{I_{1},I_{2}}}&=\mathop{\rm
tr}(\mathop{\rm Frob}\nolimits_{{\cal Z}_{I}^{0}/{\Bbb F}_{q}},
R\Gamma (\overline{{\Bbb F}}_{p}\otimes_{{\Bbb F}_{q}}{\cal Z}_{I}^{0},
{\cal L}_{I_{1},I_{2}}))\cr
&= \sum_{\rho =0}^{r}(-1)^{r-\rho}\mathop{\rm tr}(\mathop{\rm
Frob}\nolimits_{{\cal Z}_{I_{1},I_{2};\rho}^{0}/{\Bbb F}_{q}},
R\Gamma (\overline{{\Bbb F}}_{p}\otimes_{{\Bbb F}_{q}}{\cal
Z}_{I_{1},I_{2};\rho}^{0},{\Bbb Q}_{\ell})).\cr}
$$
\hfill\hfill$\square$
\vskip 5mm

\centerline{2.  LES FIBRES DE SPRINGER COMME REV\^{E}TEMENTS}
\centerline{DE JACOBIENNES COMPACTIFI\'{E}ES}
\vskip 2mm

\section{2.1}{La courbe $C_{I}$ et son sch\'{e}ma de Picard
compactifi\'{e} de $C_{I}$}

Dans la situation de la section 1.1, le sch\'{e}ma formel $\mathop{\rm
Spf}(A_{I})$ est un germe formel de courbe plane dont la famille des
branches irr\'{e}ductibles est $(\mathop{\rm Spf}(A_{i}))_{i\in I}$ et
dont le normalis\'{e} est le sch\'{e}ma formel semi-local $\mathop{\rm
Spf}(\widetilde{A}_{I})= \coprod_{i\in I}\mathop{\rm Spf}({\cal
O}_{E_{i}})$. {\it On suppose dans la suite que le nombre
d'\'{e}l\'{e}ments du corps $k$ est au moins \'{e}gal au nombre
d'\'{e}l\'{e}ments de l'ensemble fini $I$.}

\thm PROPOSITION 2.1.1
\enonce
Il existe une courbe projective et g\'{e}om\'{e}triquement int\`{e}gre
$C_{I}$ sur $k$, munie d'un $k$-point $c_{I}$, ayant les
propri\'{e}t\'{e}s suivantes:

\decale{\rm (1)} $C_{I}$ est lisse sur $k$ en dehors de
$c_{I}$,

\decale{\rm (2)} le compl\'{e}t\'{e} de l'anneau local de $C_{I}$ en
$c_{I}$ est isomorphe \`{a} $A_{I}$,

\decale{\rm (3)} la normalis\'{e}e $\widetilde{C}_{I}$ de $C_{I}$ est
isomorphe \`{a} la droite projective standard ${\Bbb P}_{k}^{1}$ sur
$k$.
\endthm

Pour une telle courbe $C_{I}$, son morphisme de normalisation
$\pi_{I}:\widetilde{C}_{I}\rightarrow C_{I}$ est un isomorphisme
au-dessus de $C_{I}\setminus\{c_{I}\}$, et $\pi_{I}^{-1}(c_{I})\subset
\widetilde{C}_{I}$ est l'ensemble des branches de $\mathop{\rm
Spf}(A_{I})$; pour chaque $i\in I$, on notera $\widetilde{c}_{i}$ le
point de $\pi_{I}^{-1}(c_{I})$ correspondant \`{a} la branche
$\mathop{\rm Spf} (A_{i})$.

\rem Preuve
\endrem
On fixe arbitrairement une injection $\iota :I\hookrightarrow k$ et,
pour chaque $i\in I$, on fixe arbitrairement une uniformisante
$\varpi_{E_{i}}$ de ${\cal O}_{E_{i}}$. On plonge $k[x]$ dans
${\cal O}_{E_{i}}$ en envoyant $x$ sur $\iota (i)+\varpi_{E_{i}}$. On
en d\'{e}duit un plongement de $k$-alg\`{e}bres de $k[x]$ dans
${\cal O}_{E_{I}}=\widetilde{A}_{I}$ qui identifie $\widetilde{A}_{I}$ au
compl\'{e}t\'{e} de l'anneau semi-local de la droite affine ${\Bbb
A}_{k}^{1}=\mathop{\rm Spec}(k[x])$ en l'ensemble fini de points
$\iota (I)$.

Consid\'{e}rons alors la $k$-alg\`{e}bre $B_{I}$ d\'{e}finie par le
carr\'{e} cart\'{e}sien
$$
\matrix{ B_{I} & \subset & A_{I} \cr
\noalign{\smallskip}
\cap &\square & \cap \cr
\noalign{\smallskip}
k[x] & \subset & \widetilde{A}_{I}&.\cr}
$$
Elle est int\`{e}gre, de type fini sur $k$ et de dimension $1$,
l'inclusion $B_{I}\hookrightarrow A_{I}$ induit un isomorphisme du
compl\'{e}t\'{e} de $B_{I}$ le long de son id\'{e}al maximal ${\frak
m}_{I}\cap B_{I}$ sur $A_{I}$ et l'inclusion $B_{I}\hookrightarrow
k[x]$ fait de $k[x]$ une $B_{I}$-alg\`{e}bre finie.  En effet, d'une
part on a
$$
\widetilde{A}_{I}=k[x]+{\frak a}_{I}^{n+1},~\forall n\in {\Bbb N},
$$
et ${\frak a}_{I}\subset A_{I}\subset\widetilde{A}_{I}$, de sorte que
$$
A_{I}=B_{I}+{\frak a}_{I}^{n+1},~\forall n\in {\Bbb N},
$$
et d'autre part
$$
k[x]\cap {\frak a}_{I}\subset B_{I}\subset k[x]
$$
est l'id\'{e}al principal engendr\'{e} par
$$
\prod_{i\in I}(x-\iota (i))^{2\delta_{i}+\sum_{j\in
I\setminus\{i\}}r_{ij}}.
$$

On peut donc effectuer la {\og}{somme amalgam\'{e}e}{\fg} de ${\Bbb
A}_{k}^{1}$ et de $\mathop{\rm Spf}(A_{I})$ le long de $\mathop{\rm
Spf}(\widetilde{A}_{I})$; c'est par d\'{e}finition le $k$-sch\'{e}ma
affine $\mathop{\rm Spec}(B_{I})$.  Bien s\^{u}r, le morphisme fini
${\Bbb A}_{k}^{1}\rightarrow \mathop{\rm Spec}(B_{I})$, induit par
l'inclusion $B_{I}\subset k[x]$, envoie le sous-ensemble fini $\iota
(I)\subset {\Bbb A}_{k}^{1}$ sur un unique $k$-point $c_{I}$ de
$\mathop{\rm Spec}(B_{I})$ et il induit un isomorphisme de ${\Bbb
A}_{k}^{1}\setminus\iota (I)$ sur $\mathop{\rm Spec}(B_{I})
\setminus\{c_{I}\}$.

On d\'{e}finit la courbe g\'{e}om\'{e}triquement int\`{e}gre et
projective $C_{I}$ sur $k$ en recollant $\mathop{\rm Spec}(B_{I})$ et
${\Bbb P}_{k}^{1} \setminus\iota (I)$ le long de leur ouvert commun
$\mathop{\rm Spec} (B_{I})\setminus\{c_{I}\}\cong {\Bbb
A}_{k}^{1}\setminus\iota (I)$.
\hfill\hfill$\square$
\vskip 3mm

Le genre arithm\'{e}tique $\mathop{\rm dim}\nolimits_{k}H^{1}(C_{I},
{\cal O}_{C_{I}})=1-\chi (C_{I},{\cal O}_{C_{I}})$ de $C_{I}$ est
\'{e}gal \`{a} $\delta_{I}$.  Pour tout ${\cal O}_{C_{I}}$-Module
coh\'{e}rent ${\cal M}$, on pose
$$
\mathop{\rm deg}({\cal M})=\chi (C_{I},{\cal M})-
\mathop{\rm rang}({\cal M})\chi (C_{I},{\cal O}_{C_{I}})
$$
o\`{u} $\mathop{\rm rang}({\cal M})$ est le rang g\'{e}n\'{e}rique de
${\cal M}$. Pour tout ${\cal O}_{C_{I}}$-Module inversible ${\cal L}$,
$\pi_{I}^{\ast}{\cal L}$ est un ${\cal O}_{ \widetilde{C}_{I}}$-Module
inversible et on a
$$
\mathop{\rm deg}({\cal L})=\mathop{\rm deg}(\pi_{I}^{\ast}{\cal L})
$$
et
$$
\mathop{\rm deg}({\cal L}\otimes_{{\cal O}_{C_{I}}}{\cal M})=
\mathop{\rm rang}({\cal M})\mathop{\rm deg}({\cal L})+ \mathop{\rm
deg}({\cal M}).
$$
\vskip 3mm

Soit $P_{I}=\mathop{\rm Pic}\nolimits_{C_{I}/k}$ le $k$-sch\'{e}ma en
groupes de Picard de $C_{I}$.  Pour toute extension $k'$ de $k$, ses
$k'$-points sont les classes d'isomorphie de ${\cal
O}_{k'\otimes_{k}C_{I}}$-Modules inversibles (avec la multiplication
d\'{e}finie par le produit tensoriel).  Ce $k$-sch\'{e}ma est lisse de
dimension $\delta_{I}$.  Ses composantes connexes sont les
sous-$k$-sch\'{e}mas $P_{I}^{d}$, $d\in {\Bbb Z}$, d\'{e}coup\'{e}s
par le degr\'{e} du Module inversible universel, et elles sont en fait
g\'{e}om\'{e}triquement connexes.  La composante neutre $P_{I}^{0}$ de
$P_{I}$ est quasi-projective.

Soit $\overline{P}_{I}=\mathop{\overline{\rm Pic}}
\nolimits_{C_{I}/k}$ le $k$-sch\'{e}ma de Picard compactifi\'{e} (cf.
[A-K 2]) dont les $k'$-points sont les classes d'isomorphie de ${\cal
O}_{k'\otimes_{k}C_{I}}$-Modules coh\'{e}rents sans torsion de rang
g\'{e}n\'{e}rique $1$.  Par d\'{e}finition, $P_{I}$ est un ouvert de
$\overline{P}_{I}$ et l'action par translation de $P_{I}$ sur
lui-m\^{e}me se prolonge en une action de $P_{I}$ sur
$\overline{P}_{I}$ (encore d\'{e}finie par produit tensoriel).  On a
aussi un d\'{e}coupage en parties ouvertes et ferm\'{e}es
$$
\overline{P}_{I}=\coprod_{d\in {\Bbb Z}}\overline{P}{}_{I}^{d}
$$
par le degr\'{e} du Module sans torsion universel, avec bien entendu
$$
P_{I}^{d}=P_{I}\cap\overline{P}{}_{I}^{d}
$$
et
$$
P_{I}^{d}\cdot\overline{P}{}_{I}^{e}=\overline{P}{}_{I}^{d+e}
$$
quels que soient les entiers $d,e$. D'apr\`{e}s Mayer et Mumford (cf.
[A-K 2] et [A-K 3]), chaque composante $\overline{P}{}_{I}^{d}$ est un
$k$-sch\'{e}ma projectif.
\vskip 2mm

\thm TH\'{E}OR\`{E}ME 2.1.2 (Altman, Iarrobino, Kleiman, [A-I-K]
Corollary (7); Rego, [Re] Theorem A)
\enonce
Chaque composante $\overline{P}{}_{I}^{d}$ de $\overline{P}{}_{I}$ est
g\'{e}om\'{e}triquement int\`{e}gre et localement d'intersection
compl\`{e}te de dimension $\delta_{I}$.
\hfill\hfill$\square$
\endthm

\section{2.2}{Lien entre fibres de Springer et sch\'{e}mas de Picard
compactifi\'{e}s}

Faisons le lien entre les $k$-sch\'{e}mas $X_{I}$ et
$\overline{P}_{I}$. On a la suite exacte
$$
0\rightarrow H^{0}(C_{I},{\Bbb G}_{{\rm m}})\rightarrow
H^{0}(\widetilde{C}_{I},{\Bbb G}_{{\rm m}})
\rightarrow H^{0}(C_{I},\pi_{I\ast}{\Bbb G}_{{\rm m}}/{\Bbb G}_{{\rm m}})
\rightarrow H^{1}(C_{I},{\Bbb G}_{{\rm m}})
\rightarrow H^{1}(\widetilde{C}_{I},{\Bbb G}_{{\rm m}})
$$
dont la fl\`{e}che de co-bord identifie les $k$-sch\'{e}mas en groupes
$G_{I}^{0}=H^{0}(C_{I},\pi_{I\ast}{\Bbb G}_{{\rm m}}/{\Bbb G}_{{\rm
m}})$ et $P_{I}^{0}=\mathop{\rm Ker}(H^{1}(C_{I},{\Bbb G}_{{\rm m}})
\rightarrow H^{1}(\widetilde{C}_{I},{\Bbb G}_{{\rm m}}))$ puisque
$\widetilde{C}_{I}$ est une droite projective sur $k$.  On prolonge
cette identification en un $k$-\'{e}pimorphisme de $k$-sch\'{e}mas en
groupes
$$
G_{I}\twoheadrightarrow P_{I}
$$
en envoyant $x\in E_{I}^{\times}/A_{I}^{\times}$ sur le ${\cal
O}_{C_{I}}$-Module inversible ${\cal L}$ obtenu en recollant ${\cal
O}_{C_{I}\setminus\{c_{I}\}}$ et $A_{I}$ le long de $\mathop{\rm Spec}
(E_{I})=(C_{I}\setminus\{c_{I}\})\times_{C_{I}}\mathop{\rm
Spec}(A_{I})$ \`{a} l'aide de la multiplication par $x$.  Cet
\'{e}pimorphisme n'est autre que la fl\`{e}che
$H_{\{c_{I}\}}^{1}(C_{I},{\Bbb G}_{{\rm m}}) \rightarrow
H^{1}(C_{I},{\Bbb G}_{{\rm m}})$ qui s'ins\`{e}re dans la suite exacte
longue
$$\eqalign{
(1)=\,&H_{\{c_{I}\}}^{0}(C_{I},{\Bbb G}_{{\rm m}})
\rightarrow H^{0}(C_{I},{\Bbb G}_{{\rm m}})\rightarrow
H^{0}(C_{I}\setminus\{c_{I}\},{\Bbb G}_{{\rm m}})\cr
\rightarrow\,&H_{\{c_{I}\}}^{1}(C_{I},{\Bbb G}_{{\rm m}})
\rightarrow H^{1}(C_{I},{\Bbb G}_{{\rm m}})\rightarrow
H^{1}(C_{I}\setminus\{c_{I}\},{\Bbb G}_{{\rm m}})=(1)\cr}
$$
et son noyau est donc le groupe discret
$$\diagram{
H^{0}(C_{I}\setminus\{c_{I}\},{\Bbb G}_{{\rm m}})/H^{0}(C_{I},{\Bbb
G}_{{\rm m}})&\subset &E_{I}^{\times}/A_{I}^{\times}\cr
||&&\downarrow\cr
H^{0}(\widetilde{C}_{I} \setminus
\pi_{I}^{-1}(c_{I}),{\Bbb G}_{{\rm m}})/ H^{0}(\widetilde{C}_{I},{\Bbb
G}_{\rm m})&\subset&E_{I}^{\times}/{\cal O}_{E_{I}}^{\times}\cr}
$$
des diviseurs de degr\'{e} $0$ sur $\widetilde{C}_{I}$ qui sont
support\'{e}s par le ferm\'{e} r\'{e}duit $\pi_{I}^{-1}(c_{I})_{{\rm
red}}=\{\widetilde{c}_{i}\mid i\in I\}$, groupe que l'on identifie
\`{a} $\Lambda_{I}^{0}$ par la fl\`{e}che $\lambda\mapsto\sum_{i\in
I}\lambda_{i} [\widetilde{c}_{i}]$.  Compte tenu de cette
identification, le plongement
$$
H^{0}(C_{I}\setminus\{c_{I}\},{\Bbb G}_{{\rm m}})/H^{0}(C_{I},{\Bbb
G}_{{\rm m}})\hookrightarrow E_{I}^{\times}/A_{I}^{\times}=G_{I}(k)
$$
d\'{e}finit un scindage $\sigma :\Lambda_{I}^{0}\rightarrow G_{I}$ de
l'extension $1\rightarrow G_{I}^{0}\rightarrow G_{I}\rightarrow
\Lambda_{I}\rightarrow 0$ au-dessus de $\Lambda_{I}^{0}\subset
\Lambda_{I}$.

On prolonge $G_{I}\twoheadrightarrow P_{I}$ en le morphisme de
$k$-sch\'{e}mas
$$
X_{I}\rightarrow\overline{P}_{I}
$$
qui envoie $M\subset E$ sur le ${\cal O}_{C_{I}}$-Module sans torsion
${\cal M}$ de rang g\'{e}n\'{e}rique $1$ obtenu en recollant ${\cal
O}_{C_{I}\setminus\{c_{I}\}}$ et $M$ le long de $\mathop{\rm Spec}
(E_{I})=(C_{I}\setminus\{c_{I}\})\times_{C_{I}}\mathop{\rm
Spec}(A_{I})$.  Pour chaque entier $d$, ce morphisme envoie la
composante connexe $X_{I}^{d}$ dans la composante connexe
$\overline{P}{}_{I}^{d}$.  Il est $G_{I}$-\'{e}quivariant pour
l'action naturelle de $G_{I}$ sur $X_{I}$ et l'action de $G_{I}$ sur
$P_{I}$ induite par celle de $P_{I}$ sur $\overline{P}_{I}$ et
l'\'{e}pimorphisme $G_{I}\twoheadrightarrow P_{I}$ ci-dessus, et il
passe au quotient en un morphisme de $k$-sch\'{e}mas projectifs
$$
Z_{I}=X_{I}/\sigma (\Lambda_{I}^{0})\rightarrow\overline{P}_{I}
$$
qui envoie, pour chaque entier $d$, la composante connexe $Z_{I}^{d}$
dans la composante connexe $\overline{P}{}_{I}^{d}$.  Le morphisme
$Z_{I}\rightarrow\overline{P}_{I}$ est birationnel puisqu'il induit un
isomorphisme de l'ouvert $Z_{I}^{\circ}=X_{I}^{\circ}/\Lambda_{I}^{0}
\subset Z_{I}$ sur l'ouvert $P_{I}\subset \overline{P}_{I}$.

\thm PROPOSITION 2.2.1
\enonce
Le $k$-morphisme birationnel $Z_{I}\rightarrow\overline{P}_{I}$
ci-dessus est un hom\'{e}omorphisme universel, c'est-\`{a}-dire est
fini, radiciel et surjectif.
\endthm

\rem Preuve
\endrem
Il suffit de voir que les fibres g\'{e}om\'{e}triques du morphisme
$Z_{I}\rightarrow\overline{P}_{I}$ sont toutes r\'{e}duites \`{a} un
point, avec \'{e}ventuellement des nilpotents.  Or tout ${\cal
O}_{C_{I}}$-Module ${\cal M}$ sans torsion de rang g\'{e}n\'{e}rique
$1$ s'obtient par recollement de ${\cal O}_{C_{I}\setminus\{c_{I}\}}$
et d'un $A_{I}$-r\'{e}seau $M\subset E_{I}$, le couple form\'{e} de
$M$ et de la donn\'{e}e de recollement \'{e}tant uniquement
d\'{e}termin\'{e} modulo l'action de $\sigma (\Lambda_{I}^{0})=
H^{0}(C_{I}\setminus\{c_{I}\},{\Bbb G}_{{\rm m}})/H^{0}(C_{I},{\Bbb
G}_{{\rm m}})$.
\hfill\hfill$\square$
\vskip 3mm

Compte tenu du th\'{e}or\`{e}me 2.1.2, on d\'{e}duit de cette
proposition:

\thm COROLLAIRE 2.2.2
\enonce
Le $k$-sch\'{e}ma $Z_{I}$ est irr\'{e}ductible et la $G_{I}$-orbite
$X_{I}^{\circ}$ de $M=A_{I}$ est dense dans $X_{I}$.
\hfill\hfill$\square$
\endthm

Bezrukavnikov [Be] a donn\'{e} une formule tr\`{e}s g\'{e}n\'{e}rale
pour la dimension des fibres de Springer, formule qui contient bien
entendu l'\'{e}galit\'{e} $\mathop{\rm dim}(X_{I})=\delta_{I}$. Par
contre, sa m\'{e}thode ne permet pas de d\'{e}montrer que $X_{I}$
n'admet pas de composantes irr\'{e}ductibles de dimension
$<\delta_{I}$.

\thm COROLLAIRE 2.2.3
\enonce
Le rev\^{e}tement \'{e}tale Galoisien $X_{I}\rightarrow Z_{I}$ de
groupe de Galois $\Lambda_{I}^{0}\cong\sigma (\Lambda_{I}^{0})$
provient par le changement de base $Z_{I}\rightarrow\overline{P}_{I}$
d'un rev\^{e}tement \'{e}tale Galoisien
$$
\overline{\varphi}_{I}:\overline{P}_{I}^{\,\natural}\rightarrow
\overline{P}_{I}
$$
dont la description au niveau des $k'$-points est la suivantes:
$\overline{P}_{I}^{\,\natural}(k')$ est l'ensemble des couples $({\cal
M},\iota )$, o\`{u} ${\cal M}$ est un ${\cal
O}_{k'\otimes_{k}C_{I}}$-Module sans torsion de rang g\'{e}n\'{e}rique
$1$ et $\iota :{\cal M}_{|k'\otimes_{k}(C_{I}\setminus\{c_{I}\})}
\buildrel\sim\over\longrightarrow {\cal
O}_{k'\otimes_{k}(C_{I}\setminus\{c_{I}\})}$ est une trivialisation de
la restriction de ${\cal M}$ \`{a}
$k'\otimes_{k}(C_{I}\setminus\{c_{I}\})$, et $\overline{\varphi}_{I}$
est le morphisme d'oubli de $\iota$.
\hfill\hfill$\square$
\endthm

On notera encore $P_{I}^{\natural}=G_{I}$ et
$\varphi_{I}:P_{I}^{\natural}\twoheadrightarrow P_{I}$ le
$k$-\'{e}pimorphisme d\'{e}fini plus haut.

\rem Exemples $2.2.4$
\endrem
Pour $|I|=1$ et $P_{I}(T)=T^{2}-\varpi_{F}^{3}$, $C_{I}$ est la
cubique n'ayant pour seule singularit\'{e} qu'un cusp ordinaire et
l'hom\'{e}omorphisme $Z_{I}^{0}\rightarrow\overline{P}_{I}^{0}$ est
naturellement isomorphe au morphisme de normalisation
$\pi_{I}:\widetilde{C}_{I}\rightarrow C_{I}$.

Pour $I=\{1,2\}$, $P_{1}(T)=T-\varpi_{F}$ et $P_{2}(T)=T+\varpi_{F}$,
$C_{I}$ est la cubique n'ayant pour seule singularit\'{e} qu'un point
double ordinaire, $\overline{\varphi}_{I}:Z_{I}\rightarrow
\overline{P}_{I}$ est un isomorphisme et la composante de degr\'{e}
$0$ du rev\^{e}tement $\overline{\varphi}_{I}:
\overline{P}_{I}^{\,\natural}\rightarrow\overline{P}_{I}$ est
naturellement isomorphe au rev\^{e}tement $C_{I}^{\natural}\rightarrow
C_{I}$ de groupe de Galois ${\Bbb Z}$ dont l'espace total est la
cha\^{\i}ne de droites projective index\'{e}e par ${\Bbb Z}$ obtenue
en prenant ${\Bbb Z}$ copies de la droite projective standard sur $k$
et en identifiant le point \`{a} l'infini de la $n$-\`{e}me copie
\`{a} l'origine de la $(n+1)$-\`{e}me.
\hfill\hfill$\square$
\vskip 3mm

On peut g\'{e}n\'{e}raliser la proposition $2.2.1$ comme suit.  Soit
$C$ une courbe g\'{e}om\'{e}tri\-que\-ment int\`{e}gre et projective sur
$k$ n'ayant que des singularit\'{e}s planes et rationnelles sur $k$,
et dont la source du morphisme de normalisation
$\pi_{C}:\widetilde{C}\rightarrow C$ est la droite projective standard
sur $k$.  Soit $\{c_{j}\mid j\in J\}$ l'ensemble fini des points
singuliers de $C$ et, pour chaque $j\in J$, soit $\pi_{C}^{-1}(c_{j})=
\{\widetilde{c}_{i}\mid i\in I_{j}\}$ l'ensemble des branches de $C$
en $c_{j}$.  Pour chaque $j\in J$, notons $A_{I_{j}}$ le
compl\'{e}t\'{e} de l'anneau local de $C$ en $c_{j}$ et $C_{I_{j}}$ la
courbe construite en $2.1$, de normalis\'{e}e une droite projective et
qui n'a qu'une seule singularit\'{e} (en $c_{I_{j}}$), singularit\'{e}
dont la germe formel est $\mathop{\rm Spf}(A_{I_{j}})$.

On peut former les $k$-sch\'{e}mas de Picard compactifi\'{e}s
$\overline{P}_{I_{j}}$ des $C_{I_{j}}$ et les {\og}{fibres de
Springer}{\fg} $X_{I_{j}}$ des $A_{I_{j}}$-r\'{e}seaux dans l'anneau
total des fractions $E_{I_{j}}$ de $A_{I_{j}}$, ainsi que leurs
quotients $Z_{I_{j}}=X_{I_{j}}/ \sigma_{j}(\Lambda_{I_{j}}^{0})$
o\`{u} $\Lambda_{I_{j}}^{0}=\mathop{\rm Ker}({\Bbb Z}^{I_{j}}\rightarrow
{\Bbb Z})\subset{\Bbb Z}^{I_{j}}=\Lambda_{I_{j}}$ et
$$
\sigma_{j}(\Lambda_{I_{j}}^{0})=H^{0}(C_{I_{j}}\setminus\{c_{I_{j}}\},{\Bbb
G}_{{\rm m}})/ H^{0}(C_{I_{j}},{\Bbb G}_{{\rm m}}) \subset
E_{I_{j}}^{\times}/A_{I_{j}}^{\times}.
$$

Alors, d'apr\`{e}s ce qui pr\'{e}c\`{e}de, on a un hom\'{e}omorphisme universel
$$
\prod_{j\in J}Z_{I_{j}}\rightarrow\prod_{j\in J}\overline{P}_{I_{j}}
$$
qui induit pour chaque famille d'entiers $(d_{j})_{j\in J}$ un
hom\'{e}omorphisme universel
$$
\prod_{j\in J}Z_{I_{j}}^{d_{j}}\rightarrow\prod_{j\in J}
\overline{P}_{I_{j}}^{d_{j}}.
$$

On peut aussi former le $k$-sch\'{e}ma de Picard compactifi\'{e}
$\overline{P}$ de $C$ et on a un $k$-morphisme
$$
X=\prod_{j\in J}X_{I_{j}}\rightarrow \overline{P}
$$
qui envoie $(M_{j}\subset E_{I_{j}})_{j\in J}$ sur le ${\cal
O}_{C}$-Module ${\cal M}$ sans torsion de rang g\'{e}n\'{e}rique $1$
obtenu en recollant ${\cal O}_{C\setminus\{c_{j}\mid j\in J\}}$ et
les $M_{j}$. Ce $k$-morphisme est \'{e}quivariant pour l'action du
groupe discret
$$
H^{0}(C\setminus\{c_{j}\mid j\in J\},{\Bbb G}_{{\rm m}})/
H^{0}(C,{\Bbb G}_{{\rm m}})\subset \prod_{j\in J}
E_{I_{j}}^{\times}/A_{I_{j}}^{\times}
$$
qui est l'image d'une section $\sigma$ de la projection canonique
$$
\prod_{j\in J}E_{I_{j}}^{\times}/A_{I_{j}}^{\times}\twoheadrightarrow
\prod_{j\in J}E_{I_{j}}^{\times}/{\cal O}_{E_{I_{j}}}^{\times}\cong
\prod_{j\in J}{\Bbb Z}^{I_{j}}={\Bbb Z}^{I}=\Lambda
$$
au-dessus de $\Lambda^{0}=\mathop{\rm Ker}({\Bbb Z}^{I}\rightarrow {\Bbb
Z})$ (on a bien entendu pos\'{e} $I=\coprod_{j\in J}I_{j}$). Il passe
donc au quotient en un $k$-morphisme
$$
Z=X/\sigma (\Lambda^{0})\rightarrow \overline{P}.
$$

La construction des courbes $C_{I_{j}}$, et donc des sections
$\sigma_{I_{j}}$, d\'{e}pend du choix d'uni\-for\-mi\-san\-tes
$\varpi_{i}=\varpi_{E_{i}}$ des $E_{i}$ pour les $i\in I_{j}$.  De
telles uniformisantes sont donn\'{e}es par le choix d'une carte affine
$\mathop{\rm Spec}(k[x])\subset\widetilde{C}\cong {\Bbb P}_{k}^{1}$ qui
contient tous les points $\widetilde{c}_{i}$, $i\in I$.  Si on fixe de
cette fa\c{c}on les uniformisantes, la restriction de $\sigma
:\Lambda^{0}\rightarrow \prod_{j\in J}
E_{I_{j}}^{\times}/A_{I_{j}}^{\times}$ \`{a}
$$
\prod_{j\in J}\Lambda_{I_{j}}^{0}=\prod_{j\in J}\mathop{\rm Ker}({\Bbb
Z}^{I_{j}}\rightarrow {\Bbb Z})\subset\mathop{\rm Ker}({\Bbb
Z}^{I}\rightarrow {\Bbb Z})=\Lambda^{0}
$$
est \'{e}gale \`{a} la section $\prod_{j\in J}\sigma_{I_{j}}$ et
l'application quotient $X=\prod_{j\in J}X_{I_{j}}\twoheadrightarrow Z$
se factorise par un $k$-morphisme
$$
\prod_{j\in J}Z_{I_{j}}\twoheadrightarrow Z
$$
qui induit, pour chaque famille d'entiers $(d_{j})_{j\in
J}$ de somme $d$, un $k$-isomorphisme
$$
\prod_{j\in J}Z_{I_{j}}^{d_{j}}\buildrel\sim\over\longrightarrow
Z^{d}
$$
o\`{u} $Z^{d}$ est la composante connexe de degr\'{e} $d$ de $Z$.

\thm PROPOSITION 2.2.5
\enonce
Le $k$-morphisme ci-dessus $Z\rightarrow \overline{P}$, somme
disjointe pour $d\in {\Bbb Z}$ des $k$-morphisme $Z^{d}\rightarrow
\overline{P}^{d}$, est un hom\'{e}omorphisme universel.
\hfill\hfill$\square$
\endthm

\thm COROLLAIRE 2.2.6
\enonce
Le rev\^{e}tement \'{e}tale Galoisien $X\rightarrow Z$ de
groupe de Galois $\Lambda^{0}\cong\sigma (\Lambda^{0})$
provient par le changement de base $Z\rightarrow\overline{P}$
d'un rev\^{e}tement \'{e}tale Galoisien
$$
\overline{\varphi}:\overline{P}^{\,\natural}\rightarrow
\overline{P}
$$
dont la description au niveau des $k'$-points est la suivantes:
$\overline{P}^{\,\natural}(k')$ est l'ensemble des couples $({\cal
M},\iota )$, o\`{u} ${\cal M}$ est un ${\cal
O}_{k'\otimes_{k}C}$-Module sans torsion de rang g\'{e}n\'{e}rique $1$
et $\iota :{\cal M}_{|k'\otimes_{k}(C\setminus\{c_{j}\mid j\in
J\})}\buildrel\sim\over\longrightarrow {\cal
O}_{k'\otimes_{k}(C\setminus\{c_{j}\mid j\in J\})}$ est une
trivialisation de la restriction de ${\cal M}$ \`{a}
$k'\otimes_{k}(C\setminus\{c_{j}\mid j\in J\})$, et
$\overline{\varphi}$ est le morphisme d'oubli de $\iota$.
\hfill\hfill$\square$
\endthm

\section{2.3}{Retour \`{a} la situation arithm\'{e}tique du lemme
fondamental}

Pla\c{c}ons nous de nouveau dans la situation de la section 1.3, de
sorte que $k={\Bbb F}_{q^{2}}$.  On suppose de plus que $|I|\leq q$.

Compte tenu de l'hypoth\`{e}se (1.3.1), l'automorphisme $\ast$ de
${\cal O}_{E_{I}}$ envoie $A_{I}$ dans lui-m\^{e}me. On v\'{e}rifie
alors que
$$\eqalign{
A_{I,0}&=A_{I}\cap E_{I,0}=\{x\in A_{I}\mid x=x^{\ast}\}\cr
&=\{\alpha_{0}+\alpha_{1}\gamma_{I}+\cdots
+\alpha_{n_{I}-1}\gamma_{I}^{n_{I}-1}\mid \alpha_{j}\in {\cal O}_{F}\hbox{
et }\alpha_{j}^{\ast}=(-1)^{j}\alpha_{j}\}\cr}
$$
est une ${\Bbb F}_{q}$-structure sur la ${\Bbb F}_{q^{2}}$-alg\`{e}bre
$A_{I}$ et que la normalis\'{e}e de la ${\Bbb F}_{q}$-alg\`{e}bre
$A_{I,0}$ est $E_{I,0}$.  Le morphisme de normalisation
$$
\mathop{\rm Spf}({\cal O}_{E_{I}})\rightarrow\mathop{\rm Spf}(A_{I})
$$
se descend donc de ${\Bbb F}_{q^{2}}$ \`{a} ${\Bbb F}_{q}$ en le
morphisme de normalisation
$$
\mathop{\rm Spf}({\cal O}_{E_{I,0}})\rightarrow\mathop{\rm
Spf}(A_{I,0}).
$$
Bien s\^{u}r, on peut aussi introduire de la m\^{e}me fa\c{c}on une
${\Bbb F}_{q}$-structure $A_{i,0}$ sur la ${\Bbb
F}_{q^{2}}$-alg\`{e}bre $A_{i}$ pour chaque $i\in I$.

On peut reprendre la construction de $C_{I}$ donn\'{e}e en $2.1$ en
choisissant $\iota :I\hookrightarrow {\Bbb F}_{q^{2}}=k$ de telles
sorte que son image soit contenue dans ${\Bbb F}_{q}\subset {\Bbb
F}_{q^{2}}$ et en choisissant chaque uniformisante $\varpi_{E_{i}}$
dans $E_{i,0}\subset E_{i}$.  La courbe $C_{I}$ obtenue se descend
alors en une courbe projective et g\'{e}om\'{e}triquement int\`{e}gre
$C_{I,0}$ sur ${\Bbb F}_{q}$, n'ayant qu'un seul point singulier
rationnel sur ${\Bbb F}_{q}$, dont le germe formel en ce point
singulier est $\mathop{\rm Spf}(A_{I,0})$ et a des branches
$\mathop{\rm Spf}(A_{i,0})$ qui sont toutes rationnelles sur ${\Bbb
F}_{q}$, et dont la normalis\'{e}e est ${\Bbb P}_{{\Bbb F}_{q}}^{1}$.

On munit alors le ${\Bbb F}_{q^{2}}$-sch\'{e}ma de Picard
compactifi\'{e} $\overline{P}_{I}$ de la ${\Bbb F}_{q}$-structure
rationnelle $\overline{{\cal P}}_{I}$ dont l'endomorphisme de
Frobenius envoie, pour toute extension $k'$ de ${\Bbb F}_{q^{2}}$, un
$k'$-point ${\cal M}$ sur le $k'\otimes_{k}C_{I}$-Module sans torsion
de rang g\'{e}n\'{e}rique $1$,
$$
\mathop{\rm Frob}\nolimits_{\overline{{\cal P}}_{I}/{\Bbb
F}_{q}}({\cal M})=(k'\otimes_{k}\mathop{\rm
Frob}\nolimits_{C_{I,0}/{\Bbb F}_{q}})^{\ast}({\cal
M}^{\vee})=((k'\otimes_{k}\mathop{\rm Frob}\nolimits_{C_{I,0}/{\Bbb
F}_{q}})^{\ast}{\cal M})^{\vee},
$$
o\`{u} $\mathop{\rm Frob}\nolimits_{C_{I,0}/{\Bbb F}_{q}}:C_{I}
\rightarrow C_{I}$ est l'endomorphisme de Frobenius pour la ${\Bbb
F}_{q}$-structure rationnelle $C_{I,0}$ sur $C_{I}$ et o\`{u} ${\cal
M}^{\vee}=\mathop{{\cal H}{\it om}}\nolimits_{k'\otimes_{k}{\cal
O}_{C_{I}}}({\cal M},k'\otimes_{k}{\cal O}_{C_{I}})$ est le dual de
${\cal M}$.  Si ${\cal M}$ est inversible, il en est de m\^{e}me de
$\mathop{\rm Frob}\nolimits_{\overline{{\cal P}}_{I}/{\Bbb
F}_{q}}({\cal M})$, de sorte que $\mathop{\rm
Frob}\nolimits_{\overline{{\cal P}}_{I}/{\Bbb F}_{q}}$ induit un
endomorphisme de Frobenius sur l'ouvert $P_{I}\subset
\overline{P}_{I}$, d'o\`{u} une ${\Bbb F}_{q}$-structure rationnelle
$\overline{{\cal P}}_{I}$ sur $P_{I}$.

\thm LEMME 2.3.1
\enonce
Les morphismes $G_{I}\rightarrow P_{I}$ et $X_{I}\rightarrow
\overline{P}_{I}$ sont compatibles aux ${\Bbb F}_{q}$-structures
rationnelles consid\'{e}r\'{e}es en $1.3$ et ci-dessus; ils se descendent
donc en des ${\Bbb F}_{q}$-morphismes ${\cal G}_{I}\rightarrow {\cal
P}_{I}$ et ${\cal X}_{I}\rightarrow \overline{{\cal P}}_{I}$.
\endthm

\rem Preuve
\endrem
Il suffit de d\'{e}montrer que, pour toute extension finie $k'$ de $k$
et tout sous-$A_{I}'$-module $M\subset E_{I}'(=k'\otimes_{k}E_{I})$
partout de rang g\'{e}n\'{e}rique $1$, on a
$$
(M^{\vee})^{\ast}=(M^{\ast})^{\vee}=M^{\perp_{I}}
$$
o\`{u}
$$
M^{\vee}=\mathop{\rm Hom}\nolimits_{A_{I}'}(M,A_{I}')=\{y\in
E_{I}'\mid yM\subset A_{I}'\}.
$$

Or, on a les inclusions \'{e}videntes
$$
(M^{\vee})^{\ast}\subset (M^{\ast})^{\vee}\subset M^{\perp_{I}},
$$
inclusions qui sont toutes les deux des \'{e}galit\'{e}s pour
$M=A_{I}'$, et on a
$$
[M^{\vee}:N^{\vee}]=[M^{\ast}:N^{\ast}]
=[M^{\perp_{I}}:N^{\perp_{I}}]=-[M:N]
$$
pour tout couple $(M,N)$ de sous-$A_{I}'$-modules de $E_{I}'$ partout
de rang g\'{e}n\'{e}rique $1$.
\hfill\hfill$\square$
\vskip 3mm

*

\section{2.4}{Auto-dualit\'{e} des jacobiennes compactifi\'{e}es
d'apr\`{e}s Esteves, Gagn\'{e} et Kleiman}

Dans [E-G-K], Esteves, Gagn\'{e} et Kleiman ont d\'{e}montr\'{e} un
th\'{e}or\`{e}me d'auto-dualit\'{e} pour les jacobiennes
compactifi\'{e}es des courbes projectives et int\`{e}gres dont toutes
les singularit\'{e}s sont planes et de multiplicit\'{e} $2$.  Ce
th\'{e}or\`{e}me g\'{e}n\'{e}ralise l'\'{e}nonc\'{e} classique
d'auto-dualit\'{e} des jacobiennes des courbes projectives et
lisses.

Nous n'aurons besoin que de la partie {\og}{facile}{\fg} de ce
th\'{e}or\`{e}me, partie qui vaut en fait sans l'hypoth\`{e}se
restrictive de multiplicit\'{e} $2$ et que nous allons rappeler
maintenant.
\vskip 2mm

On se place dans la situation naturelle pour ce r\'{e}sultat.  Soient
donc $S$ un sch\'{e}ma n{\oe}th\'{e}rien et $C\rightarrow S$ un
morphisme projectif et plat, dont toutes les fibres
g\'{e}om\'{e}triques sont int\`{e}gres et de dimension $1$.  Pour
simplifier, on supposera qu'il existe une section globale de $C$ sur
$S$ dont l'image est contenue dans le lieu de lissit\'{e} de $C$ sur
$S$ et on fixera une fois pour toute une telle section $\infty
:S\rightarrow C$.  On notera $[\infty]$ le diviseur de Cartier relatif
sur $C$ image de cette section.

Pour chaque entier $d$, soit $P^{d}=\mathop{\rm Pic}
\nolimits_{C/S}^{d}$ la composante du $S$-sch\'{e}ma de Picard de la
courbe $C/S$ qui param\`{e}tre les classes d'isomorphie de Modules
inversibles sur $C$ qui sont de degr\'{e} $d$ fibre \`{a} fibre du
morphisme $C\rightarrow S$, et soit $\overline{P}{}^{d}=
\mathop{\overline{\rm Pic}}{}_{\!C/S}^{d}$ la compactification
relative de $P^{d}$ qui param\`{e}tre les Modules coh\'{e}rents sur
$C$ qui sont plats sur $S$ et, fibre \`{a} fibre, sans torsion, de rang
g\'{e}n\'{e}rique $1$ et de degr\'{e} $d$.  La torsion ${\cal
M}\rightarrow {\cal M}(-(d+1)[\infty ])$ par le diviseur de Cartier
$-(d+1)[\infty ]$ identifie $P^{d}$ et $\overline{P}{}^{d}$ \`{a}
$P^{-1}$ et $\overline{P}{}^{-1}$ respectivement.

On supposera que, pour un entier $d$ ou, ce qui revient au m\^{e}me,
pour tout entier $d$, le foncteur de Picard de $\overline{P}{}^{d}/S$ est
repr\'{e}sentable par un $S$-sch\'{e}ma
$$
\mathop{\rm Pic}\nolimits_{\overline{P}{}^{d}/S}
$$
qui est une r\'{e}union disjointe de $S$-sch\'{e}mas quasi-projectifs.
C'est le cas d'apr\`{e}s Grothendieck ([Gr 2] Th\'{e}or\`{e}me 3.1, et
aussi [B-L-R] \S 8.2, Theorem 1) si toutes les fibres de $C/S$ ont au
pire des singularit\'{e}s planes, car le $S$-sch\'{e}ma
$\overline{P}{}^{d}$ est alors projectif, plat et \`{a} fibres
g\'{e}om\'{e}triques int\`{e}gres d'apr\`{e}s Altman, Iarrobino et
Kleiman, et Rego (cf.  [A-I-K], [Re] et notre section 2.1).  C'est
aussi le cas si $S$ est le spectre d'un corps d'apr\`{e}s Murre et
Oort (cf.  [B-L-R] \S8.2, Theorem 3).

Si $S$ est le spectre d'un corps $k$, le $k$-sch\'{e}ma en groupes
$\mathop{\rm Pic}\nolimits_{\overline{P}{}^{d}/k}$ admet une
composante neutre $\mathop{\rm Pic}
\nolimits_{\overline{P}{}^{d}/k}^{0}$ et on d\'{e}finit
$$
\mathop{\rm Pic}\nolimits_{\overline{P}{}^{d}/k}^{\tau}=\bigcup_{n>0}
[n]^{-1} \mathop{\rm Pic}\nolimits_{\overline{P}{}^{d}/k}^{0},
$$
o\`{u} $[n]:\mathop{\rm Pic}\nolimits_{\overline{P}{}^{d}/k}
\rightarrow\mathop{\rm Pic}\nolimits_{\overline{P}{}^{d}/k}$ est la
multiplication par l'entier $n$.  Pour $S$ arbitraire, on note
$\mathop{\rm Pic}\nolimits_{\overline{P}{}^{d}/S}^{0}$ (resp.
$\mathop{\rm Pic}\nolimits_{\overline{P}{}^{d}/S}^{\tau}$) le
sous-foncteur de $\mathop{\rm Pic}\nolimits_{\overline{P}{}^{d}/S}$
form\'{e} des classes d'isomorphie de Modules inversibles sur
$\overline{P}{}^{d}$ dont la restriction \`{a} chaque fibre
$\overline{P}{}_{s}^{d}$ de $\overline{P}{}^{d}\rightarrow S$ est dans
$\mathop{\rm Pic}\nolimits_{\overline{P}{}_{s}^{d}/\kappa (s)}^{0}$
(resp.  $\mathop{\rm Pic}\nolimits_{\overline{P}{}_{s}^{d}/\kappa
(s)}^{\tau}$).  Le sous-foncteur $\mathop{\rm Pic}
\nolimits_{\overline{P}{}^{d}/S}^{\tau}$ est repr\'{e}sentable par une
partie ouverte et ferm\'{e}e de $\mathop{\rm Pic}
\nolimits_{\overline{P}{}^{d}/S}$ (cf.  [B-L-R] \S 8.4, Theorem 4);
par contre, le sous-foncteur $\mathop{\rm Pic}
\nolimits_{\overline{P}{}^{d}/S}^{0}$ n'est pas repr\'{e}sentable en
g\'{e}n\'{e}ral.

On a l'application d'Abel
$$
A_{-1}:C\rightarrow\overline{P}{}^{-1}
$$
d\'{e}finie par l'Id\'{e}al de la diagonale $C\subset C\times_{S}C$:
la premi\`{e}re projection $C\times_{S}C\rightarrow C$ est un
changement de base de $C\rightarrow S$ et cet Id\'{e}al est un ${\cal
O}_{C\times_{S}C}$-Module coh\'{e}rent qui est plat sur $C$ et, fibre
\`{a} fibre, sans torsion, de rang g\'{e}n\'{e}rique $1$ et de
degr\'{e} $-1$.  Pour tout entier $d$, on d\'{e}finit
$$
A_{d}:C\rightarrow\overline{P}{}^{d}
$$
comme le compos\'{e} de $A_{-1}$ et de l'isomorphisme
$\overline{P}{}^{-1}\buildrel\sim\over\longrightarrow
\overline{P}{}^{d}$ de torsion par $(d+1)[\infty ]$.  Le morphisme
$A_{d}$ induit un homomorphisme de $S$-sch\'{e}mas en groupes
$$
A_{d}^{\ast}:\mathop{\rm Pic}\nolimits_{\overline{P}{}^{d}/S}
\rightarrow\mathop{\rm Pic}\nolimits_{C/S}.
$$

Dans [E-G-K], Esteves, Gagn\'{e} et Kleiman construisent un inverse
\`{a} droite de l'homomorphisme $A_{d}^{\ast}$ sur $P^{0}=\mathop{\rm
Pic}\nolimits_{C/S}^{0}\subset\mathop{\rm Pic}\nolimits_{C/S}$ en
utilisant le d\'{e}terminant de la cohomologie.

Plus pr\'{e}cis\'{e}ment, notons
$$\diagram{
C\times_{S}\overline{P}&\kern -1mm\smash{\mathop{\hbox to
10mm{\leftarrowfill}} \limits^{\scriptstyle \mathop{\rm pr}
\nolimits_{12}}}\kern -1mm&
C\times_{S}\overline{P}\times_{S}P^{0} &\kern -1mm
\smash{\mathop{\hbox to 10mm{\rightarrowfill}} \limits^{\scriptstyle
\mathop{\rm pr}\nolimits_{13}}}\kern -1mm& C\times_{S}P^{0}\cr
\llap{$\scriptstyle \mathop{\rm pr}\nolimits_{2}$}\left\downarrow
\vbox to 4mm{}\right.\rlap{}&&\llap{}\left\downarrow \vbox to
4mm{}\right.\rlap{$\scriptstyle \mathop{\rm
pr}\nolimits_{23}$}&&\llap{}\left\downarrow \vbox to
4mm{}\right.\rlap{$\scriptstyle \mathop{\rm pr}\nolimits_{2}$}\cr
\overline{P}&\kern -6mm\smash{\mathop{\hbox to
20mm{\leftarrowfill}} \limits_{\scriptstyle \mathop{\rm pr}
\nolimits_{1}}}\kern -6mm& \overline{P}{}\times_{S}P^{0}&\kern
-6mm\smash{\mathop{\hbox to 20mm{\rightarrowfill}}
\limits_{\scriptstyle \mathop{\rm pr}\nolimits_{2}}}\kern -6mm&
P^{0}\cr}
$$
les projections canoniques et ${\cal L}^{{\rm univ}}$ et ${\cal
M}^{{\rm univ}}$ les Modules universels sur $C\times_{S}P^{0}$ et
$C\times_{S}\overline{P}$ respectivement, rigidifi\'{e}s le long des
sections $P^{0}\rightarrow C\times_{S}P^{0}$ et $\overline{P}
\rightarrow C\times_{S}\overline{P}$ induites par la section $\infty
:S\rightarrow C$.  Alors, on peut former le ${\cal O}_{\overline{P}
\times_{S}P^{0}}$-Module inversible
$$
(\det R\mathop{\rm pr}\nolimits_{23,\ast}(\mathop{\rm
pr}\nolimits_{12}^{\ast}{\cal M}^{{\rm univ}}\otimes\mathop{\rm
pr}\nolimits_{13}^{\ast}{\cal L}^{{\rm univ}}))^{\otimes -1}
\otimes\det R\mathop{\rm pr} \nolimits_{23,\ast}\mathop{\rm
pr}\nolimits_{12}^{\ast}{\cal M}^{{\rm univ}}
$$
sur $\overline{P}\times_{S}P^{0}$, Module inversible qui d\'{e}finit
un morphisme
$$
\beta=\prod_{d\in {\Bbb Z}}\beta_{d}:P^{0}\rightarrow \mathop{\rm Pic}
\nolimits_{\overline{P}/S}=\prod_{d\in {\Bbb Z}}\mathop{\rm Pic}
\nolimits_{\overline{P}{}^{d}/S}.
$$

\thm PROPOSITION 2.4.1 (Esteves, Gagn\'{e} et Kleiman; [E-G-K],
Proposition (2.2))
\enonce
Pour chaque entier $d$, le morphisme $\beta_{d}:P^{0}
\rightarrow\mathop{\rm Pic}\nolimits_{\overline{P}{}^{d}/S}$ est un
homomorphisme de $S$-sch\'{e}mas en groupes dont l'image ensembliste
est contenue dans le sous-foncteur $\mathop{\rm
Pic}\nolimits_{\overline{P}{}^{d}/S}^{0}$, et donc dans l'ouvert et
ferm\'{e} $\mathop{\rm Pic}\nolimits_{\overline{P} {}^{d}/S}^{\tau}$,
et la formation de $\beta_{d}$ commute \`{a} tout changement de base
$S'\rightarrow S$.  De plus, le compos\'{e}
$A_{d}^{\ast}\circ\beta_{d}$ est l'identit\'{e} de $P^{0}$.
\hfill\hfill$\square$
\endthm

\rem Remarque $2.4.2$
\endrem
Si on note
$$
\mu :\overline{P}\times_{S}P^{0}\rightarrow\overline{P},~({\cal M},
{\cal L})\mapsto {\cal M}\otimes_{{\cal O}_{C}}{\cal L},
$$
l'action naturelle de $P^{0}$ sur $\overline{P}$, on a le
carr\'{e} cart\'{e}sien
$$\diagram{
C\times_{S}\overline{P}\times_{S}P^{0}&\kern -1mm\smash{\mathop{\hbox
to 16mm{\rightarrowfill}} \limits^{\scriptstyle \mathop{\rm
Id}\nolimits_{C}\times \mu}}\kern -1mm&C\times_{S}\overline{P}\cr
\llap{$\scriptstyle \mathop{\rm pr}\nolimits_{23}$}\left\downarrow
\vbox to 4mm{}\right.\rlap{}&\square&\llap{}\left\downarrow
\vbox to 4mm{}\right.\rlap{$\scriptstyle \mathop{\rm pr}\nolimits_{2}$}\cr
\overline{P}\times_{S}P^{0}&\kern -8mm\smash{\mathop{\hbox to
26mm{\rightarrowfill}} \limits_{\scriptstyle \mu}}\kern
-8mm&\overline{P}\cr}
$$
et on a un isomorphisme canonique
$$
\mathop{\rm pr}\nolimits_{12}^{\ast}{\cal M}^{{\rm univ}}
\otimes\mathop{\rm pr}\nolimits_{13}^{\ast}{\cal L}^{{\rm univ}}\cong
(\mathop{\rm Id}\nolimits_{C}\times\mu)^{\ast}{\cal M}^{{\rm univ}}.
$$
Le th\'{e}or\`{e}me de changement de base assure alors que
$$
R\mathop{\rm pr}\nolimits_{23,\ast}(\mathop{\rm pr}
\nolimits_{12}^{\ast}{\cal M}^{{\rm univ}}\otimes\mathop{\rm pr}
\nolimits_{13}^{\ast}{\cal L}^{{\rm univ}})\cong\mu^{\ast}
R\mathop{\rm pr}\nolimits_{2,\ast}{\cal M}^{{\rm univ}}
$$
et le morphisme $\beta$ d'Esteves, Gagn\'{e} et Kleiman est encore
d\'{e}fini par le Module inversible
$$
(\mu^{\ast}\det R\mathop{\rm pr}\nolimits_{2,\ast}{\cal M}^{{\rm
univ}})^{\otimes -1}\otimes \det R\mathop{\rm
pr}\nolimits_{2,\ast}{\cal M}^{{\rm univ}}.\eqno{\square}
$$

\section{2.5}{Sch\'{e}mas de Picard et dualit\'{e} pour les tores}

Nous aurons besoin dans la section suivante et la section $4.3$ d'un
r\'{e}sultat g\'{e}n\'{e}ral de Grothendieck qui nous a \'{e}t\'{e}
communiqu\'{e} par Raynaud.
\vskip 2mm

Soient $f:Z\rightarrow S$ un morphisme propre, plat et de
pr\'{e}sentation finie de sch\'{e}mas tel que $f_{\ast}{\cal
O}_{Z}={\cal O}_{S}$ et $T$ un $S$-tore (plat et de pr\'{e}sentation
finie), de faisceau des caract\`{e}res $X^{\ast}(T)= \mathop{{\cal
H}{\it om}}\nolimits_{S\hbox{\sevenrm -sch.gr.}}(T,{\Bbb G}_{{\rm
m},S})$.  Pour chaque changement de base
$$\diagram{
Z'&\kern -1mm\smash{\mathop{\hbox to 8mm{\rightarrowfill}}
\limits^{\scriptstyle }}\kern -1mm&Z\cr
\llap{$\scriptstyle f'$}\left\downarrow
\vbox to 4mm{}\right.\rlap{}&\square&\llap{}\left\downarrow
\vbox to 4mm{}\right.\rlap{$\scriptstyle f$}\cr
S'&\kern -1mm\smash{\mathop{\hbox to 8mm{\rightarrowfill}}
\limits_{\scriptstyle }}\kern -1mm&S&,\cr}\kern 5mm
$$
chaque section globale $t'$ du $S'$-tore $T'=S'\times_{S}T$ et chaque
$f'^{\ast}X^{\ast}(T')$-torseur $Y'$ sur $Z'$, on note $\langle
Y',t'\rangle$ le ${\Bbb G}_{{\rm m},Z'}$-torseur sur $Z'$ obtenu en
poussant $Y'$ par le caract\`{e}re image r\'{e}ciproque par $f'$ du
caract\`{e}re $X^{\ast}(T')\rightarrow {\Bbb G}_{{\rm
m},S'},~\chi'\mapsto \chi'(t')$.

\thm PROPOSITION 2.5.1 (cf. [Ra] Proposition (6.2.1))
\enonce
Notons
$$
\mathop{\rm Pic}\nolimits_{Z/S}=R^{1}f_{\ast}{\Bbb G}_{{\rm
m},Z}
$$
le foncteur de Picard relatif. Alors, l'homomorphisme de faisceaux
\'{e}tales sur $S$
$$
R^{1}f_{\ast}f^{\ast}X^{\ast}(T)\rightarrow\mathop{{\cal H}{\it om}}
\nolimits_{S\hbox{\sevenrm -sch.gr.}}(T,\mathop{\rm Pic}
\nolimits_{Z/S})
$$
qui, quel que soit le $S$-sch\'{e}ma $S'$, envoie la classe d'un
$f'^{\ast}X^{\ast}(T')$-torseur $Y'\rightarrow Z'$ sur l'homomorphisme
$\langle Y', \cdot\rangle :T' \rightarrow\mathop{\rm
Pic}\nolimits_{Z'/S'}$, est un isomorphisme.
\endthm

\rem Preuve
\endrem
On raisonne comme le fait Raynaud pour prouver la Proposition (6.2.1)
de [Ra].  Pour tout faisceau fppf en groupes commutatifs ${\cal F}$
sur $S$, on consid\`{e}re le foncteur
$$
{\cal G}\mapsto H({\cal G})=f_{\ast}\mathop{{\cal H}{\it
om}}(f^{\ast}{\cal F},{\cal G})
$$
sur la cat\'{e}gorie des faisceaux fppf en groupes commutatifs sur
$Z$. On a deux suite spectrale
$$
E_{2}^{pq}=R^{q}f_{\ast}\mathop{{\cal E}{\it
xt}}\nolimits^{p}(f^{\ast}{\cal F},{\cal G})\Rightarrow R^{p+q}H({\cal
G})
$$
et
$$
E_{2}^{pq}=\mathop{{\cal E}{\it xt}}\nolimits^{q}({\cal
F},R^{p}f_{\ast}{\cal G})\Rightarrow R^{p+q}H({\cal G})
$$
et donc deux suites exactes courtes des termes de bas degr\'{e}s, qui 
s'\'{e}crivent
pour ${\cal F}=T$ et ${\cal G}={\Bbb G}_{{\rm m},Z}$,
$$
0\rightarrow R^{1}f_{\ast}f^{\ast}X^{\ast}(T)\rightarrow R^{1}H({\Bbb G}_{{\rm
m},Z})\rightarrow f_{\ast}\mathop{{\cal E}{\it
xt}}\nolimits^{1}(f^{\ast}T,{\Bbb G}_{{\rm m},Z})
$$
et
$$\eqalign{
0\rightarrow &\mathop{{\cal E}{\it xt}}\nolimits^{1} (T,f_{\ast}{\Bbb
G}_{{\rm m},Z})\rightarrow R^{1}H({\Bbb G}_{{\rm m},Z})\rightarrow
\mathop{{\cal H}{\it om}}(T,R^{1}f_{\ast}{\Bbb G}_{{\rm m},Z})\cr
\rightarrow &\mathop{{\cal E}{\it xt}}\nolimits^{2} (T,f_{\ast}{\Bbb
G}_{{\rm m},Z})\rightarrow R^{2}H({\Bbb G}_{{\rm m},Z}).\cr}
$$
L'homomorphisme de la proposition est alors l'homomorphisme compos\'{e}
$$
R^{1}f_{\ast}f^{\ast}X^{\ast}(T)\rightarrow R^{1}H({\Bbb G}_{{\rm
m},Z})\rightarrow \mathop{{\cal H}{\it
om}}(T,R^{1}f_{\ast}{\Bbb G}_{{\rm m},Z}).
$$

On a $f_{\ast}{\Bbb G}_{{\rm m},Z}={\Bbb G}_{{\rm m},S}$ par
hypoth\`{e}se et on sait que $\mathop{{\cal E}{\it xt}}\nolimits^{1}
(f^{\ast}T,{\Bbb G}_{{\rm m},Z})=(0)$ et $\mathop{{\cal E}{\it
xt}}\nolimits^{1}(T,{\Bbb G}_{{\rm m},S})=(0)$.  Par suite,
l'homomorphisme compos\'{e} ci-dessus est injectif et, pour
d\'{e}montrer qu'il est bijectif, il suffit de v\'{e}rifier que
l'application $\mathop{{\cal E}{\it xt}}\nolimits^{2}(T,f_{\ast}{\Bbb
G}_{{\rm m},Z})\rightarrow R^{2}H({\Bbb G}_{{\rm m},Z})$ est
injective.  Comme cette propri\'{e}t\'{e} d'injectivit\'{e} est locale
pour la topologie fppf sur $S$, on peut supposer que $f$ admet une
section globale $z:S\rightarrow Z$, et donc que l'on a une fl\`{e}che $\iota
:f_{\ast}{\cal G}\rightarrow z^{\ast}{\cal G}$ fonctorielle en ${\cal
G}$.  Or, la fl\`{e}che $\iota$ induit un morphisme $R^{2}H({\cal
G})\rightarrow \mathop{{\cal E}{\it xt}}\nolimits^{2}({\cal
F},z^{\ast}{\cal G})$ tel que l'homomorphisme compos\'{e}
$\mathop{{\cal E}{\it xt}}\nolimits^{2}({\cal F},f_{\ast}{\cal
G})\rightarrow R^{2}H({\cal G})\rightarrow \mathop{{\cal E}{\it
xt}}\nolimits^{2}({\cal F},z^{\ast}{\cal G})$ soit \'{e}gal \`{a}
$\mathop{{\cal E}{\it xt}}\nolimits^{2}({\cal F},\iota )$, et elle est
un isomorphisme pour ${\cal G}={\Bbb G}_{{\rm m},Z}$ par
hypoth\`{e}se.  Notre assertion d'injectivit\'{e} est donc
v\'{e}rifi\'{e}e et la proposition est d\'{e}montr\'{e}e.
\hfill\hfill$\square$

\section{2.6}{Auto-dualit\'{e} des jacobiennes compactifi\'{e}es et
fibres de Springer}

Revenons maintenant \`{a} la situation de la section 2.1.  Nous avons
donc la courbe int\`{e}gre et projective $C=C_{I}$ sur $k$, de
normalis\'{e}e $\pi =\pi_{I}:\widetilde{C}=\widetilde{C}_{I}
\rightarrow C$ une droite projective, avec son unique point singulier
$c=c_{I}$ en lequel la singularit\'{e} est plane.  Le compl\'{e}t\'{e}
de l'anneau local de $C$ en $c$ est notre anneau $A=A_{I}$ de la
section 1.1 et l'ensemble des branches $\pi^{-1}(c)=
\{\widetilde{c}_{i}\mid i\in I\}$ est index\'{e} par $I$.

On consid\`{e}re le $k$-sch\'{e}ma de Picard compactifi\'{e}
$\overline{P}=\overline{P}_{I}$ de $C$ et son rev\^{e}tement \'{e}tale
Galoisien
$$
\overline{P}^{\,\natural}\rightarrow\overline{P},
$$
de groupe de Galois
$$\eqalign{
\Lambda^{0}&=H^{0}(C\setminus\{c\},{\Bbb G}_{{\rm m}})/
H^{0}(C,{\Bbb G}_{{\rm m}})\cr
&= H^{0}(\widetilde{C}\setminus \{\widetilde{c}_{i}\mid i\in I\},
{\Bbb G}_{{\rm m}})/ H^{0}(\widetilde{C},{\Bbb G}_{{\rm m}})
=\mathop{\rm Ker}({\Bbb Z}^{I} \rightarrow {\Bbb Z}),\cr}
$$
d\'{e}fini dans la section $2.2$.  Rappelons que les $k$-points de
$\overline{P}^{\,\natural}$ sont les couples $({\cal M},\iota )$
o\`{u} ${\cal M}$ est un ${\cal O}_{C}$-Module sans torsion de rang
g\'{e}n\'{e}rique $1$ et $\iota :{\cal M}_{|C\setminus\{c\}}
\buildrel\sim\over\longrightarrow {\cal O}_{C\setminus\{c\}}$ est une
rigidification de ${\cal M}$ sur le lieu de lissit\'{e}
$C\setminus\{c\}$ de $C$, et que la fl\`{e}che du rev\^{e}tement est
l'oubli de la rigidification.

Par construction, la restriction du rev\^{e}tement
$\overline{P}^{\,\natural}\rightarrow\overline{P}$ par le morphisme
radiciel
$$
Z=Z_{I}\rightarrow\overline{P}
$$
d\'{e}fini en $2.2$ est le rev\^{e}tement
$$
X\rightarrow X/\Lambda^{0}=Z
$$
o\`{u} $X=X_{I}$ est la fibre de Springer de la premi\`{e}re partie.

Nous allons utiliser les r\'{e}sultats g\'{e}n\'{e}raux des sections
2.4 et 2.5 pour donner une autre d\'{e}finition du rev\^{e}tement
$\overline{P}^{\,\natural}\rightarrow\overline{P}$.
\vskip 2mm

Soit $T$ le tore maximal de la composante neutre $P^{0}$ du
$k$-sch\'{e}ma de Picard de $C$.  On a $T={\Bbb G}_{{\rm m},k}^{I}/
{\Bbb G}_{{\rm m},k}$ et le groupe des caract\`{e}res de $T$ est
naturellement isomorphe au noyau
$$
\mathop{\rm Ker}({\Bbb Z}^{I}\rightarrow {\Bbb Z})
$$
de l'homomorphisme somme, avec pour accouplement
$$
\mathop{\rm Ker}({\Bbb Z}^{I}\rightarrow {\Bbb Z})\times T\rightarrow
{\Bbb G}_{{\rm m},k},~(\lambda ,t)\mapsto t^{\lambda}=\prod_{i\in
I}t_{i}^{\lambda_{i}}.
$$

L'homomorphisme de $k$-sch\'{e}mas en groupes
$$
T\rightarrow\mathop{\rm Pic}\nolimits_{\overline{P}/k}
$$
obtenu en composant l'inclusion $T\subset P^{0}$ et l'homomorphisme
canonique $\beta$ d'Esteves, Gagn\'{e} et Kleiman (cf.  2.4),
d\'{e}finit donc d'apr\`{e}s la proposition 2.5.1 un
$\Lambda^{0}$-torseur $\overline{P}{}^{\dagger}
\rightarrow\overline{P}$.

\thm PROPOSITION 2.6.1
\enonce
Le $\Lambda^{0}$-torseur $\overline{P}{}^{\dagger}\rightarrow
\overline{P}$ n'est autre que le rev\^{e}tement
$\overline{P}^{\,\natural}\rightarrow\overline{P}$.
\endthm

\rem Preuve
\endrem
Il s'agit de d\'{e}montrer que, pour chaque entier $d$, la restriction
\`{a} $Z^{d}$ du $\Lambda^{0}$-torseur $\overline{P}{}^{\dagger}
\rightarrow \overline{P}$ est isomorphe \`{a} $X^{d}\rightarrow
Z^{d}=X^{d}/\Lambda^{0}$ et, bien s\^{u}r, il suffit de le faire pour
$d=0$.

Si l'on interpr\`{e}te le $k$-sch\'{e}ma de Picard $\mathop{\rm Pic}
\nolimits_{Z^{0}/k}$ de $Z^{0}$ comme le $k$-sch\'{e}ma de Picard
$\Lambda^{0}$-\'{e}quivariant $\mathop{\rm Pic}
\nolimits_{X^{0}/k}^{\Lambda^{0}}$ de $X^{0}$, le
$\Lambda^{0}$-torseur $X^{0}\rightarrow Z^{0}$ correspond de
mani\`{e}re tautologique par la proposition 2.5.1 au $k$-homomorphisme
$$
T\rightarrow \mathop{\rm Pic}\nolimits_{X^{0}/k}^{\Lambda^{0}}
=\mathop{\rm Pic}\nolimits_{Z^{0}/k}
$$
qui envoie $t\in T(k)$ sur le ${\cal O}_{X^{0}}$-Module inversible trivial
muni de l'action de $\Lambda^{0}$ donn\'{e}e par le caract\`{e}re
$$
\chi_{t}:\Lambda^{0}\rightarrow {\Bbb G}_{{\rm m},k},~\lambda\mapsto
t^{\lambda},
$$
et on veut montrer que ce $k$-homomorphisme co\"{\i}ncide avec le
$k$-homomorphisme compos\'{e}
$$
T\subset P^{0}\,\smash{\mathop{\hbox to 8mm{\rightarrowfill}}
\limits^{\scriptstyle \beta_{0}}}\,\mathop{\rm Pic}
\nolimits_{\overline{P}{}^{0}/k}\rightarrow \mathop{\rm Pic}
\nolimits_{Z^{0}/k}
$$
o\`{u} la derni\`{e}re fl\`{e}che est la fl\`{e}che de restriction par
le morphisme radiciel $Z^{0}\rightarrow \overline{P}{}^{0}$.  Il
suffit donc de d\'{e}montrer que, pour tout $t\in T(k)$, la
restriction \`{a} $Z^{0}$ du ${\cal O}_{\overline{P}{}^{0}}$-Module
inversible $\beta_{0}(t)$ est le ${\cal O}_{X^{0}}$-Module inversible
trivial muni de l'action de $\Lambda^{0}$ donn\'{e}e par le
caract\`{e}re $\chi_{t}$.

Rappelons que l'anneau total des fractions de $A$ est
$E=E_{I}=\prod_{i\in I}E_{i}$, que $\widetilde{C}$ est la droite
projective standard ${\Bbb P}_{k}^{1}$, que le point $\infty$ de
$\widetilde{C}$ est distinct des $\widetilde{c}_{i}$, que pour chaque
$i\in I$, on peut prendre pour uniformisante
$\varpi_{i}=\varpi_{E_{i}}$ de $E_{i}$ la fonction
$x-\widetilde{c}_{i}$ o\`{u} $x$ est une coordonn\'{e}e affine sur
$\widetilde{C}\setminus\{\infty\}\cong {\Bbb A}_{k}^{1}$, que
$\widetilde{A}=\prod_{i\in I}k[[\varpi_{i}]]\subset\prod_{i\in I}
k((\varpi_{i}))\cong E_{I}$ et que $\Lambda^{0}$ est l'ensemble des
fractions rationnelles de la forme $\prod_{i\in I}
(x-\widetilde{c}_{i})^{\lambda_{i}}$ avec $\lambda\in {\Bbb Z}^{I}$ et
$\sum_{i\in I}\lambda_{i}=0$.

Rappelons de plus que $X^{0}$ est le $k$-sch\'{e}ma des
$A$-r\'{e}seaux $M\subset E$ d'indice $0$ relativement \`{a} $A\subset
E$, que l'image ${\cal M}\in \overline{P}{}^{0}(k)$ de $M\in X^{0}(k)$
est le ${\cal O}_{C}$-Module obtenu en recollant ${\cal
O}_{C\setminus\{c\}}$ et $M$, et que l'action du groupe de Galois
$\Lambda^{0}$ est donn\'{e}e par
$$
(\lambda ,M)\mapsto\lambda\cdot M=(\varpi_{i}^{\lambda_{i}})_{i\in I}M.
$$

Rappelons enfin que le $k$-sch\'{e}ma $X^{0}$ est localement de type
fini, r\'{e}union d'une suite croissante
$$
X_{0}^{0}\subset\cdots\subset X_{n}^{0}\subset X_{n+1}^{0}\subset
\cdots\subset X^{0}
$$
de ferm\'{e}s de type fini stables sous l'action par translation de
$P^{0}(k)=\widetilde{A}^{\times}/A^{\times}$; par exemple, on peut
prendre pour $X_{n}^{0}$ le ferm\'{e} form\'{e} des $M$ qui sont
contenus dans $(\varpi_{i}^{-n})_{i\in I}\widetilde{A}\subset E$.

Il suffit donc de construire, un syst\`{e}me compatible de
trivialisations des restrictions du Module inversible $\beta_{0}(t)$
aux $X_{n}^{0}$ et de montrer que, pour tout $M\in X^{0}(k)$ et tout
$\lambda\in\Lambda^{0}$, l'isomorphisme de la fibre en $M$ de la
restriction de $\beta_{0}(t)$ \`{a} $X^{0}$ sur celle en $\lambda\cdot
M$ donn\'{e} par l'action de $\lambda$ s'exprime dans ce syst\`{e}me
de trivialisations comme la multiplication par $t^{\lambda}$.

La fibre de $\beta_{0}(t)$ en l'image ${\cal M}\in\overline{P}{}^{0}
(k)$ de $M\in X_{n}^{0}(k)$ peut se calculer comme suit.  Soit
$R=H^{0}(C\setminus\{c\},{\cal O}_{C})\subset E$ la $k$-alg\`{e}bre
des fonctions rationnelles sur $\widetilde{C}$ qui sont
r\'{e}guli\`{e}res en dehors de $\{\widetilde{c}_{i}\mid i\in I\}$.
On a une suite exacte
$$
0\rightarrow H^{0}(C,{\cal M})\rightarrow R\rightarrow E/M\rightarrow
H^{1}(C,{\cal M})\rightarrow 0
$$
et il existe $N\geq n$, ind\'{e}pendant de $M\in X_{n}^{0}(k)$, tel
que la restriction de la fl\`{e}che surjective $E/M\twoheadrightarrow
H^{1}(C,{\cal M})$ \`{a} $(\varpi_{i}^{-N})_{i\in I}\widetilde{A}/M\subset
E/M$ soit encore surjective.  La suite induite
$$
0\rightarrow H^{0}(C,{\cal M})\rightarrow R\cap (\varpi_{i}^{-N})_{i\in I}
\widetilde{A}\rightarrow (\varpi_{i}^{-N})_{i\in I}\widetilde{A}/M
\rightarrow H^{1}(C,{\cal M})\rightarrow 0.
$$
est donc aussi exacte.

La classe d'isomorphie fix\'{e}e $t\in T(k)\subset P^{0}(k)$ est
repr\'{e}sent\'{e}e par le ${\cal O}_{C}$-Module inversible ${\cal L}$
obtenu en recollant ${\cal O}_{C\setminus\{c\}}$ et le r\'{e}seau
$L=(t_{i})_{i\in I}A\subset E$, et on peut remplacer dans les suites
exactes ci-dessus ${\cal M}$ par ${\cal L}\otimes_{{\cal O}_{C}}{\cal
M}$ et $M$ par $t\cdot M=(t_{i})_{i\in I}M\subset
(\varpi_{i}^{-n})_{i\in I}\widetilde{A}\subset E$.  Par
cons\'{e}quent, la fibre de $\beta_{0}(t)$ en ${\cal M}$ s'identifie
canoniquement \`{a} la $k$-droite
$$
\left(\bigwedge^{{\rm max}}(\varpi_{i}^{-N})_{i\in I}\widetilde{A}
\big/t\cdot M\right) \otimes_{k}\left(\bigwedge^{{\rm max}}
(\varpi_{i}^{-N})_{i\in I}\widetilde{A}\big/M\right)^{\otimes -1},
$$
ou se qui revient au m\^{e}me \`{a} la $k$-droite
$$
\left(\bigwedge^{{\rm max}}(\varpi_{i}^{-n})_{i\in I}\widetilde{A}
\big/t\cdot M\right) \otimes_{k}\left(\bigwedge^{{\rm max}}
(\varpi_{i}^{-n})_{i\in I}\widetilde{A}\big/M\right)^{\otimes -1}.
$$
Le d\'{e}terminant de l'isomorphisme $(\varpi_{i}^{-n})_{i\in
I}\widetilde{A}/M \buildrel\sim\over\longrightarrow (\varpi_{i}^{-n})_{i\in
I}\widetilde{A}/t\cdot M$ induit par la multiplication par
$(t_{i})_{i\in I}$ d\'{e}finit un vecteur de base $e_{n,M}$ de cette
derni\`{e}re $k$-droite.  La section $M\mapsto e_{M}=\left(\prod_{i\in
I}t_{i}^{-n}\right)e_{n,M}$ de la restriction de $\beta_{0}(t)$ \`{a}
$X_{n}^{0}$ est la trivialisation cherch\'{e}e: elle est
{\og}{ind\'{e}pendant de $n$}{\fg} puisque le d\'{e}terminant de
l'automorphisme de multiplication par $(t_{i})_{i\in I}$ sur
$(\varpi_{i}^{-n-1})_{i\in I}\widetilde{A}/ (\varpi_{i}^{-n})_{i\in
I}\widetilde{A}$ est \'{e}gal \`{a} $\prod_{i\in I}t_{i}$.

Maintenant, si $M$ et $\lambda\cdot M$ sont dans $X_{n}^{0}(k)$, la
multiplication par $(\varpi_{i}^{-\lambda_{i}})_{i\in I}$ induit un
isomorphisme de la $k$-droite
$$
\left(\bigwedge^{{\rm max}}(\varpi_{i}^{-n})_{i\in I}\widetilde{A}
\big/t\cdot (\lambda\cdot M)\right)\otimes_{k}\left(\bigwedge^{{\rm
max}}(\varpi_{i}^{-n})_{i\in I}\widetilde{A}\big/\lambda\cdot
M\right)^{\otimes -1}
$$
sur la $k$-droite
$$
\left(\bigwedge^{{\rm max}}(\varpi_{i}^{-\lambda_{i}-n})_{i\in I}
\widetilde{A}\big/t\cdot M\right) \otimes_{k}\left(\bigwedge^{{\rm
max}}(\varpi_{i}^{-\lambda_{i}-n})_{i\in I}\widetilde{A}\big/
M\right)^{\otimes -1}
$$
qui envoie le vecteur de base $e_{n,M}$ sur le d\'{e}terminant $e'$ de
l'isomorphisme
$$
(\varpi_{i}^{-\lambda_{i}-n})_{i\in I}\widetilde{A}/M\buildrel\sim\over
\longrightarrow (\varpi_{i}^{-\lambda_{i}-n})_{i\in I}\widetilde{A}/t\cdot M
$$
induit par la multiplication par $(t_{i})_{i\in I}$.  Or cette
derni\`{e}re $k$-droite est canoniquement isomorphe \`{a} la
$k$-droite
$$
\left(\bigwedge^{{\rm max}}(\varpi_{i}^{-n})_{i\in I}\widetilde{A}
\big/t\cdot M\right)\otimes_{k}\left(\bigwedge^{{\rm max}}
(\varpi_{i}^{-n})_{i\in I}\widetilde{A}\big/M\right)^{\otimes -1}
$$
par un isomorphisme qui envoie $e'$ sur $t^{-\lambda}e_{n,M}$ (choisir
arbitrairement un entier $m\geq \lambda_{i}+n$ quel que soit $i\in I$
et utiliser les inclusions $(\varpi_{i}^{-\lambda_{i}-n})_{i\in I}
\widetilde{A}\subset (\varpi_{i}^{-m})_{i\in I}\widetilde{A}\supset
(\varpi_{i}^{-n})_{i\in I} \widetilde{A}$).  On en d\'{e}duit donc que
l'action de $-\lambda$ envoie $e_{\lambda\cdot M}$ sur
$t^{-\lambda}e_{M}$, ce que l'on voulait d\'{e}montrer.
\hfill\hfill$\square$
\vskip 3mm

\rem Remarque $2.6.2$
\endrem
Pour tout $\widetilde{a}\in \widetilde{A}^{\times}/A^{\times}$, la
multiplication par $\widetilde{a}$ induit un isomorphisme de
$(\varpi_{i}^{-n})_{i\in I}\widetilde{A}/M$ sur $(\varpi_{i}^{-n})_{i\in
I}\widetilde{A}/\widetilde{a}M$ dont le d\'{e}terminant ne d\'{e}pend
que du {\og}{terme constant}{\fg} $\widetilde{a}(0)$ de
$\widetilde{a}$, c'est-\`{a}-dire de l'image de $\widetilde{a}$ par la
projection canonique $\widetilde{A}^{\times}/
A^{\times}\twoheadrightarrow (k^{\times})^{I}/k^{\times}$.  Les
m\^{e}mes calculs de d\'{e}terminant de la cohomologie que ceux de la
preuve de la proposition montrent donc que la fl\`{e}che compos\'{e}e
$$
P^{0}\,\smash{\mathop{\hbox to 8mm{\rightarrowfill}}
\limits^{\scriptstyle\beta}}\,\mathop{\rm Pic}\nolimits_{\overline{P}}
\rightarrow \mathop{\rm Pic}\nolimits_{Z}
$$
est triviale sur la composante unipotente de $P^{0}$.
\hfill\hfill$\square$
\vskip 3mm

\rem Remarque $2.6.3$
\endrem
On a bien s\^{u}r une variante de la proposition 2.6.1 pour le
rev\^{e}tement du corollaire 2.2.6.
\hfill\hfill$\square$
\vskip 5mm

\centerline{3. D\'{E}FORMATIONS DE COURBES (RAPPELS)}
\vskip 2mm

\section{3.1}{D\'{e}formations miniverselles : g\'{e}n\'{e}ralit\'{e}s}

On note $\mathop{\rm Art}\nolimits_{k}$ la cat\'{e}gorie des
$k$-alg\`{e}bres locales artiniennes de corps r\'{e}siduel $k$. Tous
les $k$-sch\'{e}mas consid\'{e}r\'{e}s dans la suite sont suppos\'{e}s
s\'{e}par\'{e}s et localement n{\oe}th\'{e}riens.

Soit $X_{k}$ un $k$-sch\'{e}ma. On a le foncteur des {\it
d\'{e}formations}
$$
\mathop{\rm Def}\nolimits_{X_{k}}:\mathop{\rm Art}\nolimits_{k}
\rightarrow \mathop{\rm Ens}
$$
qui associe \`{a} une $k$-alg\`{e}bre locale artinienne $R$ l'ensemble
des classes d'isomorphie de $R$-sch\'{e}mas plats $X_{R}$ munis d'un
isomorphisme de $k$-sch\'{e}mas
$$
\iota :X_{k}\buildrel\sim\over\longrightarrow k\otimes_{R}X_{R}.
$$

Si $X_{k}=\mathop{\rm Spec}(A_{k})$ est affine, toute d\'{e}formation
$X_{R}$ de $X_{k}$ sur $R\in\mathop{\rm ob}\mathop{\rm Art}
\nolimits_{k}$ est aussi un sch\'{e}ma affine $X_{R}=\mathop{\rm
Spec}(A_{R})$ o\`{u} $A_{R}$ est une $R$-alg\`{e}bre plate munie d'un
isomorphisme de $k$-alg\`{e}bres $k\otimes_{R}A_{R}\buildrel\sim
\over\longrightarrow A_{k}$. On identifiera donc $\mathop{\rm Def}
\nolimits_{X_{k}}$ au foncteur $\mathop{\rm Def}\nolimits_{A_{k}}$ qui
envoie $R$ sur l'ensemble des classes d'isomorphie des $A_{R}$.

On d\'{e}finit aussi le foncteur des d\'{e}formations
$$
\mathop{\rm Def}\nolimits_{X_{k}}^{{\rm top}}=\mathop{\rm
Def}\nolimits_{A_{k}}^{{\rm top}}
$$
d'un $k$-sch\'{e}ma formel affine $X_{k}=\mathop{\rm Spf}(A_{k})$ en
tenant compte de la topologie de $A_{k}$.

Si $U_{k}$ est un ouvert de $X_{k}$ et si $x\in U_{k}$, on a des
morphismes de foncteurs \'{e}vidents
$$
\mathop{\rm Def}\nolimits_{X_{k}}\rightarrow\mathop{\rm Def}
\nolimits_{U_{k}}\rightarrow \mathop{\rm Def}\nolimits_{{\cal
O}_{X_{k},x}}\rightarrow \mathop{\rm Def} \nolimits_{\widehat{{\cal
O}}_{X_{k},x}}^{{\rm top}}.
$$

Une {\it d\'{e}formation formelle} $({\cal R},{\cal X})$ de $X_{k}$
est un couple form\'{e} d'une $k$-alg\`{e}bre locale compl\`{e}te
n{\oe}th\'{e}rienne ${\cal R}$ de corps r\'{e}siduel $k$ et d'un
couple ${\cal X}=(X_{\bullet},\iota_{\bullet})$ o\`{u} $X_{\bullet}$
est une suite de d\'{e}formations $X_{n}$ de $X_{k}$ sur $R_{n}={\cal
R}/{\frak m}_{{\cal R}}^{n+1}$ avec $X_{0}=X_{k}$ et o\`{u}
$\iota_{\bullet}$ est une suite d'isomorphismes
$$
\iota_{n}:X_{n}\buildrel\sim\over\longrightarrow
R_{n}\otimes_{R_{n+1}}X_{n+1},~n\in {\Bbb N}.
$$
On verra dans la suite ${\cal X}$ comme un ${\cal
S}$-sch\'{e}ma formel o\`{u} ${\cal S}=\mathop{\rm Spf}({\cal R})$.

Une telle d\'{e}formation d\'{e}finit un morphisme de foncteurs
$$
{\cal S}\rightarrow\mathop{\rm Def}\nolimits_{X_{k}}
$$
qui envoie le $R$-point $\varphi :{\cal R}\rightarrow R$ de ${\cal S}$
sur $(R\otimes_{\varphi ,{\cal R}}{\cal X},\iota )$ o\`{u} on a
pos\'{e}
$$
R\otimes_{\varphi ,{\cal R}}{\cal X}=R\otimes_{\varphi_{n},R_{n}}
X_{n}
$$
quel que soit l'entier $n$ assez grand pour que $\varphi$ se factorise
en ${\cal R}\twoheadrightarrow R_{n}\,\smash{\mathop{\hbox to
6mm{\rightarrowfill}} \limits^{\scriptstyle \varphi_{n}}}\,R$.

\thm D\'{E}FINITION 3.1.1
\enonce
Une d\'{e}formation formelle $({\cal R},{\cal X})$ de $X_{k}$ est dite
{\rm miniverselle} si le morphisme de foncteurs $\mathop{\rm Hom}
\nolimits_{k{\rm -alg}}({\cal R},\cdot )\rightarrow \mathop{\rm Def}
\nolimits_{X_{k}}$ associ\'{e} est formellement lisse et si
l'application entre les espaces tangents
$$
\mathop{\rm Hom}\nolimits_{k}({\frak m}_{{\cal R}}/{\frak m}_{{\cal
R}}^{2},k)\buildrel\sim\over\longleftarrow {\cal S}(k[\varepsilon])
\rightarrow\mathop{\rm Def}\nolimits_{X_{k}}(k[\varepsilon ])
$$
est bijective.
\endthm

On d\'{e}finit de m\^{e}me une d\'{e}formation formelle miniverselle
d'une $k$-alg\`{e}bre $A_{k}$ et une d\'{e}formation formelle
miniverselle topologique d'une $k$-alg\`{e}bre topologique ${\cal
A}_{k}$. Une d\'{e}formation formelle miniverselle de $X_{k}$ ou
$A_{k}$ (resp. une d\'{e}formation formelle miniverselle topologique
$A_{k}$) est unique \`{a} isomorphisme pr\`{e}s.

\thm TH\'{E}OR\`{E}ME 3.1.2 (Schlessinger, cf. [Ri] Theorem 4.5)
\enonce
Soit $X_{k}$ un sch\'{e}ma s\'{e}par\'{e} et de type fini sur $k$.
Supposons de plus que
\vskip 1mm

\itemitem{-}soit $X_{k}$ est propre sur $k$,
\vskip 1mm

\itemitem{-} soit $X_{k}$ est affine avec un lieu singulier
$X_{k}^{{\rm sing}}$ fini sur $k$.
\vskip 1mm

Alors, $X_{k}$ admet une d\'{e}formation formelle miniverselle $({\cal
R},{\cal X})$.
\hfill\hfill$\square$
\endthm

\thm TH\'{E}OR\`{E}ME 3.1.3 (Rim, [Ri] Theorem 4.11, Corollary 4.13)
\enonce
Soit $X_{k}$ un sch\'{e}ma s\'{e}par\'{e} et de type fini sur $k$
dont le lieu singulier $X_{k}^{{\rm sing}}$ est fini sur $k$. Pour chaque $x\in
X_{k}^{{\rm sing}}$, soit $U_{k,x}$ un ouvert affine de $X_{k}$ qui ne
rencontre $X_{k}^{{\rm sing}}$ qu'en $x$. On a alors les
propri\'{e}t\'{e}s suivantes.

\decale{\rm (i)} Si $H^{2}(X_{k},\mathop{{\cal H}{\it om}}
\nolimits_{{\cal O}_{X_{k}}}(\Omega_{X_{k}/k}^{1},{\cal O}_{X_{k}}))
=(0)$, le morphisme de foncteurs
$$
\mathop{\rm Def}\nolimits_{X_{k}}\rightarrow\prod_{x\in X_{k}^{{\rm
sing}}}\mathop{\rm Def}\nolimits_{U_{k,x}}
$$
est formellement lisse.

\decale{\rm (ii)} Pour chaque $x\in X_{k}^{{\rm sing}}$, le morphisme
de foncteurs
$$
\mathop{\rm Def}\nolimits_{U_{k,x}}\rightarrow\mathop{\rm
Def}\nolimits_{\widehat{{\cal O}}_{X,x}}^{{\rm top}}
$$
est formellement lisse et induit une bijection entre les espaces
tangents de Zariski
$$
\mathop{\rm Def}\nolimits_{U_{k,x}}(k[\varepsilon ])
\buildrel\sim\over\longrightarrow \mathop{\rm Def}
\nolimits_{\widehat{{\cal O}}_{X,x}}^{{\rm top}}(k[\varepsilon ]);
$$
en particulier, si $({\cal R},{\cal U}_{x})$ est une d\'{e}formation
formelle miniverselle du $k$-sch\'{e}ma affine $U_{k,x}$, alors $({\cal
R},\widehat{{\cal O}}_{{\cal U}_{x},x})$ est une d\'{e}formation
formelle miniverselle topologique de $\widehat{{\cal O}}_{X,x}$.

\decale{\rm (iii)} Si $H^{2}(X_{k},\mathop{{\cal H}{\it om}}
\nolimits_{{\cal O}_{X_{k}}}(\Omega_{X_{k}/k}^{1},{\cal
O}_{X_{k}}))=(0)$ et si $X_{k}$ est localement d'intersection
compl\`{e}te, tous les foncteurs de d\'{e}formations
consid\'{e}r\'{e}s dans {\rm (i)} et {\rm (ii)} sont formellement
lisses sur $k$; si de plus le lieu de lissit\'{e} est dense dans
$X_{k}$, la suite exacte entre les espaces tangents associ\'{e}e au
morphisme de foncteurs
$$
\mathop{\rm Def}\nolimits_{X_{k}}\rightarrow\prod_{x\in X_{k}^{{\rm
sing}}}\mathop{\rm Def}\nolimits_{\widehat{{\cal O}}_{X_{k},x}}^{{\rm top}}
$$
n'est autre que la suite exacte courte
$$\displaylines{
\qquad 0\rightarrow H^{1}(X_{k},\mathop{{\cal H}{\it om}}
\nolimits_{{\cal O}_{X_{k}}}(\Omega_{X_{k}/k}^{1},{\cal O}_{X_{k}}))
\rightarrow \mathop{\rm Ext}\nolimits_{{\cal O}_{X_{k}}}^{1}
(\Omega_{X_{k}/k}^{1},{\cal O}_{X_{k}})
\hfill\cr\hfill
\rightarrow H^{0}(X_{k}, \mathop{{\cal E}{\it xt}}\nolimits_{{\cal
O}_{X_{k}}}^{1}(\Omega_{X_{k}/k}^{1},{\cal O}_{X_{k}}))\rightarrow
0.\qquad\square}
$$
\endthm

\section{3.2}{D\'{e}formations miniverselles: le cas des courbes
\`{a} singularit\'{e}s planes}

Soit $C_{k}$ une courbe int\`{e}gre, projective et localement
d'intersection compl\`{e}te sur $k$; $C_{k}$ est g\'{e}n\'{e}riquement
lisse sur $k$ puisque l'on a suppos\'{e} $k$ alg\'{e}briquement clos.

Le $k$-sch\'{e}ma $C_{k}$ admet une d\'{e}formation formelle
miniverselle $({\cal R},{\cal C})$ d'apr\`{e}s le th\'{e}or\`{e}me
3.1.2 de Schlessinger.

\thm D\'{E}FINITION 3.2.1
\enonce
Soit $c$ un point ferm\'{e} de $C_{k}$.  Nous dirons que $C_{k}$ admet
en $c$ une {\rm singularit\'{e} plane isol\'{e}e} si le
compl\'{e}t\'{e} $\widehat{{\cal O}}_{C_{k},c}$ de l'anneau local de
$C_{k}$ en $c$ est une $k$-alg\`{e}bre topologique isomorphe \`{a}
$k[[x,y]]/(f)$ pour une s\'{e}rie formelle $f=f(x,y)\in (x,y)\subset
k[[x,y]]$ telle que le $k$-espace vectoriel $k[[x,y]]/(\partial_{x}f,
\partial_{y}f)$ soit de dimension finie.
\endthm

Bien entendu, si le point $c$ est une singularit\'{e} plane isol\'{e}e
de $C_{k}$ au sens de la d\'{e}finition ci-dessus, il l'est aussi au
sens usuel puisque le $k$-espace vectoriel
$k[[x,y]]/(f,\partial_{x}f,\partial_{y}f)$ est a fortiori de dimension
finie.  La r\'{e}ciproque est vraie si $k$ est de caract\'{e}ristique
nulle, puisque $f$ est alors entier sur l'id\'{e}al $(\partial_{x}f,
\partial_{y}f)$ (cf. [Te 1] 1.1).
\vskip 2mm

Supposons dans la suite que les seules singularit\'{e}s de $C_{k}$
sont des singularit\'{e}s planes isol\'{e}es.

Comme $C_{k}$ est de dimension $1$ et localement d'intersection
compl\`{e}te, on a
$$
H^{2}(C_{k},\mathop{{\cal H}{\it om}}\nolimits_{{\cal O}_{C_{k}}}
(\Omega_{C_{k}/k}^{1},{\cal O}_{C_{k}})) =(0)
$$
et il n'y a pas d'obstruction \`{a} d\'{e}former $C_{k}$.  Par suite
la base ${\cal R}$ de la d\'{e}formation miniverselle est une
$k$-alg\`{e}bre de s\'{e}ries formelles en $\tau (C_{k})$ variables,
o\`{u} $\tau (C_{k})$ est la dimension de l'espace tangent de Zariski
$$
\mathop{\rm Ext}\nolimits_{{\cal O}_{C_{k}}}^{1}
(\Omega_{C_{k}/k}^{1},{\cal O}_{C_{k}}).
$$
De plus, le $k$-espace vectoriel sous-jacent \`{a} l'alg\`{e}bre de
Lie du $k$-sch\'{e}ma en groupes $\mathop{\rm Aut}\nolimits_{k}
(C_{k})$ des $k$-automorphismes de $C_{k}$ est \'{e}gal \`{a}
$\mathop{\rm Hom} \nolimits_{{\cal O}_{C_{k}}}(\Omega_{C_{k}/k}^{1},
{\cal O}_{C_{k}})$.

De m\^{e}me, pour chaque point singulier $c\in C_{k}^{{\rm sing}}$, le
foncteur des d\'{e}formations topologiques de $\widehat{{\cal
O}}_{C_{k},c}\cong k[[x,y]]/(f)$ est formellement lisse avec pour
espace tangent
$$
\mathop{\rm Def}\nolimits_{\widehat{{\cal O}}_{C_{k},c}}^{{\rm
top}}(k[\varepsilon]) =\mathop{\rm Ext}\nolimits_{\widehat{{\cal
O}}_{C_{k},c}}^{1} (\Omega_{C_{k},c}^{1,{\rm top}},\widehat{{\cal
O}}_{C_{k},c})\cong k[[x,y]]/(f,\partial_{x}f,\partial_{y}f)
$$
(on a
$$
\Omega_{C_{k},c}^{1,{\rm top}}=\mathop{\rm Coker}(k[[x,y]]/(f)
\hookrightarrow (k[[x,y]]/(f)){\rm d}x\oplus (k[[x,y]]/(f)){\rm d}y)
$$
o\`{u} la fl\`{e}che envoie $1$ sur ${\rm d}f=(\partial_{x}f){\rm
d}x+(\partial_{y}f){\rm d}y$). De plus, on obtient une d\'{e}formation
miniverselle topologique
$$
\mathop{\rm Spf}(k[[x,y,z_{1},\ldots ,z_{\tau_{c}(C_{k}))}]]/
(\widetilde{f}))\rightarrow\mathop{\rm Spf}(k[[z_{1},\ldots ,
z_{\tau_{c}(C_{k})}]])
$$
du germe formel de $C_{k}$ en $c$ en choisissant arbitrairement des
s\'{e}rie formelles
$$
g_{1},\ldots ,g_{\tau_{c}(C_{k})}\in k[[x,y]]
$$
dont les r\'{e}ductions modulo l'id\'{e}al $(f,\partial_{x}f,
\partial_{y}f)$ forment une base de $\mathop{\rm Def}
\nolimits_{\widehat{{\cal O}}_{C_{k},c}}^{{\rm top}}
(k[\varepsilon])$, et en posant
$$
\widetilde{f}=\widetilde{f}(x,y,z)=f(x,y)
+\sum_{t=1}^{\tau_{c}(C_{k})}z_{t}g_{t}(x,y).
$$
Bien s\^{u}r, on peut prendre $g_{\tau_{c}(C_{k})}=1$, de sorte que
$k[[x,y,z_{1},\ldots ,z_{\tau_{c}(C_{k})}]]/(\widetilde{f})$ est
isomorphe \`{a} la $k$-alg\`{e}bre de s\'{e}ries formelles
$k[[x,y,z_{1},\ldots ,z_{\tau_{c}(C_{k})-1}]]$.
\vskip 2mm

Il r\'{e}sulte de cette discussion que le morphisme propre et plat de
$k$-sch\'{e}mas formels ${\cal C}\rightarrow \mathop{\rm Spf}({\cal
R})$ est localement d'intersection compl\`{e}te.  On peut donc
consid\'{e}rer son ${\cal O}_{{\cal C}}$-Module dualisant relatif
$\omega_{{\cal C}/{\cal S}}$; la restriction de ce Module inversible
\`{a} $C_{k}$ est bien entendu le ${\cal O}_{C_{k}}$-Module dualisant
$\omega_{C_{k}/k}$.

\thm TH\'{E}OR\`{E}ME 3.2.2
\enonce
Soient $S=\mathop{\rm Spec}({\cal R})$ et $s=\mathop{\rm Spec}(k)$
l'unique point ferm\'{e} de $S$.  La d\'{e}formation formelle
miniverselle $({\cal R},{\cal C})$ de $C_{k}$ est alg\'{e}brisable :
il existe un $S$-sch\'{e}ma propre et plat $C$ dont le
compl\'{e}t\'{e} formel pour la topologie ${\frak m}_{{\cal
R}}$-adique est le ${\cal R}$-sch\'{e}ma formel ${\cal C}$.
\hfill\hfill$\square$
\endthm

\rem Preuve
\endrem
Le degr\'{e} de $\omega_{C_{k}/k}$ est $2p(C_{k})-2$ o\`{u}
$p(C_{k})=\mathop{\rm dim} \nolimits_{k}H^{1}(C_{k},{\cal
O}_{C_{k}})$ est le genre arithm\'{e}tique de $C_{k}$.

Si $p(C_{k})\geq 2$, $\omega_{C_{k}/k}$ est ample et on peut
appliquer directement le th\'{e}or\`{e}me d'alg\'{e}brisation de
Grothendieck ([Gr 3] Th\'{e}or\`{e}me (5.4.5)).

Si $p(C_{k})=0\hbox{ ou }1$, on fixe arbitrairement un $k$-point
$\infty$ du lieu de lissit\'{e} de $C_{k}$ et un rel\`{e}vement de ce
$k$-point en un ${\cal R}$-point de ${\cal C}$ not\'{e} encore
$\infty$.  Le ${\cal O}_{C_{k}}$-Module inversible
$\omega_{C_{k}/k}((3-2p(C_{k}))[\infty ])$ est ample et on peut
appliquer de nouveau th\'{e}or\`{e}me d'alg\'{e}brisation de
Grothendieck.
\hfill\hfill$\square$
\vskip 3mm

Bien s\^{u}r, le morphisme propre et plat de sch\'{e}mas $C\rightarrow
S$ du th\'{e}or\`{e}me est projectif et sa fibre $C_{t}$ en tout point
g\'{e}om\'{e}trique $t$ de $S$ est int\`{e}gre (cf.  [Gr 4]
Th\'{e}or\`{e}me (12.2.4)).

\thm PROPOSITION 3.2.3
\enonce
Soient $C_{k}$ une courbe int\`{e}gre et projective sur $k$, qui
n'admet que des singularit\'{e}s planes, et $C\rightarrow S$ une
alg\'{e}brisation d'une d\'{e}formation formelle miniverselle de
$C_{k}$.  Alors, le sch\'{e}ma $C$ est formellement lisse sur $k$ et
la fibre g\'{e}n\'{e}rique du morphisme $C\rightarrow S$ est lisse.
\endthm

\rem Preuve
\endrem
La premi\`{e}re assertion est locale en chaque point singulier $c\in
C_{k}$ et r\'{e}sulte imm\'{e}diatement des \'{e}critures
$$
\mathop{\rm Spf}(k[[x,y]]/(f))\hbox{ et }\mathop{\rm Spf}
(k[[x,y,z_{1},\ldots ,z_{\tau_{c}(C_{k})}]]/ (\widetilde{f}))
$$
pour les compl\'{e}t\'{e}s formels de $C_{k}$ et $C$ en $c$, avec
$$
\widetilde{f}=\widetilde{f}(x,y,z)=f(x,y)
+\sum_{t=1}^{\tau_{c}(C_{k})}z_{t}g_{t}(x,y).
$$
et $g_{\tau_{c}(C_{k})}=1$.

Passons \`{a} la deuxi\`{e}me assertion.  On sait d\'{e}j\`{a} que la
fibre $C_{\eta}$ de $C\rightarrow S$ au point g\'{e}n\'{e}rique $\eta$
de $S$ est g\'{e}om\'{e}triquement int\`{e}gre et donc
g\'{e}n\'{e}riquement lisse.

Maintenant, si $C_{\eta}$ admettait une singularit\'{e} en un point
ferm\'{e} $c_{\eta}$, on pourrait localiser la situation au voisinage
du point ferm\'{e} $c$ de $C_{s}$ qui sp\'{e}cialise $c_{\eta}$.  Pour
conclure, il suffit donc de v\'{e}rifier que, pour $\widetilde{f}$
comme ci-dessus, le lieu critique de la projection canonique
$$
\mathop{\rm Spf}(k[[x,y,z_{1},\ldots ,z_{\tau_{c}(C_{k})}]]/
(\widetilde{f}))\rightarrow\mathop{\rm Spf}(k[[z_{1},\ldots ,
z_{\tau_{c}(C_{k})}]])={\cal S}\leqno{(\ast )}
$$
est fini sur ${\cal S}$ et son discriminant, c'est-\`{a}-dire l'image
dans ${\cal S}$ du lieu critique, est un ferm\'{e} strict de ${\cal
S}$.

Or, le lieu critique de $(\ast )$ est le ferm\'{e} de $\mathop{\rm
Spf} (k[[x,y,z_{1},\ldots ,z_{\tau_{c}(C_{k})}]])$ d'\'{e}quations
$$
\widetilde{f}=\partial_{x}\widetilde{f}=\partial_{y}\widetilde{f}=0.
$$
Il est donc bien fini sur ${\cal S}$ puisque $k[[x,y]]/(f,
\partial_{x}f,\partial_{y}f)$ est de dimension finie sur $k$, et le
discriminant de $(\ast )$ est le ferm\'{e} d\'{e}fini par le
$0$-\`{e}me id\'{e}al de Fitting du $k[[z_{1},\ldots ,
z_{\tau_{c}(C_{k})}]]$-module de type fini $k[[x,y,z_{1},\ldots ,
z_{\tau_{c}(C_{k})}]]/(\widetilde{f},\partial_{x}\widetilde{f},
\partial_{y}\widetilde{f})$.

Si le discriminant de $(\ast )$ \'{e}tait ${\cal S}$ tout entier, il
contiendrait en particulier tout l'axe des $z_{1}$.  Or la formation
de ce discriminant commute \`{a} tout changement de base (cf.  [Te~2]
\S 5).  Sa restriction \`{a} l'axe des $z_{1}$ est donc le
discriminant du morphisme
$$
\mathop{\rm Spf}(k[[x,y]])\rightarrow\mathop{\rm Spf}(k[[z_{1}]])
$$
qui envoie $(x,y)$ sur $z_{1}=-f(x,y)$.  Mais, il est facile de voir
(cf.  [Te 2] \S 2.6) que ce dernier discriminant est le ferm\'{e} de
$\mathop{\rm Spf}(k[[z_{1}]])$ d\'{e}fini par l'\'{e}quation
$z_{1}^{\mu (f)}$ o\`{u}
$$
\mu (f)=\mathop{\rm dim}\nolimits_{k}(k[[x,y]]/
(\partial_{x}f,\partial_{y}f)),
$$
d'o\`{u} la conclusion.
\hfill\hfill$\square$
\vskip 3mm

La caract\'{e}ristique d'Euler-Poincar\'{e} $\chi (R\mathop{\rm
Hom}\nolimits_{{\cal O}_{C_{k}}}(\Omega_{C_{k}/k}^{1},{\cal
O}_{C_{k}}))$, c'est-\`{a}-dire la dif\-f\'{e}\-rence $\mathop{\rm
dim}_{k}\mathop{\rm Lie}(\mathop{\rm Aut}\nolimits_{k}(C_{k})) -\tau
(C_{k})$, a \'{e}t\'{e} calcul\'{e}e par Deligne (cf.  [De 2] Proposition
2.32): on a
$$
\tau (C_{k})=3p(C_{k})-3+\mathop{\rm dim}\nolimits_{k}
\mathop{\rm Lie} (\mathop{\rm Aut}\nolimits_{k}(C_{k})).\leqno{(3.2.4)}
$$

\rem Variante $3.2.5$
\endrem
Fixons un ensemble fini $\Sigma$ de $k$-points dans le lieu de non
lissit\'{e} de $C_{k}$ dont le cardinal est suffisamment grand pour
que $\mathop{\rm Hom}\nolimits_{k}(\Omega_{C_{k}/k}^{1}(\Sigma ),{\cal
O}_{C_{k}})=(0)$.  On peut remplacer $C/S$ par une alg\'{e}brisation
$(C_{\Sigma}/S_{\Sigma})$ de la d\'{e}formation formelle universelle
du couple $(C_{k},\Sigma )$ (il n'y a plus d'automorphismes
infinit\'{e}simaux et le foncteur des d\'{e}formations est donc
pro-repr\'{e}sentable).  La dimension de $S_{\Sigma}$ est
$3p(C_{k})-3+|\Sigma |$ et le morphisme naturel $S_{\Sigma}\rightarrow
S$ est formellement lisse de dimension relative $|\Sigma |-\mathop{\rm
dim}\nolimits_{k} \mathop{\rm Lie} (\mathop{\rm
Aut}\nolimits_{k}(C_{k}))$.
\hfill\hfill$\square$

\section{3.3}{Normalisation en famille et constance de l'invariant
$\delta$, d'apr\`{e}s Teissier}

Nous noterons dans la suite $\pi_{X}:\widetilde{X} \rightarrow X$
le morphisme de normalisation de tout sch\'{e}ma int\`{e}gre $X$.
\vskip 1mm

Pour toute courbe $C_{\kappa}$ g\'{e}om\'{e}triquement int\`{e}gre sur
un corps $\kappa$, on a la suite exacte
$$
0\rightarrow {\cal O}_{C_{\kappa}}\rightarrow\pi_{\kappa,\ast}{\cal
O}_{\widetilde{C_{\kappa}}}\rightarrow\bigoplus_{c\in C_{\kappa}^{{\rm
sing}}}(\pi_{\kappa,\ast}{\cal O}_{\widetilde{C_{\kappa}}}/{\cal
O}_{C_{\kappa}})_{c}\rightarrow 0
$$
o\`{u} $C_{\kappa}^{{\rm sing}}$ est le lieu singulier (fini) de
$C_{\kappa}$, et on pose
$$
\delta (C_{\kappa})=\sum_{c\in C_{\kappa}}\delta_{c}(C_{\kappa})
$$
o\`{u}, pour chaque $c\in C_{\kappa}^{{\rm sing}}$,
$\delta_{c}(C_{\kappa})$ est la dimension du $\kappa$-espace vectoriel
$(\pi_{\kappa,\ast} {\cal O}_{\widetilde{C}_{\kappa}}/{\cal
O}_{C_{\kappa}})_{c}$.
\vskip 2mm

Soient $S$ un sch\'{e}ma local complet, int\`{e}gre et
n{\oe}th\'{e}rien, de point ferm\'{e} $s$, $\varphi :C\rightarrow S$
une courbe relative plate \`{a} fibres g\'{e}om\'{e}triquement
r\'{e}duites et $D\subset C$ un diviseur effectif fini et plat sur $S$
tel $\varphi$ soit lisse en dehors du support de $D$.

Soit ${\frak a}\subset {\cal O}_{C}$ l'Id\'{e}al annulateur du ${\cal
O}_{C}$-Module coh\'{e}rent $\pi_{C,\ast}{\cal
O}_{\widetilde{C}}/{\cal O}_{C}$ et, pour chaque $t\in S$, soit
${\frak a}_{t}\subset {\cal O}_{C_{t}}$ l'Id\'{e}al annulateur du
${\cal O}_{C_{t}}$-Module coh\'{e}rent $\pi_{C_{t},\ast}{\cal
O}_{\widetilde{C_{t}}}/{\cal O}_{C_{t}}$, o\`{u} bien entendu $C_{t}$
est la fibre de $\varphi$ en $t$. Il n'est pas vrai en g\'{e}n\'{e}ral
que ${\frak a}_{t}$ soit la restriction de ${\frak a}$ \`{a} $C_{t}$.

\thm LEMME 3.3.1
\enonce
On peut trouver un diviseur effectif $D'\subset C$, fini et plat sur
$S$, tel que l'on ait les inclusions
$$
{\cal O}_{C_{t}}(-D_{t}')\subset {\frak a}_{t}\subset {\cal O}_{C_{t}}
$$
quel que soit $t\in S$.
\endthm

\rem Preuve
\endrem
Par hypoth\`{e}se, on a un diviseur effectif $D\subset C$ fini et plat
sur $S$ tel que $\varphi$ soit lisse en dehors du support de $D$.

Pour chaque point $t$ de $S$, $D$ induit un diviseur $D_{t}\subset
C_{t}$ et, comme $\pi_{C_{t}}$ est un isomorphisme en dehors du
support de $D_{t}$, le support du ferm\'{e} de $C_{t}$ d\'{e}fini par
${\frak a}_{t}$ est contenu dans celui de $D_{t}$. Il existe donc un
entier $n\geq 0$ tel que ${\cal O}_{C_{t}}(-nD_{t})\subset {\frak
a}_{t}\subset {\cal O}_{C_{t}}$; notons $n_{t}$ le plus petit entier
$n\geq 0$ ayant cette propri\'{e}t\'{e}.

La fonction $S\rightarrow {\Bbb N},~t\mapsto n_{t}$, est
constructible. En effet, par induction n{\oe}th\'{e}rienne, il suffit de
montrer que cette fonction est localement constante sur un ouvert
dense de $S$. Mais cela r\'{e}sulte de l'existence d'un ouvert dense
normal $U$ de $S$ tel que, pour tout $t\in U$, $\pi_{C_{t}}$ soit la
fibre en $t$ de $\pi_{C}$ et ${\frak a}_{t}$ soit la restriction de
${\frak a}$ \`{a} $C_{t}$.

On peut donc choisir un entier $n$ tel que $n\geq n_{t}$ quel que soit
$t\in S$ et le diviseur $D'=nD$ r\'{e}pond \`{a} la question.
\hfill\hfill$\square$
\vskip 3mm

\thm PROPOSITION 3.3.2 (Teissier, [Te 3] I, 1.3.2)
\enonce
La fonction $S\rightarrow {\Bbb Z},~ t\mapsto\delta (C_{t})$, est
semi-continue sup\'{e}rieurement. Elle est m\^{e}me constante si l'on
suppose de plus que le morphisme compos\'{e} $\widetilde{\varphi}=
\varphi\circ\pi_{C}:\widetilde{C} \rightarrow S$ est lisse.

Inversement, si $S$ est {\rm normal} et si la fonction $S\rightarrow
{\Bbb Z},~ t\mapsto\delta (C_{t})$, est constante sur $S$, le
morphisme $\widetilde{\varphi}$ est lisse.
\endthm

Chaque fois que le morphisme compos\'{e} $\varphi\circ\pi_{C}:
\widetilde{C}\rightarrow S$ sera lisse, on verra le morphisme de
normalisation $\pi_{C}:\widetilde{C}\rightarrow C$ comme une {\it
normalisation en famille} de la courbe relative $\varphi :C\rightarrow
S$.
\vskip 2mm

Teissier d\'{e}montre dans un premier temps le cas particulier plus
pr\'{e}cis suivant:

\thm LEMME 3.3.3
\enonce
Supposons de plus que $S=\mathop{\rm Spec}(V)$ soit un trait {\rm
(}$V$ est donc un anneau de valuation discr\`{e}te{\rm )}, de point
g\'{e}n\'{e}rique $\eta$. Alors, le morphisme compos\'{e}
$\widetilde{\varphi}=\varphi\circ\pi_{C}: \widetilde{C}\rightarrow S$
est plat et on a la relation
$$
\delta (C_{s})-\delta (C_{\eta})=\delta ((\widetilde{C})_{s})\geq 0
$$
o\`{u} $(\widetilde{C})_{s}$ est la fibre sp\'{e}ciale de
$\widetilde{\varphi}$.
\endthm

\rem Preuve du lemme $3.3.3$
\endrem
Il r\'{e}sulte du lemme $3.3.1$ que l'image r\'{e}ciproque de la
restriction de ${\frak a}$ \`{a} $C_{s}\subset C$ par le morphisme de
normalisation $\pi_{C_{s}}:\widetilde{C_{s}}\rightarrow C_{s}$ est non
nulle. La propri\'{e}t\'{e} universelle de la normalisation $\pi_{C}
:\widetilde{C}\rightarrow C$ nous donne alors une factorisation
$$
\widetilde{C_{s}}\longrightarrow\widetilde{C}\,\smash{\mathop{\hbox to
8mm{\rightarrowfill}} \limits^{\scriptstyle \pi_{C}}}\,C
$$
de $\widetilde{C_{s}}\,\smash{\mathop{\hbox to 8mm{\rightarrowfill}}
\limits^{\scriptstyle \pi_{C_{s}}}}\,C_{s}\subset C$, ou ce qui
revient au m\^{e}me une factorisation
$$
\widetilde{C_{s}}\,\smash{\mathop{\hbox to 8mm{\rightarrowfill}}
\limits^{\scriptstyle \rho}}\,(\widetilde{C})_{s}\,\smash{\mathop{\hbox to
8mm{\rightarrowfill}} \limits^{\scriptstyle \pi_{C,s}}}\,C_{s}
$$
de $\pi_{C_{s}}$ par la fibre sp\'{e}ciale $\pi_{C,s}$ de $\pi_{C}$.

Fixons une uniformisante $v$ de $V$. On v\'{e}rifie que
$$
(v\cdot \pi_{C,\ast}{\cal O}_{\widetilde{C}})\cap {\cal O}_{C}=v\cdot
{\cal O}_{C},
$$
et donc que le ${\cal O}_{S}$-Module coh\'{e}rent $\pi_{C,\ast}{\cal
O}_{\widetilde{C}}/{\cal O}_{C}$ est sans $v$-torsion et que les
homomorphismes de $k$-Alg\`{e}bres
$$
{\cal O}_{C_{s}}\rightarrow (\pi_{C,s})_{\ast}{\cal
O}_{(\widetilde{C})_{s}}\rightarrow \pi_{C_{s},\ast}{\cal
O}_{\widetilde{C_{s}}}
$$
correspondant \`{a} la factorisation ci-dessus sont injectifs. On
v\'{e}rifie aussi que la fibre g\'{e}n\'{e}rique $\pi_{C,\eta}:
(\widetilde{C})_{\eta} \rightarrow C_{\eta}$ de $\pi_{C}$ est la
normalisation de $C_{\eta}$. Par suite, le ${\cal O}_{S}$-Module
$\varphi_{\ast}(\pi_{C,\ast}{\cal O}_{\widetilde{C}}/{\cal O}_{C})$
est libre de rang fini \'{e}gal \`{a}
$$\eqalign{
\mathop{\rm dim}\nolimits_{\kappa (s)}\varphi_{s,\ast}
((\pi_{C,s})_{\ast}{\cal O}_{(\widetilde{C})_{s}}/{\cal
O}_{C_{s}})&=\mathop{\rm dim}\nolimits_{\kappa (s)}
(\varphi_{\ast}(\pi_{C,\ast}{\cal O}_{\widetilde{C}}/{\cal
O}_{C}))_{s}\cr
&=\mathop{\rm dim}\nolimits_{\kappa (\eta )}
(\varphi_{\ast}(\pi_{C,\ast}{\cal O}_{\widetilde{C}}/{\cal
O}_{C}))_{\eta}\cr
&=\mathop{\rm dim}\nolimits_{\kappa (\eta )}\varphi_{\eta ,\ast}
((\pi_{C,\eta})_{\ast}{\cal O}_{(\widetilde{C})_{\eta}}/{\cal
O}_{C_{\eta}})=\delta (C_{\eta}),\cr}
$$
le morphisme $\widetilde{\varphi}$ est plat, sa fibre
sp\'{e}ciale $(\widetilde{C})_{s}$ est int\`{e}gre, $\rho$ est le
morphisme de normalisation de $(\widetilde{C})_{s}$ et on a bien
$$
\delta (C_{s})=\delta ((\widetilde{C})_{s})+\mathop{\rm dim}
\nolimits_{k}\varphi_{s,\ast}((\pi_{C,s})_{\ast}{\cal
O}_{(\widetilde{C})_{s}}/{\cal O}_{C_{s}})=\delta
((\widetilde{C})_{s})+\delta (C_{\eta}).
$$
\hfill\hfill$\square$
\vskip 3mm

\rem Preuve de la proposition $3.3.2$
\endrem
La semi-continuit\'{e} de la fonction $t\mapsto \delta (C_{t})$ dans
le cas g\'{e}n\'{e}ral r\'{e}sulte de sa constructibilit\'{e}
laiss\'{e}e au lecteur (voir la preuve du lemme $3.3.1$) et du cas
d\'{e}j\`{a} trait\'{e} o\`{u} $S$ est un trait.
\vskip 1mm

Si $\widetilde{\varphi}$ est lisse, il r\'{e}sulte du lemme $3.3.3$ que
$$
\delta (C_{s})-\delta (C_{\eta})=\delta ((\widetilde{C})_{s})=0
$$
dans le cas o\`{u} $S$ est un trait. On en d\'{e}duit
imm\'{e}diatement que la fonction $t\mapsto\delta (C_{t})$ est
constante pour $S$ arbitraire.
\vskip 1mm

Inversement, supposons que $t\mapsto \delta (C_{t})$ est constante de
valeur $\delta$ et fixons un diviseur effectif $D'\subset S$ comme dans
le lemme $3.3.1$. Le ${\cal O}_{S}$-Module $\varphi_{\ast}({\cal
O}_{C}(D')/{\cal O}_{C})$ est localement libre de rang fini
$d\geq\delta$. Soit $G\rightarrow S$ la grassmannienne des quotients
localement libres de rang $d-\delta$ de ce ${\cal O}_{S}$-Module. On
d\'{e}finit une section {\it ensembliste} $\sigma :S\rightarrow G$ de
$G$ sur $S$ en envoyant le point $t$ sur le conoyau de l'inclusion
$$
\Gamma (C_{t},\pi_{C_{t},\ast}{\cal O}_{\widetilde{C_{t}}}/{\cal
O}_{C_{t}})\subset\Gamma (C_{t},{\cal O}_{C_{t}}(D_{t}')/{\cal O}_{C_{t}})
$$
(on a ${\cal O}_{C_{t}}(-D_{t}')\cdot \pi_{C_{t},\ast}{\cal
O}_{\widetilde{C_{t}}}\subset {\frak a}_{t}\cdot \pi_{C_{t},\ast}{\cal
O}_{\widetilde{C_{t}}}\subset {\cal O}_{C_{t}}$ et donc
$\pi_{C_{t},\ast}{\cal O}_{\widetilde{C_{t}}}\subset {\cal
O}_{C_{t}}(D_{t}')$).

Soit $Z\subset G$ l'ensemble des points $\sigma (t)$ pour $t$
parcourant $S$. On v\'{e}rifie par induction n{\oe}th\'{e}rienne
(comme dans la preuve du lemme $3.3.1$) que l'ensemble $Z$ est
constructible. On v\'{e}rifie aussi que les constructions de ${\cal
F}$, $G$ et $Z$ commutent aux changements de bases $S'\rightarrow S$.
On v\'{e}rifie enfin, \`{a} l'aide du lemme $3.3.3$, que dans le cas
o\`{u} $S$ est un trait, $Z$ est l'image ensembliste d'une section
alg\'{e}brique de $G$ sur $S$.

On d\'{e}duit de ces propri\'{e}t\'{e}s que $Z$ est l'ensemble des
points d'un ferm\'{e} r\'{e}duit de $G$, not\'{e} encore $Z$, et que
la restriction de la projection canonique $G\rightarrow S$ \`{a} ce
ferm\'{e} est un hom\'{e}omorphisme de $Z$ sur $S$.

Supposons de plus que $S$ soit normal. Cet hom\'{e}omorphisme est
alors n\'{e}cessairement un isomorphisme et $\sigma$ est en fait une
section alg\'{e}brique de $G\rightarrow S$. En d'autres termes, les
espaces vectoriels $\Gamma (C_{t},\pi_{C_{t},\ast}{\cal
O}_{\widetilde{C_{t}}}/{\cal O}_{C_{t}})$ pour $t\in S$ sont les
fibres fibr\'{e} vectoriel ${\cal E}$ de rang $\delta$ sur $S$ qui est
un sous-${\cal O}_{S}$-Module localement facteur direct de
$\varphi_{\ast}({\cal O}_{C}(D')/{\cal O}_{C})$.

Consid\'{e}rons maintenant le ${\cal O}_{C}$-Module coh\'{e}rent
$$
{\cal M}=({\cal O}_{C}\oplus {\cal O}_{C}(D'))/{\cal O}_{C}
$$
o\`{u} ${\cal O}_{C}$ est plong\'{e} diagonalement dans ${\cal
O}_{C}\oplus {\cal O}_{C}(D')$. On a d'une part une suite exacte
\'{e}vidente
$$
0\rightarrow {\cal O}_{C}\rightarrow {\cal M}\rightarrow {\cal
O}_{C}(D')/{\cal O}_{C}\rightarrow 0.
$$
On a d'autre part le plongement
$$
{\cal M}\hookrightarrow {\cal O}_{C}(D'),~(f\oplus f')+{\cal
O}_{C}\mapsto f-f',
$$
que l'on peut composer avec le plongement canonique de ${\cal
O}_{C}(D')$ dans l'Anneau total des fractions ${\cal K}_{C}$ de ${\cal
O}_{C}$.

Notons ${\cal F}\subset {\cal M}$ l'image r\'{e}ciproque de ${\cal E}$
par la surjection ${\cal M}\twoheadrightarrow {\cal O}_{C}(D')/{\cal
O}_{C}$. On peut voir ${\cal F}\subset {\cal M}$ comme un
sous-${\cal O}_{C}$-Module de ${\cal K}_{C}$. Comme ${\cal
O}_{C_{t}}(-D_{t}')\subset {\frak a}_{t}$, on v\'{e}rifie facilement
que la restriction de ${\cal F}$ \`{a} $C_{t}$ est une sous-${\cal
O}_{C_{t}}$-Alg\`{e}bre de l'Anneau total des fractions ${\cal
K}_{C_{t}}$ de ${\cal O}_{C_{t}}$ quel que soit $t\in S$. Il s'en suit
que ${\cal F}$ est une sous-${\cal O}_{C}$-Alg\`{e}bre de ${\cal
K}_{C}$, et donc que
$$
{\cal O}_{C}\subset {\cal F}\subset\pi_{C,\ast}{\cal
O}_{\widetilde{C}}\subset {\cal K}_{C}
$$
puisque $\pi_{C,\ast}{\cal O}_{\widetilde{C}}$ est la cl\^{o}ture
int\'{e}grale de ${\cal O}_{C}$ dans ${\cal K}_{C}$ et que ${\cal F}$
est \'{e}videmment un ${\cal O}_{C}$-Module coh\'{e}rent.

Comme on l'a vu au cours de la preuve du lemme $3.3.3$, pour chaque
$t\in S$, on a les inclusions
$$
{\cal O}_{C_{t}}\subset (\pi_{C,t})_{\ast}{\cal
O}_{(\widetilde{C})_{t}}\subset\pi_{C_{t},\ast}{\cal
O}_{\widetilde{C_{t}}}\subset {\cal K}_{C_{t}}
$$
et donc $(\pi_{C,t})_{\ast}{\cal O}_{(\widetilde{C})_{t}}$ est contenu
dans la restriction de ${\cal F}$ \`{a} $C_{t}$. Il s'en suit que
$\pi_{C,\ast}{\cal O}_{\widetilde{C}}\subset {\cal F}\subset {\cal
K}_{C}$ et donc que ${\cal F}=\pi_{C,\ast}{\cal O}_{\widetilde{C}}$.

Comme ${\cal O}_{C}$ et ${\cal O}_{C}(D')/{\cal O}_{C}$ sont $S$-plats,
il en est de m\^{e}me de ${\cal M}$. De plus, comme ${\cal E}$ est
aussi $S$-plat, il en est de m\^{e}me de ${\cal F}=\pi_{C,\ast}{\cal
O}_{\widetilde{C}}$.  On a donc bien montr\'{e} que $\widetilde{\varphi}$
est plat et qu'en outre
$$
(\widetilde{C})_{t}=\widetilde{C_{t}}
$$
pour tout $t\in S$.
\hfill\hfill$\square$

\thm COROLLAIRE 3.3.4
\enonce
Les points $t\in S$ pour lesquels $\delta (C_{t})=\delta (C_{s})$ sont
les points points d'un ferm\'{e} r\'{e}duit $S^{\delta}\subset S$,
{\rm la strate \`{a} $\delta$-constant}.

La courbe plate relative $C_{\widetilde{S^{\delta}}}=
\widetilde{S^{\delta}}\times_{S}C \rightarrow\widetilde{S^{\delta}}$
d\'{e}duite de $\varphi :C\rightarrow S$ par le changement de base par
le morphisme
$$
\widetilde{S^{\delta}}\twoheadrightarrow S^{\delta}\hookrightarrow S,
$$
compos\'{e} du morphisme de normalisation $\pi_{S^{\delta}}$ de
$S^{\delta}$ et de l'inclusion de $S^{\delta}$ dans $S$, admet une
normalisation en famille $\pi_{C_{\widetilde{S^{\delta}}}}:
\widetilde{C_{\widetilde{S^{\delta}}}}\rightarrow
C_{\widetilde{S^{\delta}}}$.
\hfill\hfill$\square$
\endthm

\section{3.4}{La strate \`{a} $\delta$ constant}

Soit $A=k[[x,y]]/(f)$ l'anneau formel d'un germe de courbe plane \`{a}
singularit\'{e} isol\'{e}e (le $k$-espace vectoriel
$k[[x,y]]/(f,\partial_{x}f,\partial_{y}f)$ est donc suppos\'{e} de
dimension finie).

On note $K$ l'anneau total des fractions de $A$ et
$\widetilde{A}\subset K$ la normalisation de $A$ dans $K$. On peut
d\'{e}composer $K$ en un produit fini de corps $K=\prod_{i\in I}K_{i}$
et $\widetilde{A}$ en le produit correspondant d'anneaux int\`{e}gres
$\widetilde{A}=\prod_{i\in I}\widetilde{A}_{i}$. On pose $\delta
(A)=\mathop{\rm dim}\nolimits_{k}(\widetilde{A}/A)$.

On fixe dans la suite une uniformisante $t_{i}$ de $\widetilde{A}_{i}$
pour chaque $i\in I$, de sorte que $\widetilde{A}_{i}=k[[t_{i}]]$. Le
plongement $A\hookrightarrow\widetilde{A}$ est donn\'{e} par une
famille de couples $(x_{i}(t_{i}),y_{i}(t_{i}))\in k[[t_{i}]]^{2}$
index\'{e}e par $i\in I$.

On rappelle (cf.  [A-K 1], Chapter 8) que:
\vskip 1mm

\itemitem{-} le module dualisant $\omega_{A}$ est le $A$-module libre
de rang $1$ d\'{e}fini par
$$
\omega_{A}=\mathop{\rm Ext}\nolimits_{k[[x,y]]}^{1}(A,
\omega_{k[[x,y]]})
$$
o\`{u} $\omega_{k[[x,y]]}=\Omega_{k[[x,y]]/k}^{2}$ est un
$k[[x,y]]$-module libre de rang $1$,
\vskip 1mm

\itemitem{-} pour tout $A$-module $M$ et tout entier $i$, on a un
isomorphisme canonique de $A$-modules
$$
\mathop{\rm Ext}\nolimits_{A}^{i}(M,\omega_{A})\cong
\mathop{\rm Ext}\nolimits_{k[[x,y]]}^{i+1}(M,\omega_{k[[x,y]]})
$$
o\`{u} $M$ est vu comme un $k[[x,y]]$-module via l'\'{e}pimorphisme
canonique $k[[x,y]]\twoheadrightarrow A$,
\vskip 1mm

\itemitem{-} le module dualisant $\omega_{A}$ est donn\'{e}
concr\`{e}tement par
$$
\omega_{A}=A\textstyle{{\rm d}x\wedge {\rm d}y\over {\rm d}f}=A{{\rm
d}x\over\partial_{y}f}= A\bigl(-{{\rm d}y\over\partial_{x}f}\bigr)
\subset\Omega_{K/k}^{1},
$$
et aussi par
$$
\omega_{A}=\{\alpha\in\Omega_{K/k}^{1}\mid \mathop{\rm Res}(A\alpha
)=(0)\}
$$
o\`{u}
$$
\mathop{\rm Res}=\sum_{i\in I}\mathop{\rm
Res}\nolimits_{i}:\Omega_{K/k}^{1}=\bigoplus_{i\in
I}\Omega_{K_{i}/k}^{1}\rightarrow k,~\oplus_{i\in I}a_{i}(t_{i}){{\rm
d}t_{i}\over t_{i}}\rightarrow\sum_{i\in I}a_{i}(0),
$$
est la somme des homomorphismes r\'{e}sidus,
\vskip 1mm

\itemitem{-} on a
$$
\omega_{\widetilde{A}}=\Omega_{\widetilde{A}/k}^{1}\subset\omega_{A}
$$
et l'accouplement
$$
(\widetilde{A}/A)\times (\omega_{A}/\omega_{\widetilde{A}})\rightarrow
k,~(\widetilde{a}+A,\alpha +\omega_{\widetilde{A}})\rightarrow
\mathop{\rm Res}(\widetilde{a}\alpha ),
$$
est un accouplement parfait entre deux $k$-espaces vectoriels de
dimension $\delta (A)$,
\vskip 1mm

\itemitem{-} le conducteur
$$
{\frak a}=\{a\in A\mid a\widetilde{A}\subset A\}
$$
est aussi le conducteur
$$
{\frak a}=\{a\in A\mid a\omega_{A}\subset\omega_{\widetilde{A}}\}.
$$
\vskip 1mm

En particulier, on a:

\thm LEMME 3.4.1
\enonce
L'id\'{e}al de $A$ engendr\'{e} par les classes de $\partial_{x}f$ et
$\partial_{y}f$ modulo $(f)$ est contenu dans le conducteur ${\frak
a}$.
\hfill\hfill$\square$
\endthm

On consid\`{e}re maintenant le foncteur de d\'{e}formations
$$
\mathop{\rm Def}\nolimits_{A}^{{\rm top}}:\mathop{\rm Art}
\nolimits_{k}\rightarrow\mathop{\rm Ens}
$$
(cf.  3.1) et on s'int\'{e}resse plus particuli\`{e}rement aux
d\'{e}formations de $A$ \`{a} $\delta$ constant.  Pour cela, on va
\'{e}tudier les d\'{e}formations de l'homomorphisme de
$k$-alg\`{e}bres $A\hookrightarrow \widetilde{A}$.  On consid\`{e}re
donc le foncteur
$$
\mathop{\rm Def}\nolimits_{A\hookrightarrow \widetilde{A}}^{{\rm
top}}:\mathop{\rm Art} \nolimits_{k}\rightarrow\mathop{\rm Ens}
$$
qui envoie $R\in\mathop{\rm ob}\mathop{\rm Art}\nolimits_{k}$ sur
l'ensemble des classes d'isomorphie d'homomorphismes de
$R$-alg\`{e}bres
$$
A_{R}\rightarrow\widetilde{A}_{R}
$$
dont la r\'{e}duction modulo ${\frak m}_{R}$ est l'inclusion $A\subset
\widetilde{A}$, o\`{u} $A_{R}$ est une d\'{e}formation plate sur $R$
de $A$ et $\widetilde{A}_{R}$ est une d\'{e}formation plate sur $R$ de
$\widetilde{A}$. On a bien s\^{u}r un morphisme d'oubli
$$
\mathop{\rm Def}\nolimits_{A\hookrightarrow\widetilde{A}}^{{\rm top}}
\rightarrow \mathop{\rm Def}\nolimits_{A}^{{\rm top}}.
$$

\thm LEMME 3.4.2
\enonce
Tout homomorphisme de $R$-alg\`{e}bres $(A_{R}\rightarrow
\widetilde{A}_{R})\in \mathop{\rm Def}\nolimits_{A\hookrightarrow
\widetilde{A}}^{{\rm top}}(R)$ est n\'{e}cessairement injectif et son
conoyau est automatiquement $R$-plat.
\endthm

\rem Preuve
\endrem
Notons $N_{R}$, $I_{R}$ et $C_{R}$ les noyau, image et conoyau de
l'homomorphisme $A_{R}\rightarrow\widetilde{A}_{R}$. On a les suites
exactes
$$
0\rightarrow \mathop{\rm Tor}\nolimits_{1}^{R}(C_{R},k)\rightarrow
I_{R}\otimes_{R}k\rightarrow \widetilde{A}_{R}\otimes_{R}k\rightarrow
C_{R}\otimes_{R}k\rightarrow 0
$$
et
$$
0\rightarrow \mathop{\rm Tor}\nolimits_{1}^{R}(I_{R},k)\rightarrow
N_{R}\otimes_{R}k\rightarrow A_{R}\otimes_{R}k\rightarrow
I_{R}\otimes_{R}k\rightarrow 0
$$
et l'homomorphisme compos\'{e}
$$
A_{R}\otimes_{R}k\twoheadrightarrow I_{R}\otimes_{R}k\rightarrow
\widetilde{A}_{R}\otimes_{R}k
$$
est par hypoth\`{e}se l'inclusion $A\subset\widetilde{A}$; par suite,
la surjection $A_{R}\otimes_{R}k\twoheadrightarrow I_{R}\otimes_{R}k$
est n\'{e}cessairement un isomorphisme, la fl\`{e}che $I_{R}\otimes_{R}k
\rightarrow \widetilde{A}_{R}\otimes_{R}k$ est injective et
$\mathop{\rm Tor}\nolimits_{1}^{R}(C_{R},k)=(0)$, de sorte que $C_{R}$
est bien $R$-plat; mais alors $I_{R}$ est aussi $R$-plat puisque
$\widetilde{A}_{R}$ l'est, et on a $N_{R}\otimes_{R}k=(0)$. Il s'en
suit que $N_{R}=(0)$ et le lemme est d\'{e}montr\'{e}.
\hfill\hfill$\square$
\vskip 3mm

\thm LEMME 3.4.3
\enonce
Soit $M$ un $A$-module de type fini sans torsion. On a
$$
\mathop{\rm Ext}\nolimits_{A}^{i}(M,A)=(0),~\forall i\not=0,
$$
et il existe un entier $n\geq 0$ tel que $M$, vu comme
$k[[x,y]]$-module via l'\'{e}pimorphisme canonique
$k[[x,y]]\twoheadrightarrow A$, admette une r\'{e}solution
$$
0\rightarrow k[[x,y]]^{n}\rightarrow k[[x,y]]^{n}\rightarrow
M\rightarrow 0.
$$
Si on suppose de plus que $M$ est de rang g\'{e}n\'{e}rique $1$, il
existe m\^{e}me un tel entier $n\leq\delta (A)+1$.
\endthm

\rem Preuve
\endrem
Comme le $A$-module $\omega_{A}$ est libre de rang $1$, pour tout
$A$-module $M$ on a
$$
\mathop{\rm Ext}\nolimits_{A}^{i}(M,A)\cong \mathop{\rm
Ext}\nolimits_{A}^{i}(M,\omega_{A})\cong \mathop{\rm
Ext}\nolimits_{k[[x,y]]}^{i+1}(M,\omega_{k[[x,y]]})
$$
et, comme la $k$-alg\`{e}bre $k[[x,y]]$ est r\'{e}guli\`{e}re de
dimension $2$, il s'en suit que
$$
\mathop{\rm Ext}\nolimits_{A}^{i}(M,A)=(0)
$$
quel que soit $i\not=0,1$.  Si $M$ est sans torsion, on a de plus la
suite exacte
$$
\mathop{\rm Ext}\nolimits_{A}^{1}(K\otimes_{A}M,A)\rightarrow
\mathop{\rm Ext}\nolimits_{A}^{1}(M,A)\rightarrow \mathop{\rm
Ext}\nolimits_{A}^{2}(K\otimes_{A}M/M,A)
$$
o\`{u}
$$
\mathop{\rm Ext}\nolimits_{A}^{1}(K\otimes_{A}M,A)\cong \mathop{\rm
Ext}\nolimits_{K}^{1}(K\otimes_{A}M,K)=(0)
$$
et $\mathop{\rm Ext}\nolimits_{A}^{2}(K\otimes_{A}M/M,A)=(0)$, et donc
on a aussi $\mathop{\rm Ext}\nolimits_{A}^{1}(M,A)=(0)$ et a fortiori
$$
\mathop{\rm Ext}\nolimits_{k[[x,y]]}^{2}(M,k[[x,y]])\cong
\mathop{\rm Ext}\nolimits_{k[[x,y]]}^{2}(M,\omega_{k[[x,y]]})=(0).
$$
Par suite, si $M$ est de type fini et sans torsion, le noyau de tout
\'{e}pimorphisme de $k[[x,y]]$-modules $k[[x,y]]^{n}\twoheadrightarrow
M$ est n\'{e}cessairement libre de rang fini et donc non canoniquement
isomorphe \`{a} $k[[x,y]]^{n}$ puisqu'en tant que $k[[x,y]]$-module,
$M$ est de rang g\'{e}n\'{e}rique $0$.

Montrons enfin que, si $M$ un $A$-module de type fini, sans torsion et
de rang g\'{e}n\'{e}rique $1$, $M$ peut \^{e}tre engendr\'{e} par
$\delta (A)+1$ \'{e}l\'{e}ments.  Comme $\widetilde{A}\otimes_{A}M$
est un $\widetilde{A}$-module libre de rang $1$ et que l'homomorphisme
canonique $M\rightarrow \widetilde{A}\otimes_{A}M$ est injectif, on
peut supposer $M\subset\widetilde{A}$ et que $\widetilde{A}M=
\widetilde{A}$.  Mais alors, $M$ contient au moins un \'{e}l\'{e}ment
inversible $m$ de $\widetilde{A}$ et, \`{a} isomorphisme pr\`{e}s, on
peut supposer que $m=1$, c'est-\`{a}-dire que
$$
A\subset M\subset\widetilde{A}.
$$
Sous ces conditions, on a $\mathop{\rm
dim}\nolimits_{k}(M/A)\leq\delta (A)$ et on conclut en remarquant que
$M$ est engendr\'{e} sur $A$ par $\{1,m_{1},\ldots ,m_{n}\}$ o\`{u}
$\{m_{1},\ldots ,m_{n}\}\subset M$ repr\'{e}sente une base de $M/A$
sur $k$.
\hfill\hfill$\square$
\vskip 3mm

Consid\'{e}rons maintenant les foncteurs
$$
\mathop{\rm Def}\nolimits_{k[[x,y]]\twoheadrightarrow A}^{{\rm
top}},~\mathop{\rm Def}\nolimits_{k[[x,y]]
\rightarrow\widetilde{A}}^{{\rm top}}:\mathop{\rm Art}
\nolimits_{k}\rightarrow\mathop{\rm Ens}
$$
qui envoie $R\in\mathop{\rm ob}\mathop{\rm Art}\nolimits_{k}$ sur
l'ensemble des classes d'isomorphie d'homomorphismes de
$R$-alg\`{e}bres
$$
R[[x,y]]\rightarrow A_{R}\hbox{ et }R[[x,y]]\rightarrow \widetilde{A}_{R}
$$
dont les r\'{e}ductions modulo ${\frak m}_{R}$ sont l'\'{e}pimorphisme
canonique $k[[x,y]]\twoheadrightarrow A$ et l'homomorphisme
compos\'{e} $k[[x,y]]\twoheadrightarrow A\hookrightarrow
\widetilde{A}$, o\`{u} bien entendu $A_{R}$ et $\widetilde{A}_{R}$
sont des d\'{e}formations plates sur $R$ de $A$ et $\widetilde{A}$
respectivement. On remarque que $R[[x,y]]\rightarrow A_{R}$ est
automatiquement un \'{e}pimorphisme d'apr\`{e}s le lemme de Nakayama.
On a les morphismes de foncteurs \'{e}vidents
$$\diagram{
&&\mathop{\rm
Def}\nolimits_{A\hookrightarrow\widetilde{A}}^{{\rm top}}\cr
&&\llap{}\left\downarrow\vbox to 4mm{}\right.\rlap{$\scriptstyle
$}\cr
\mathop{\rm Def}\nolimits_{k[[x,y]]\twoheadrightarrow A}^{{\rm top}}& \kern
-2mm\smash{\mathop{\hbox to 8mm{\rightarrowfill}}
\limits_{\scriptstyle }}\kern -4mm&\mathop{\rm Def}\nolimits_{A}^{{\rm
top}}\cr}
$$
et on pose
$$
\mathop{\rm Def}\nolimits_{k[[x,y]]\twoheadrightarrow A
\hookrightarrow\widetilde{A}}^{{\rm top}}=\mathop{\rm Def}
\nolimits_{k[[x,y]] \twoheadrightarrow A}^{{\rm top}}
\times_{\mathop{\rm Def}\nolimits_{A}^{{\rm top}}}\mathop{\rm Def}
\nolimits_{A\hookrightarrow \widetilde{A}}^{{\rm top}}.
$$
Le foncteur $\mathop{\rm Def}\nolimits_{A}^{{\rm top}}$, le morphisme
de foncteurs $\mathop{\rm Def} \nolimits_{k[[x,y]]\twoheadrightarrow
A}^{{\rm top}}\rightarrow \mathop{\rm Def}\nolimits_{A}^{{\rm top}}$
et donc aussi le foncteur $\mathop{\rm Def}\nolimits_{k[[x,y]]
\twoheadrightarrow A}^{{\rm top}}$, sont formellement lisses.

\thm TH\'{E}OR\`{E}ME 3.4.4
\enonce
Le morphisme de composition $\mathop{\rm Def}\nolimits_{k[[x,y]]
\twoheadrightarrow A\hookrightarrow\widetilde{A}}^{{\rm
top}}\rightarrow \mathop{\rm Def}\nolimits_{k[[x,y]]
\rightarrow\widetilde{A}}^{{\rm top}}$ est un isomorphisme.
\endthm

\rem Preuve
\endrem
Nous allons expliciter un inverse pour ce morphisme de composition.
Suivant une suggestion de J.-B. Bost, nous utiliserons pour cela un
cas tr\`{e}s simple de la construction {\og}{$\mathop{\rm Div}$}{\fg}
de Mumford (cf.  [M-F] Chapter 5, \S 3).

Soit $R[[x,y]]\rightarrow\widetilde{A}_{R}$ une d\'{e}formation de
$k[[x,y]]\rightarrow\widetilde{A}$. On voit $\widetilde{A}$ et
$\widetilde{A}_{R}$ comme des modules sur $k[[x,y]]$ et $R[[x,y]]$
\`{a} l'aide de ces homomorphismes. Fixons arbitrairement une
r\'{e}solution
$$
0\rightarrow k[[x,y]]^{n}\,\smash{\mathop{\hbox to
8mm{\rightarrowfill}}\limits^{\scriptstyle F}}\,
k[[x,y]]^{n}\rightarrow\widetilde{A}\rightarrow 0
$$
de $\widetilde{A}$ en tant que $k[[x,y]]$-module (il en existe
d'apr\`{e}s le lemme pr\'{e}c\'{e}dent).  Si $F^{\ast}$ est la
matrice des co-facteurs de $F$, on a
$F^{\ast}F=FF^{\ast}=\det F$ et $\det F$
annule $\widetilde{A}$.  Inversement si une fonction $g\in k[[x,y]]$
annule $\widetilde{A}$, c'est-\`{a}-dire est telle que
$gk[[x,y]]^{n}\subset Fk[[x,y]]^{n}$, il existe une matrice
carr\'{e}e $G$ de taille $n\times n$ \`{a} coefficient dans
$k[[x,y]]^{n}$ telle que $g=FG$ et $g^{n}=\det F\cdot\det G$.
Comme $k[[x,y]]$ est r\'{e}duit, il s'en suit que $\det F$ et $f$
engendre le m\^{e}me id\'{e}al de $k[[x,y]]$ et que l'on peut demander
de plus que $\det F=f$.

Comme $\widetilde{A}_{R}$ est plat sur $R$, on peut relever la
r\'{e}solution ci-dessus en une r\'{e}solution
$$
0\rightarrow R[[x,y]]^{n}\,\smash{\mathop{\hbox to
8mm{\rightarrowfill}}\limits^{\scriptstyle F_{R}}}\,R[[x,y]]^{n}
\rightarrow\widetilde{A}_{R}\rightarrow 0.
$$
En consid\'{e}rant la matrice des co-facteurs de $F_{R}$, on montre
comme ci-dessus que $f_{R}=\mathop{\rm det}F_{R}$ dans le noyau de
l'homomorphisme $R[[x,y]]\rightarrow \widetilde{A}_{R}$, ou ce qui
revient au m\^{e}me que l'on a une factorisation
$$
R[[x,y]]\twoheadrightarrow R[[x,y]]/(f_{R})\rightarrow
\widetilde{A}_{R}
$$
qui rel\`{e}ve la factorisation $k[[x,y]]\twoheadrightarrow A
\hookrightarrow \widetilde{A}$. Alors, $A_{R}:=R[[x,y]]/(f_{R})$ est
un rel\`{e}vement plat sur $R$ de $A$ et la fl\`{e}che
$A_{R}\rightarrow \widetilde{A}_{R}$ est n\'{e}cessairement injective
(\`{a} conoyau $R$-plat) d'apr\`{e}s le lemme 3.4.2. La factorisation
$R[[x,y]]\twoheadrightarrow A_{R}\hookrightarrow\widetilde{A}_{R}$
est donc la factorisation canonique par l'image et la construction de
Mumford produit bien un inverse au morphisme de foncteurs $\mathop{\rm
Def}\nolimits_{k[[x,y]] \twoheadrightarrow A\hookrightarrow
\widetilde{A}}^{{\rm top}}\rightarrow \mathop{\rm
Def}\nolimits_{k[[x,y]] \rightarrow\widetilde{A}}^{{\rm top}}$.
\hfill\hfill$\square$
\vskip 3mm

Si $J$ est un id\'{e}al de carr\'{e} nul dans $R\in \mathop{\rm ob}
\mathop{\rm Art}\nolimits_{k}$ et si $\overline{R}=R/J$, il n'y a pas
d'obstruction \`{a} relever un homomorphisme de
$\overline{R}$-alg\`{e}bres $\overline{R}[[x,y]]\rightarrow
\widetilde{A}_{\overline{R}}$ en un homomorphisme de
$\overline{R}$-alg\`{e}bres $R[[x,y]]\rightarrow\widetilde{A}_{R}$.
En particulier, le foncteur $\mathop{\rm Def}\nolimits_{k[[x,y]]
\rightarrow \widetilde{A}}^{{\rm top}}$ est formellement lisse et son
espace tangent est le $k$-espace vectoriel $\widetilde{A}\oplus
\widetilde{A}$.  Le th\'{e}or\`{e}me admet donc le corollaire suivant:

\thm COROLLAIRE 3.4.5
\enonce
Le foncteur $\mathop{\rm Def}\nolimits_{k[[x,y]]\rightarrow
A\hookrightarrow \widetilde{A}}^{{\rm top}}\,$, et donc aussi le
foncteur $\mathop{\rm Def}\nolimits_{A\hookrightarrow
\widetilde{A}}^{{\rm top}}\,$, est formellement lisse.

\hfill\hfill$\square$
\endthm

\rem Remarque $3.4.6$
\endrem
En \'{e}valuant l'isomorphisme du th\'{e}or\`{e}me sur
$R=k[\varepsilon ]$ avec $\varepsilon^{2}=0$, on trouve que, pour tout
$(\dot{x}(t_{i}), \dot{y}(t_{i}))_{i\in I}\in\prod_{i\in
I}(k[[t_{i}]]\times k[[t_{i}]])=\widetilde{A}\oplus\widetilde{A}$, il
existe $\dot{f}(x,y)\in k[[x,y]]$ tel que
$$
(f+\varepsilon \dot{f})(x(t_{i})+\varepsilon \dot{x}(t_{i}),
y(t_{i})+\varepsilon \dot{y}(t_{i}))\equiv 0,~\forall i\in I,
$$
c'est-\`{a}-dire tel que
$$
\dot{f}(x(t_{i}),y(t_{i}))+\partial_{x}f(x(t_{i}),
y(t_{i}))\dot{x}(t_{i}) +\partial_{y}f(x(t_{i}),
y(t_{i}))\dot{y}(t_{i})\equiv 0,~\forall i\in I,
$$
ce qui est une reformulation du lemme 3.4.1.
\hfill\hfill$\square$
\vskip 3mm

Nous allons maintenant d\'{e}terminer l'espace tangent au foncteur
$\mathop{\rm Def}\nolimits_{A\hookrightarrow \widetilde{A}}^{{\rm
top}}$.  Plus g\'{e}n\'{e}ralement \'{e}tudions le probl\`{e}me de
rel\`{e}vement d'une d\'{e}formation $A_{\overline{R}}\rightarrow
\widetilde{A}_{\overline{R}}$ de $A\hookrightarrow\widetilde{A}$ sur
$\overline{R}=R/J$, o\`{u} $J$ est un id\'{e}al de carr\'{e} nul dans
$R$, en une d\'{e}formation $A_{R}\rightarrow \widetilde{A}_{R}$ de
$A\hookrightarrow\widetilde{A}$ sur $R$.

Commen\c{c}ons par fixer un rel\`{e}vement $R$-plat
$\widetilde{A}_{R}$ de $\widetilde{A}_{\overline{R}}$ sur $R$.  On a
donc un diagramme
$$\diagram{
R&\twoheadrightarrow&\overline{R}\cr
&&\llap{$\scriptstyle $}\left\downarrow
\vbox to 4mm{}\right.\rlap{}\cr
&&A_{\overline{R}}\cr
&&\llap{$\scriptstyle $}\left\downarrow
\vbox to 4mm{}\right.\rlap{}\cr
\widetilde{A}_{R}&\kern -1mm\smash{\mathop{\hbox to
8mm{\rightarrowfill}} \limits_{\scriptstyle }}\kern
-1mm&\widetilde{A}_{\overline{R}}\cr}
$$
que l'on cherche \`{a} compl\'{e}ter en un diagramme
$$\diagram{
R&\twoheadrightarrow&\overline{R}\cr
\llap{$\scriptstyle $}\left\downarrow
\vbox to 4mm{}\right.\rlap{}&&\llap{}\left\downarrow
\vbox to 4mm{}\right.\rlap{$\scriptstyle $}\cr
A_{R}&\twoheadrightarrow&A_{\overline{R}}\cr
\llap{$\scriptstyle $}\left\downarrow
\vbox to 4mm{}\right.\rlap{}&&\llap{}\left\downarrow
\vbox to 4mm{}\right.\rlap{$\scriptstyle $}\cr
\widetilde{A}_{R}&\twoheadrightarrow &\widetilde{A}_{\overline{R}}\cr}
$$
o\`{u} $A_{R}$ est une d\'{e}formation plate de $A_{\overline{R}}$ sur
$R$. D'apr\`{e}s Illusie ([Il] Chapitre III, \S 2.3), ce probl\`{e}me
de rel\`{e}vement est contr\^{o}l\'{e} par un complexe $T$ qui
s'ins\`{e}re dans un triangle distingu\'{e}
$$
T\rightarrow \mathop{\rm RHom}\nolimits_{A_{\overline{R}}}
(L_{A_{\overline{R}}/\overline{R}},J\otimes_{\overline{R}}
A_{\overline{R}})\rightarrow \mathop{\rm RHom}
\nolimits_{\widetilde{A}_{\overline{R}}}(\widetilde{A}_{\overline{R}}
\otimes_{A_{\overline{R}}}L_{A_{\overline{R}}/\overline{R}},
J\otimes_{\overline{R}}\widetilde{A}_{\overline{R}})\rightarrow
$$
o\`{u} $L_{A_{\overline{R}}/\overline{R}}$ est le complexe cotangent
de la $\overline{R}$-alg\`{e}bre $A_{\overline{R}}$.  Plus
pr\'{e}cis\'{e}ment, il y a une obstruction \`{a} relever dans
$H^{2}(T)$ et, si cette obstruction est nulle, l'ensemble des classes
d'isomorphie de rel\`{e}vements est un torseur sous $H^{1}(T)$ et le
groupe des automorphismes d'un rel\`{e}vement donn\'{e} s'identifie
canoniquement \`{a} $H^{0}(T)$.  Comme la d\'{e}formation plate
$A_{\overline{R}}$ de $A$ sur $\overline{R}$ est n\'{e}cessairement de
la forme $A_{\overline{R}}=\overline{R}[[x,y]]/(f_{\overline{R}})$
pour $f_{\overline{R}}\in\overline{R}[[x,y]]$, le complexe cotangent
$$
L_{A_{\overline{R}}/\overline{R}}=[A_{\overline{R}}\rightarrow
A_{\overline{R}}{\rm d}x\oplus A_{\overline{R}}{\rm d}y]
$$
est concentr\'{e} en degr\'{e}s $[-1,0]$ et son unique
diff\'{e}rentielle non triviale envoie $1$ sur
$df_{\overline{R}}=(\partial_{x}f_{\overline{R}}){\rm d}x \oplus
(\partial_{y}f_{\overline{R}}){\rm d}y$.  Comme de plus la fl\`{e}che
$A_{\overline{R}}\rightarrow \widetilde{A}_{\overline{R}}$ est
n\'{e}cessairement injective \`{a} conoyau $\overline{R}$-plat
d'apr\`{e}s le lemme 3.4.2, le complexe
$$
T=[J\otimes_{\overline{R}}(\widetilde{A}_{\overline{R}}/A_{\overline{R}})
\partial_{x}\oplus J\otimes_{\overline{R}}
(\widetilde{A}_{\overline{R}}/A_{\overline{R}})\partial_{y}
\rightarrow J\otimes_{\overline{R}}(\widetilde{A}_{\overline{R}}/
A_{\overline{R}})]
$$
est concentr\'{e} en degr\'{e}s $[1,2]$ et son unique
diff\'{e}rentielle envoie $\widetilde{a}\partial_{x}\oplus
\widetilde{b}\partial_{y}$ sur $\widetilde{a}\partial_{x}
f_{\overline{R}} + \widetilde{b}\partial_{y}f_{\overline{R}}$.  On a
donc $H^{0}(T)=(0)$ et la suite exacte
$$
0\rightarrow H^{1}(T)\rightarrow J\otimes_{\overline{R}}
(\widetilde{A}_{\overline{R}}/A_{\overline{R}})
\partial_{x}\oplus J\otimes_{\overline{R}}
(\widetilde{A}_{\overline{R}}/A_{\overline{R}})\partial_{y}
\rightarrow J\otimes_{\overline{R}}(\widetilde{A}_{\overline{R}}/
A_{\overline{R}})\rightarrow H^{2}(T)\rightarrow 0.
$$
L'obstruction au rel\`{e}vement se calcule comme suit. On se donne des
rel\`{e}vements arbitraires $f_{R}\in R[[x,y]]$ de $f_{\overline{R}}$
et $R[[x,y]]\rightarrow\widetilde{A}_{R}$ de l'homomorphisme de
$\overline{R}$-alg\`{e}bres compos\'{e} $\overline{R}[[x,y]]
\twoheadrightarrow A_{\overline{R}}\hookrightarrow
\widetilde{A}_{\overline{R}}$. Alors l'image de $f_{R}$ par
l'homomorphisme de $R$-alg\`{e}bres $R[[x,y]]\rightarrow
\widetilde{A}_{R}$ est dans $J\widetilde{A}_{R}\subset
\widetilde{A}_{R}$. L'obstruction cherch\'{e}e est l'image de
l'\'{e}l\'{e}ment ainsi construit de $J\widetilde{A}_{R}$ par
l'\'{e}pimorphisme compos\'{e}
$$
J\widetilde{A}_{R}\cong J\otimes_{\overline{R}}
\widetilde{A}_{\overline{R}}\twoheadrightarrow
J\otimes_{\overline{R}}(\widetilde{A}_{\overline{R}}/A_{\overline{R}})
\twoheadrightarrow H^{2}(T).
$$
Bien entendu, cette obstruction ne d\'{e}pend pas des choix faits et
elle est nulle d'apr\`{e}s le corollaire 3.4.5.

En particulier, l'espace tangent relatif au morphisme de foncteurs
$\mathop{\rm Def}\nolimits_{\widetilde{A}}^{{\rm top}}\rightarrow
\mathop{\rm Def}\nolimits_{A\hookrightarrow\widetilde{A}}^{{\rm top}}$
est le noyau de l'application $k$-lin\'{e}aire
$$
(\widetilde{A}/A)\partial_{x}\oplus
(\widetilde{A}/A)\partial_{y}\rightarrow
(\widetilde{A}/A),~\widetilde{a}\partial_{x}\oplus
\widetilde{b}\partial_{y}\mapsto \widetilde{a}\partial_{x}f+
\widetilde{b}\partial_{y}f,
$$
application qui est identiquement nulle d'apr\`{e}s le lemme 3.4.1.
Cet espace tangent est donc \'{e}gal au $k$-espace vectoriel
$(\widetilde{A}/A)\partial_{x}\oplus (\widetilde{A}/A)\partial_{y}$ de
dimension $2\delta (A)$.

Comme $\widetilde{A}=\prod_{i\in I}k[[t_{i}]]$, le foncteur
$\mathop{\rm Def}\nolimits_{\widetilde{A}}^{{\rm top}}$ est trivial.
Plus pr\'{e}cis\'{e}ment, toute d\'{e}formation plate de
$\widetilde{A}$ sur $R\in\mathop{\rm ob}\mathop{\rm Art}\nolimits_{k}$
est isomorphe \`{a} $\prod_{i\in I}R[[t_{i}]]$ et, en termes de
rel\`{e}vements, pour tout id\'{e}al $J$ de carr\'{e} nul dans $R$ et
pour $\overline{R}=R/J$, il n'y a pas d'obstruction \`{a} relever
$\widetilde{A}_{\overline{R}} =\prod_{i\in I}\overline{R}[[t_{i}]]$
\`{a} $R$, il n'y a qu'une seule classe d'isomorphie de tels
rel\`{e}vements, \`{a} savoir celle de $\widetilde{A}_{R}=\prod_{i\in
I}R[[t_{i}]]$, et le groupe des automorphismes d'un rel\`{e}vement
arbitraire dans cette classe est le groupe
$$
\mathop{\rm Hom}\nolimits_{\widetilde{A}_{\overline{R}}}
(\Omega_{\widetilde{A}_{\overline{R}}/\overline{R}}^{1},
J\otimes_{\overline{R}}\widetilde{A}_{\overline{R}})=\prod_{i\in
I}(J\otimes_{\overline{R}}\overline{R}[[t_{i}]])
\partial_{t_{i}}.
$$

On a donc montr\'{e}:

\thm PROPOSITION 3.4.7
\enonce
Soient $R\in\mathop{\rm ob}\mathop{\rm Art}\nolimits_{k}$, $J\subset
R$ un id\'{e}al de carr\'{e} nul et $\overline{R}=R/J$.

\decale{\rm (i)} Tout objet $A_{R}\hookrightarrow\widetilde{A}_{R}$ de
$\mathop{\rm Def}\nolimits_{A\hookrightarrow\widetilde{A}}^{{\rm
top}}(R)$ est de la forme
$$
R[[x,y]]/(f_{R})\hookrightarrow \prod_{i\in I}R[[t_{i}]],~x\mapsto
(x_{R,i}(t_{i}))_{i\in I},~y\mapsto
(y_{R,i}(t_{i}))_{i\in I},
$$
pour des s\'{e}ries $f_{R}\in R[[x,y]]$ et $x_{R,i}(t_{i}),
y_{R,i}(t_{i})\in R[[t_{i}]]$, qui rel\`{e}vent les s\'{e}ries
$f$ et $x_{i}(t),y_{i}(t)$, et qui v\'{e}rifient bien s\^{u}r
$$
f_{R}(x_{R,i}(t_{i}),y_{R,i}(t_{i}))\equiv 0,~\forall i\in I.
$$

\decale{\rm (ii)} Soit $(A_{\overline{R}}\hookrightarrow
\widetilde{A}_{\overline{R}})\in \mathop{\rm Def}
\nolimits_{A\hookrightarrow \widetilde{A}}^{{\rm top}}(\overline{R})$
isomorphe \`{a}
$$
\overline{R}[[x,y]]/(f_{\overline{R}})\hookrightarrow\prod_{i\in
I}\overline{R}[[t_{i}]],~x\mapsto (x_{\overline{R},i}(t_{i}))_{i\in
I},~y\mapsto (y_{\overline{R},i}(t_{i}))_{i\in I}.
$$
Il n'y a pas d'obstruction \`{a} relever cet objet \`{a} $R$,
l'ensemble des classes d'isomorphie des rel\`{e}vements est le conoyau
de la fl\`{e}che
$$
\bigoplus_{i\in I}(J\otimes_{\overline{R}}\overline{R}[[t_{i}]])
\partial_{t_{i}}\rightarrow \mathop{\rm Ker}(J \otimes_{\overline{R}}
(\widetilde{A}_{\overline{R}}/A_{\overline{R}}) \partial_{x}\oplus
J\otimes_{\overline{R}} (\widetilde{A}_{\overline{R}}/
A_{\overline{R}})\partial_{y}\rightarrow
J\otimes_{\overline{R}}(\widetilde{A}_{\overline{R}}/
A_{\overline{R}}))
$$
qui envoie $\oplus_{i\in I}a_{i}(t_{i})\partial_{t_{i}}$ sur
$$
\bigl((a_{i}(t_{i})(\partial_{t_{i}}x_{\overline{R},i})(t_{i}))_{i\in
I}+ J\otimes_{\overline{R}}A_{\overline{R}}\bigr) \partial_{x}\oplus
\bigl((a_{i}(t_{i})(\partial_{t_{i}}y_{\overline{R},i})(t_{i}))_{i\in
I}+ J\otimes_{\overline{R}} A_{\overline{R}}\bigr)\partial_{y}
$$
et que le groupe des automorphismes d'un rel\`{e}vement arbitraire
dans cette classe est son noyau.
\hfill\hfill$\square$
\endthm

\thm COROLLAIRE 3.4.8
\enonce
L'espace tangent au foncteur formellement lisse $\mathop{\rm
Def}\nolimits_{A\hookrightarrow \widetilde{A}}^{{\rm top}}$ est le
conoyau de l'application $k$-lin\'{e}aire
$$
\bigoplus_{i\in I}k[[t_{i}]]\partial_{t_{i}}\rightarrow
(\widetilde{A}/A)\partial_{x}\oplus (\widetilde{A}/A)
\partial_{y}
$$
qui envoie $\oplus_{i\in I}a_{i}(t_{i})\partial_{t_{i}}$ sur
$$
\bigl((a_{i}(t_{i})(\partial_{t_{i}}x_{i})(t_{i}))_{i\in
I}+A\bigr)\partial_{x}\oplus
\bigl((a_{i}(t_{i})(\partial_{t_{i}}y_{i})(t_{i}))_{i\in I}
+A\bigr)\partial_{y}.
$$

De plus, l'application tangente au morphisme de foncteur $\mathop{\rm
Def}\nolimits_{A\hookrightarrow \widetilde{A}}^{{\rm top}}\rightarrow
\mathop{\rm Def}\nolimits_{A}^{{\rm top}}$ est induite par
l'application $k$-lin\'{e}aire
$$
(\widetilde{A}/A)\partial_{x}\oplus
(\widetilde{A}/A)\partial_{y}\rightarrow A/(\partial_{x}f,\partial_{y}f)
$$
qui envoie la classe d'un \'{e}l\'{e}ment $\widetilde{a}\partial_{x}\oplus
\widetilde{b}\partial_{y}\in
\widetilde{A}\partial_{x}\oplus
\widetilde{A}\partial_{y}$ sur la classe de
l'\'{e}l\'{e}ment $-(\widetilde{a}\partial_{x}f+
\widetilde{b}\partial_{y}f)\in A\subset \widetilde{A}$.

\hfill\hfill$\square$
\endthm

L'\'{e}nonc\'{e} suivant est d\^{u} \`{a} Diaz et Harris Il devrait
\^{e}tre v\'{e}rifi\'{e} pour $k$ de caract\'{e}ristique arbitraire.
Faute de r\'{e}f\'{e}rence nous nous limiterons \`{a} la
caract\'{e}ristique nulle.

\thm TH\'{E}OR\`{E}ME 3.4.9 (Diaz et Harris, [D-H] Proposition (4.17)
et Theorem (4.15))
\enonce
Supposons de plus que $k$ est de caract\'{e}ristique nulle.  Alors, le
morphisme de $k$-sch\'{e}mas formels $\mathop{\rm Def}
\nolimits_{A\hookrightarrow \widetilde{A}}^{{\rm top}}
\rightarrow\mathop{\rm Def}\nolimits_{A}^{{\rm top}}$ est fini, son
image sch\'{e}matique
$$
\mathop{\rm Def}\nolimits_{A}^{{\rm top},\delta}\subset\mathop{\rm
Def}\nolimits_{A}^{{\rm top}}
$$
est une ferm\'{e} int\`{e}gre de codimension $\delta (A)$ et le
morphisme canonique $\mathop{\rm Def}\nolimits_{A\hookrightarrow
\widetilde{A}}^{{\rm top}}\rightarrow \mathop{\rm
Def}\nolimits_{A}^{{\rm top},\delta}$ est le morphisme de
normalisation.

De plus, le c\^{o}ne tangent de $\mathop{\rm Def}
\nolimits_{A}^{{\rm top},\delta}$ est le sous-$k$-espace vectoriel
$$
V(A)={\frak a}/(\partial_{x}f,\partial_{y}f)\subset A/(\partial_{x}f,
\partial_{y}f)=\mathop{\rm Def}\nolimits_{A}^{{\rm top}}
(k[\varepsilon])
$$
de l'espace tangent de $\mathop{\rm Def}\nolimits_{A}^{{\rm top}}$,
o\`{u} ${\frak a}\subset A$ est le conducteur du normalis\'{e}
$\widetilde{A}$ de $A$ dans $A$ {\rm (}cf.  Lemme $3.4.1${\rm )}.
\hfill\hfill$\square$
\endthm

Soient maintenant $C_{k}$ une courbe projective, int\`{e}gre et \`{a}
singularit\'{e}s planes isol\'{e}es sur $k$, et $\varphi :C\rightarrow
S$ une alg\'{e}brisation d'une d\'{e}formation formelle miniverselle
de $C_{k}$ comme dans la section $3.2$.  Toutes les fibres de
$\varphi$ sont g\'{e}om\'{e}triquement int\`{e}gre et le lieu
singulier $\bigcup_{t\in S}C_{t}^{{\rm sing}}$ est fini sur $S$.  On
peut donc appliquer le corollaire $3.3.4$ \`{a} cette courbe et on
obtient la strate \`{a} $\delta$ contant
$$
S^{\delta}\subset S.
$$

Pour chaque point singulier $c$ de $C_{k}$, on peut aussi
consid\'{e}rer une d\'{e}formation miniverselle ${\cal
C}_{c}=\mathop{\rm Spf}({\cal A}_{c})\rightarrow {\cal
S}_{c}=\mathop{\rm Spf}({\cal R}_{c})$ de ${\cal C}_{c,s}=\mathop{\rm
Spf}(\widehat{{\cal O}}_{C_{k},c})$.  Il r\'{e}sulte du
th\'{e}or\`{e}me 3.1.3 que l'on a un morphisme canonique formellement
lisse
$$
{\cal S}\rightarrow\prod_{c\in C_{k}^{{\rm sing}}}{\cal S}_{c}
$$
o\`{u} ${\cal S}=\mathop{\rm Spf}({\cal O}_{S,s})$ est le
compl\'{e}t\'{e} de $S$ en son point ferm\'{e} $s$.

Pour chaque $c\in C_{k}^{{\rm
sing}}$, soit
$$
{\cal S}_{c}^{\delta}\subset {\cal S}_{c}
$$
la strate \`{a} $\delta$ constant.

\thm LEMME 3.4.10
\enonce
Le compl\'{e}t\'{e} de $S^{\delta}$ au point ferm\'{e} $s$ est le
ferm\'{e} de ${\cal S}$ image r\'{e}ciproque du ferm\'{e} $\prod_{c\in
C_{k}^{{\rm sing}}}{\cal S}_{c}^{\delta}$ de $\prod_{c\in C_{k}^{{\rm
sing}}}{\cal S}_{c}$ par le morphisme canonique ci-dessus.
\hfill\hfill$\square$
\endthm

L'\'{e}nonc\'{e} suivant est une variante globale du th\'{e}or\`{e}me
3.4.9.  Tout comme ce dernier th\'{e}or\`{e}me, il devrait \^{e}tre
v\'{e}rifi\'{e} pour $k$ de caract\'{e}ristique arbitraire.  Faute de
r\'{e}f\'{e}rence nous nous limiterons de nouveau \`{a} la
caract\'{e}ristique nulle.  La derni\`{e}re assertion est bien connue
mais ne semble \^{e}tre d\'{e}montr\'{e}e nulle part; elle r\'{e}sulte
du th\'{e}or\`{e}me (1.3) de [D-H] dans le cas o\`{u} $C_{k}$ est une
courbe plane.

\thm TH\'{E}OR\`{E}ME 3.4.11 (Diaz et Harris)
\enonce
Supposons de plus que $k$ est de caract\'{e}ristique nulle.  Alors, la
strate $S^{\delta}$ est irr\'{e}ductible de codimension dans $S$
\'{e}gale \`{a} $\delta (C_{k})$.  Le sch\'{e}ma normalis\'{e} de
$S^{\delta}$ est formellement lisse sur $k$.

De plus, le c\^{o}ne tangent \`{a} l'origine de  $S^{\delta}$ est le
sous-$k$-espace vectoriel
$$
V(C_{k})\subset T_{s}S
$$
de l'espace tangent \`{a} l'origine $s$ de $S$ obtenu par image
inverse du sous-$k$-espace vectoriel
$$
\bigoplus_{c\in C_{k}^{{\rm sing}}}V(\widehat{{\cal O}}_{C_{s},c})
\subset\bigoplus_{c\in C_{k}^{{\rm sing}}} \mathop{\rm Def}
\nolimits_{\widehat{{\cal O}}_{C_{s},c}}^{{\rm top}}(k[\varepsilon])
$$
par l'\'{e}pimorphisme naturel
$T_{s}S\twoheadrightarrow\bigoplus_{c\in C_{k}^{{\rm sing}}}
\mathop{\rm Def}\nolimits_{\widehat{{\cal O}}_{C_{s},c}}^{{\rm
top}}(k[\varepsilon])$.

En outre, la fibre de $C\rightarrow S$ en tout point
g\'{e}om\'{e}trique g\'{e}n\'{e}rique de $S^{\delta}\subset S$ est une
courbe n'ayant comme seules singularit\'{e}s que $\delta (C_{k})$
points doubles ordinaires.
\hfill\hfill$\square$
\endthm
\vskip 5mm

\centerline{4. {\og}{D\'{E}FORMATIONS}{\fg} DES FIBRES DE SPRINGER}
\vskip 2mm

\section{4.1}{Retour aux fibres de Springer}

Revenons \`{a} la situation de la premi\`{e}re partie.  On rappelle
(cf.  \S 2.1) qu'on a introduit une courbe int\`{e}gre, projective
$C=C_{I}$ sur $k$ qui n'a qu'un seul point singulier $c=c_{I}$ pour
lequel le compl\'{e}t\'{e} $\widehat{{\cal O}}_{C,c}$ de l'anneau
local de $C$ en $c$ est isomorphe \`{a}
$$
A_{I}\cong {\cal O}_{F}[\gamma_{I}]\cong k[[\varpi_{F},T]]/(P_{I}(T)).
$$
Par hypoth\`{e}se de s\'{e}parabilit\'{e} de $P_{I}(T)$, on sait que
l'id\'{e}al
$$
(\partial_{T}P_{I}(T),P_{I}(T))\subset k[[\varpi_{F}]][T]
$$
est de codimension finie. Il en est donc de m\^{e}me des id\'{e}aux
$$
(\partial_{T}P_{I}(T),P_{I}(T))\hbox{ et }
(\partial_{\varpi_{F}}P_{I}(T),\partial_{T}P_{I}(T),P_{I}(T))\subset
k[[\varpi_{F},T]].
$$
On supposera dans la suite que l'id\'{e}al
$$
(\partial_{\varpi_{F}}P_{I}(T),\partial_{T}P_{I}(T))\subset
k[[\varpi_{F},T]]
$$
est lui aussi de codimension finie (comme on l'a d\'{e}j\`{a} dit,
cette hypoth\`{e}se est automatiquement v\'{e}rifi\'{e}e si $k$ est de
caract\'{e}ristique nulle).
\vskip 2mm

Par construction la courbe $C$ est muni d'un point $\infty$
dans son lieu de lissit\'{e} et la normalis\'{e}e $\widetilde{C}$
de $C$ est identifi\'{e}e \`{a} la droite projective ${\Bbb
P}_{k}^{1}$ de telle sorte que $\infty$ soit le point \`{a}
l'infini (cf. la preuve de la proposition $2.1.1$).

\thm LEMME 4.1.1
\enonce
Supposons de plus que $C$ n'est pas lisse en $c$ et que la
caract\'{e}ristique de $k$ est $>|I|$.  Alors, le $k$-sch\'{e}ma en
groupes des automorphismes $(C,\infty)$ sur $k$ est soit fini, soit
isomorphe \`{a} ${\Bbb G}_{{\rm m},k}$.
\endthm

\rem Preuve
\endrem
Tout automorphisme $g$ de $(C,\infty )$ induit un automorphisme de
$\widetilde{C}={\Bbb P}_{k}^{1}$ qui fixe le point \`{a} l'infini et
$\{\widetilde{c}_{i}\mid i\in I\}\subset \widetilde{C}$ dans son
ensemble.  Par suite, $\mathop{\rm Aut}\nolimits_{k}(C,\infty)$ est un
sous-$k$-sch\'{e}ma en groupes ferm\'{e} du sous-groupe de Borel des
matrices triangulaires sup\'{e}rieures dans $\mathop{\rm PGL}(2)$.  De
plus, la puissance $g^{|I|!}$ de $g$ fixe chacun des
$\widetilde{c}_{i}$.  Par suite, si $|I|\geq 2$, $g^{|I|!}$ est
n\'{e}cessairement l'identit\'{e} et $\mathop{\rm Aut}
\nolimits_{k}(C,\infty)$ est un $k$-sch\'{e}ma en groupes fini d'ordre
premier \`{a} la caract\'{e}ristique de $k$, et si $|I|=1$,
$\mathop{\rm Aut}\nolimits_{k}(C,\infty)$ est un sous-$k$-sch\'{e}ma
en groupes du tore maximal diagonal ${\Bbb G}_{{\rm m}}\subset
\mathop{\rm PGL}(2)$
\hfill\hfill$\square$
\vskip 3mm

On s'int\'{e}resse de nouveau \`{a} la fibre de Springer $X=X_{I}$ et
\`{a} son quotient $Z=Z_{I}$ par le r\'{e}seau $\Lambda^{0}=
\Lambda_{I}^{0}$.

On a vu d'une part que $Z$ est naturellement hom\'{e}omorphe au
$k$-sch\'{e}ma de Picard compactifi\'{e} $\overline{P}=
\overline{P}_{I}$ et que le rev\^{e}tement $X\rightarrow Z$ provient
d'un rev\^{e}tement $\overline{P}^{\,\natural}=
\overline{P}_{I}^{\,\natural}\rightarrow\overline{P}$.  On a vu
d'autre part comment d\'{e}former la courbe $C$.  On se propose
maintenant d'utiliser les d\'{e}formations de $C$ pour d\'{e}former
$\overline{P}$, et aussi d'une certaine mani\`{e}re
$\overline{P}^{\,\natural}$.
\vskip 2mm

{\it On note dor\'{e}navant les objets ci-dessus par
$X_{s}=X_{I,s_{I}}$, $Z_{s}=Z_{I,s_{I}}$, $\Lambda_{s}^{0}=
\Lambda_{I,s_{I}}^{0}$, $C_{s}=C_{I,s_{I}}$, $\infty_{s}=\infty$,
$P_{s}=P_{I,s_{I}}$, $\overline{P}_{s}= \overline{P}_{I,s_{I}}$, ...,
ce qui lib\`{e}re les notation $X$, $Z$, $\Lambda$, $C$, $\infty$,
$P$, $\overline{P}$, ...  que nous allons pouvoir utiliser pour les
d\'{e}formations de ces m\^{e}mes objets.}
\vskip 2mm

Ainsi, on note simplement $(C=C_{I},\infty =\infty_{I})\rightarrow
S_{I}=S$ une alg\'{e}brisation de la d\'{e}formation miniverselle
(presque universelle compte tenu du lemme 4.1.1) de
$(C_{s},\infty_{s})$ et $s=s_{I}$ est l'unique point ferm\'{e} de $S$.
Le sch\'{e}ma $S$ est le spectre d'une $k$-alg\`{e}bre de s\'{e}ries
formelles en $\tau (C_{s},\infty_{s})$ variables, o\`{u} $\tau
(C_{s},\infty_{s})$ est la dimension sur $k$ de $\mathop{\rm Ext}
\nolimits_{{\cal C}_{s}}^{1}(\Omega_{C_{s}/k}^{1} (\infty_{s}),{\cal
O}_{C_{s}})$.  La courbe relative $C$ est projective et plate,
$\infty$ est une section dans $C/S$ \`{a} image contenue dans le lieu
de lissit\'{e} et toutes les fibres g\'{e}om\'{e}triques $C_{t}$ ($t$
un point g\'{e}om\'{e}trique de $S$) sont int\`{e}gres et \`{a}
singularit\'{e}s planes.  La fibre g\'{e}n\'{e}rique
g\'{e}om\'{e}trique de $C\rightarrow S$ est lisse (cf.  Proposition
3.2.3).
\vskip 2mm

Soient $P=P_{I}\subset\overline{P}{}_{I}=\overline{P}$ les
sch\'{e}mas de Picard et de Picard compactifi\'{e} de $C/S$.  On
rappelle que $\overline{P}$ param\`{e}tre les classes de ${\cal
O}_{C}$-Modules plats sur $S$ qui sont fibre \`{a} fibre sans torsion
et de rang g\'{e}n\'{e}rique $1$, et que $P$ est l'ouvert de
$\overline{P}$ form\'{e} des classes de Modules inversibles.

Les sch\'{e}mas $P$ et $\overline{P}$ sont purement de dimension
relative $\delta_{I}$ sur $S$, r\'{e}unions disjointes de composantes
connexes $P^{d}= P_{I}^{d}\subset\overline{P}{}_{I}^{d}=
\overline{P}{}^{d}$, respectivement quasi-projectives et projectives
sur $S$.  Les fibres $P_{t}^{d}$ et $\overline{P}{}_{t}^{d}$ de
$P^{d}\rightarrow S$ et $\overline{P}{}^{d}\rightarrow S$ en tout
point g\'{e}om\'{e}trique $t$ de $S$ sont les composantes de degr\'{e}
$d$ des sch\'{e}mas de Picard et de Picard compactifi\'{e} de la
fibre $C_{t}$ de $C\rightarrow S$ en $t$; chaque
$\overline{P}{}_{t}^{d}$ est int\`{e}gre et localement d'intersection
compl\`{e}te, d'apr\`{e}s le th\'{e}or\`{e}me de Rego et Altman,
Iarrobino et Kleiman (cf.  Th\'{e}or\`{e}me 2.1.2).

\thm TH\'{E}OR\`{E}ME 4.1.2 (Fantechi, G\"{o}ttsche et van Straten;
[F-G-S] Corollary B.2)
\enonce
Le sch\'{e}ma $\overline{P}$ est r\'{e}gulier {\rm (}formellement
lisse sur $k${\rm )} et la fibre sp\'{e}ciale $\overline{P}{}_{s}$ de
$\overline{P}\rightarrow S$ est localement d'intersection compl\`{e}te
dans $\overline{P}$.
\endthm

Pour prouver ce th\'{e}or\`{e}me, Fantechi, G\"{o}ttsche et van Straten
consid\`{e}rent le foncteur des d\'{e}formations
$$
\mathop{\rm Def}\nolimits_{C_{s},{\cal M}_{s}}:\mathop{\rm Art}
\nolimits_{k}\rightarrow \mathop{\rm Ens}
$$
du couple $(C_{s},{\cal M}_{s})$ o\`{u} ${\cal M}_{s}$ est un ${\cal
O}_{C_{s}}$-Module coh\'{e}rent sans torsion de rang g\'{e}n\'{e}rique
$1$, et sa variante locale
$$
\mathop{\rm Def}\nolimits_{\widehat{{\cal O}}_{C_{s},c},
\widehat{{\cal M}}_{s,c}}^{{\rm top}}:\mathop{\rm Art}\nolimits_{k}
\rightarrow\mathop{\rm Ens}
$$
en l'unique point singulier $c$ de $C_{s}$.  On a le carr\'{e}
commutatif de morphismes naturels de foncteurs
$$\diagram{
\mathop{\rm Def}\nolimits_{C_{s},{\cal M}_{s}} &\kern -1mm
\smash{\mathop{\hbox to 8mm{\rightarrowfill}} \limits^{\scriptstyle }}
\kern -1mm&\mathop{\rm Def}
\nolimits_{\widehat{{\cal O}}_{C_{s},c},\widehat{{\cal 
M}}_{s,c}}^{{\rm top}}\cr
\llap{$\scriptstyle $}\left\downarrow \vbox to 4mm{}\right.\rlap{}&+&
\llap{}\left\downarrow \vbox to 4mm{}\right.\rlap{$\scriptstyle $}\cr
\mathop{\rm Def}\nolimits_{C_{s}}&\kern -2mm\smash{\mathop{\hbox to
10mm{\rightarrowfill}}\limits_{\scriptstyle }}\kern -2mm&
\mathop{\rm Def}\nolimits_{\widehat{{\cal
O}}_{C_{s},c}}^{{\rm top}}\cr}
$$
et il est facile de v\'{e}rifier :

\thm LEMME 4.1.3
\enonce
Le morphisme de foncteurs
$$
\mathop{\rm Def}\nolimits_{C_{s},{\cal M}_{s}} \rightarrow\mathop{\rm
Def}\nolimits_{C_{s}}\times_{\mathop{\rm Def}
\nolimits_{\widehat{{\cal O}}_{C_{s},c}}^{{\rm top}}}
\mathop{\rm Def} \nolimits_{\widehat{{\cal O}}_{C_{s},c},
\widehat{{\cal M}}_{s,c}}^{{\rm top}}
$$
induit par le diagramme ci-dessus est formellement lisse.
\hfill\hfill$\square$
\endthm

Le th\'{e}or\`{e}me r\'{e}sulte donc de la proposition suivante:

\thm PROPOSITION 4.1.4
\enonce
Soient $f\in (x,y)\subset k[[x,y]]$ tel que l'id\'{e}al
$(\partial_{x}f,\partial_{y}f)$ soit de codimension finie dans
$k[[x,y]]$, $A= k[[x,y]]/(f)$ et $M$ un $A$-module de type fini, sans
torsion et de rang g\'{e}n\'{e}rique $1$.  Alors, le foncteur des
d\'{e}formations $\mathop{\rm Def}\nolimits_{A,M}^{{\rm top}}$ est
formellement lisse sur $k$
\endthm

\rem Preuve
\endrem
D'apr\`{e}s le lemme 3.4.3, $M$ admet, en tant que $k[[x,y]]$-module,
une r\'{e}solution
$$
0\rightarrow k[[x,y]]^{n}\,\smash{\mathop{\hbox to
8mm{\rightarrowfill}} \limits^{\scriptstyle F}}\,
k[[x,y]]^{n}\rightarrow M\rightarrow 0
$$
o\`{u} $n\in \{1,\ldots ,\delta (A)+1\}$ et $F$ est une
matrice carr\'{e} de taille $n\times n$ \`{a} coefficients dans
$k[[x,y]]$.  En raisonnant comme dans la preuve du th\'{e}or\`{e}me
3.4.4 avec la matrice des co-facteurs de $F$, on peut exiger de plus
que le d\'{e}terminant $\det F$ est \'{e}gal \`{a} $f$.

Toute d\'{e}formation topologique ($R$-plate) $M_{R}$ du
$k[[x,y]]$-module $M$ sur $R\in \mathop{\rm ob}\mathop{\rm
Art}\nolimits$ peut-\^{e}tre obtenue en d\'{e}formant la matrice $F$
en une matrice carr\'{e}e $F_{R}$ de taille $n\times n$ \`{a}
coefficients dans $R[[x,y]]= R\otimes_{k}k[[x,y]]$, de telle sorte que
$M_{R}$ admette la pr\'{e}sentation
$$
0\rightarrow R[[x,y]]^{n}\,\smash{\mathop{\hbox to 8mm
{\rightarrowfill}}\limits^{\scriptstyle F_{R}}}\,R[[x,y]]^{n}
\rightarrow M_{R}\rightarrow 0.
$$

Pour une telle d\'{e}formation topologique $M_{R}$ de $M$ en tant que
$k[[x,y]]$-module, la $R$-alg\`{e}bre quotient $A_{R}=R[[x,y]]/(\det
F_{R})$ est une d\'{e}formation topologique de $A$ sur $R$ et $M_{R}$
est de mani\`{e}re \'{e}vidente un $A_{R}$-module sans torsion.  Si
maintenant $B_{R}=R[[x,y]]/(g_{R})$ est une autre d\'{e}formation
plate de $A$ sur $R$ telle que $g_{R}M_{R}=(0)$, alors on a
$g_{R}=h_{R}\det F_{R}$ o\`{u} $h_{R}\in R[[x,y]]$ est congru modulo
l'id\'{e}al maximal de $R$ \`{a} un \'{e}l\'{e}ment inversible de
$k[[x,y]]$ et est donc inversible dans $R[[x,y]]$, de sorte que
$B_{R}=A_{R}$.

Si $\mathop{\rm Def}\nolimits_{F}^{{\rm top}}:\mathop{\rm
Art}\nolimits_{k}\nolimits \rightarrow \mathop{\rm Ens}$ est le
foncteur des d\'{e}formations de la matrice $F$ ci-dessus, on a
donc construit un morphisme formellement lisse de foncteurs
$$
\mathop{\rm Def}\nolimits_{F}^{{\rm top}}\rightarrow \mathop{\rm
Def}\nolimits_{A,M}^{{\rm top}},~F_{R}\mapsto (R[[x,y]]/(\det
F_{R}),\mathop{\rm Coker}(F_{R})).
$$
Comme le foncteur $\mathop{\rm Def}\nolimits_{F}^{{\rm top}}$ est
trivialement formellement lisse sur $k$, la proposition s'en suit.
\hfill\hfill$\square$
\vskip 3mm

Le r\'{e}sultat suivant g\'{e}n\'{e}ralise le th\'{e}or\`{e}me
4.1.2.

\thm TH\'{E}OR\`{E}ME 4.1.5 (Fantechi, G\"{o}ttsche, van Straten;
[F-G-S] Corollary B.3)
\enonce
Soit $C_{T}\rightarrow T$ une courbe relative de base $T\cong
\mathop{\rm Spec}(k[[t_{1},\ldots ,t_{m}]])$ qui provient de la courbe
universelle $C\rightarrow S$ par un changement de base local
$T\rightarrow S$.  Alors, le $T$-sch\'{e}ma
$$
\overline{P}_{T}=T\times_{S}\overline{P}
$$
de Picard compactifi\'{e} de $C_{T}/T$ est r\'{e}gulier {\rm
(}formellement lisse sur $k${\rm )} si et seulement si l'espace
tangent \`{a} l'origine de $T$ est transverse au sous-espace
$V(C_{k})\subset T_{s}S$ de l'espace tangent \`{a} l'origine de $S$
introduit dans le th\'{e}or\`{e}me $3.4.11$.
\endthm

\rem Preuve
\endrem
Reprenons les notations de la proposition 4.1.4 et de sa preuve.  On a
un morphisme d'oubli
$$
\mathop{\rm Def}\nolimits_{A,M}^{{\rm top}}\rightarrow
\mathop{\rm Def}\nolimits_{A}^{{\rm top}}
$$
et il suffit de d\'{e}montrer que l'image de l'application tangente
\`{a} ce morphisme contient le sous-espace $V(A)$ de l'espace tangent
\`{a} $\mathop{\rm Def}\nolimits_{A}^{{\rm top}}$ (cf. 3.4.9).

L'image de cette application tangente est le quotient
$I/(f,\partial_{x}f, \partial_{y}f)$ o\`{u} $I\subset k[[x,y]]$ est
l'id\'{e}al engendr\'{e}e par les entr\'{e}es de la matrice $F^{\ast}$
des co-facteurs de $F$ puisque l'on a
$$
\det (F+\varepsilon E^{ij})=\det F + \varepsilon F_{ji}^{\ast}
$$
quels que soient $0\leq i,j\leq n$, o\`{u} $E^{ij}$ est la matrice
\'{e}l\'{e}mentaire dont l'entr\'{e}e $(i,j)$ est \'{e}gale \`{a} $1$
et dont toutes les autres entr\'{e}es sont nulles.

Or, le $A$-module $M=\mathop{\rm Coker}(F)$ admet la r\'{e}solution
p\'{e}riodique, de p\'{e}riode $2$,
$$
\cdots\,\smash{\mathop{\hbox to 8mm{\rightarrowfill}}
\limits^{\scriptstyle F}}\,A^{n}\,\smash{\mathop{\hbox to
8mm{\rightarrowfill}} \limits^{\scriptstyle
F^{\ast}}}\,A^{n}\,\smash{\mathop{\hbox to 8mm{\rightarrowfill}}
\limits^{\scriptstyle F}}\,A^{n}\rightarrow M\rightarrow 0,
$$
et $M$ est encore \'{e}gal \`{a}
$$
M=\mathop{\rm Ker}(F)=\mathop{\rm Im}(F^{\ast}).
$$
Comme on a $\mathop{\rm Ext}\nolimits_{A}^{1}(N,A)=(0)$ pour tout
$A$-module sans torsion $N$, le morphisme $\mathop{\rm
Hom}\nolimits_{A}(F^{\ast},A)$ admet la factorisation canonique
$$
\mathop{\rm Hom}\nolimits_{A}(A^{n},A)\twoheadrightarrow \mathop{\rm
Hom}\nolimits_{A}(M,A)\hookrightarrow \mathop{\rm Hom}
\nolimits_{A}(A^{n},A)
$$
et l'image $I/(f)$ de $I$ dans $A$ est \'{e}gale \`{a} l'id\'{e}al
engendr\'{e} par les $\varphi (m)$ pour $\varphi$ parcourant
$\mathop{\rm Hom}\nolimits_{A}(M,A)$ et $m$ parcourant $M$.

Mais, \`{a} isomorphisme pr\`{e}s, on peut supposer que $A\subset
M\subset \widetilde{A}$ (cf. la preuve du lemme 3.4.3), et alors il
est clair que
$$
\{\varphi (m)\mid\varphi\in \mathop{\rm Hom}\nolimits_{A}(M,A),~
m\in M\}\supset {\frak a},
$$
ce qui termine la preuve du th\'{e}or\`{e}me.
\hfill\hfill$\square$
\vskip 3mm

\section{4.2}{Application \`{a} la puret\'{e} dans le cas homog\`{e}ne}

Soient $m>n\geq 1$ premiers entre eux.  On suppose que $m!$ est
inversible dans $k$.  On consid\`{e}re la fibre de Springer $X$
associ\'{e}e \`{a} l'extension totalement ramifi\'{e}e $F\subset E$ de
degr\'{e} $n$ d\'{e}finie par $\varpi_{F}=\varpi_{E}^{n}$ et \`{a}
l'\'{e}l\'{e}ment $\gamma =\varpi_{E}^{m}$, ou encore celle
associ\'{e}e \`{a} l'extension totalement ramifi\'{e}e $F\subset E$ de
degr\'{e} $m$ d\'{e}finie par $\varpi_{F}=\varpi_{E}^{m}$ et \`{a}
l'\'{e}l\'{e}ment $\gamma =\varpi_{E}^{n}$.  Dans les deux cas, le
germe formel de courbe plane correspondant est le m\^{e}me, \`{a}
savoir $\mathop{\rm Spf}(A)$ o\`{u} $A=k[[x,y]]/(x^{m}-y^{n})$.  Comme
ce germe n'a qu'une seule branche, on a $\Lambda ={\Bbb Z}$ et
$X=Z=Z^{0}\times {\Bbb Z}$.

Rappelons que, si $C_{k}$ est n'importe quelle courbe int\`{e}gre,
projective sur $k$, de normalis\'{e}e isomorphe \`{a} la droite
projective ${\Bbb P}_{k}^{1}$, et qui est munie d'un $k$-point $c$ en
dehors duquel elle est lisse sur $k$ et en lequel le germe formel de
$C_{k}$ est isomorphe \`{a} $\mathop{\rm Spf}(A)$,
alors le $k$-sch\'{e}ma $Z^{0}$ est hom\'{e}omorphe \`{a} la
composante $\mathop{\overline{\rm Pic}}\nolimits_{C_{k}/k}^{\,0}$ de
degr\'{e} $0$ du $k$-sch\'{e}ma de Picard compactifi\'{e} de $C_{k}$.

On peut construire une telle courbe $C_{k}$ de la fa\c{c}on suivante. Soit
$$
{\Bbb P}_{k}=\mathop{\rm Proj}(k[X,Y,Z])
$$
le plan projectif pond\'{e}r\'{e} o\`{u} $\mathop{\rm deg}X=n$,
$\mathop{\rm deg}Y=m$ et $\mathop{\rm deg}Z=1$, et soit $C_{k}\subset
{\Bbb P}_{k}$ le ferm\'{e} d\'{e}fini par l'\'{e}quation homog\`{e}ne
$$
F(X,Y,Z)=X^{m}-Y^{n}=0
$$
de degr\'{e} $mn$.  L'intersection de $C_{k}$ avec la carte affine
$\{Z\not=0\}=\mathop{\rm Spec}(k[x,y])$ de ${\Bbb P}_{k}$ est la courbe
d'\'{e}quation
$$
f=x^{m}-y^{n}=0;
$$
elle est donc lisse en dehors de l'origine $(0,0)$ et admet le germe
formel voulu en $(0,0)$.  L'intersection de $C_{k}$ avec le diviseur
\`{a} l'infini $\{Z=0\}$ de ${\Bbb P}_{k}$ est r\'{e}duite au point de
coordonn\'{e}es homog\`{e}nes $(1;1;0)$.  De plus le germe formel de
$C_{k}$ en ce point est isomorphe \`{a} celui en $(x=1,z=0)$ de la
courbe d'\'{e}quation
$$
x^{m}-1=0
$$
dans le plan affine $\mathop{\rm Spec}(k[x,z])$.  (La carte affine
$\{Y\not=0\}$ de ${\Bbb P}_{k}$ est le quotient de ce plan affine par le
groupe fini des racines $m$-\`{e}me de l'unit\'{e} dans $k$ pour
l'action d\'{e}finie par $(\zeta ,(x,z))\mapsto (\zeta^{n}x,\zeta z)$
et cette action est libre en dehors de l'origine.)  Par suite, $C_{k}$
est lisse au point $(1;1;0)$.  Enfin, la normalisation de $C_{k}$ est
donn\'{e}e par
$$
{\Bbb P}_{k}^{1}\rightarrow C,~(T;U)\mapsto (T^{n};T^{m};U).
$$

L'espace tangent au foncteur des d\'{e}formations plates du germe
formel de courbe $\mathop{\rm Spf}(k[[x,y]]/(f))$ est \'{e}gal \`{a}
$$
k[[x,y]]/(f,\partial_{x}f,\partial_{y}f)=k[[x,y]]/(x^{m-1},y^{n-1})
$$
et est donc de dimension
$$
\mu =(m-1)(n-1).
$$
Une d\'{e}formation miniverselle de $f$ est donn\'{e}e par
$$
\mathop{\rm Spf}(k[[\ldots ,a_{ij},\ldots ,x,y]]/(\widetilde{f}))
$$
o\`{u} les $a_{ij}$ sont des ind\'{e}termin\'{e}es sur le corps de
base $k$ et
$$
\widetilde{f}(x,y)=x^{m}-y^{n}+\sum_{{\scriptstyle 0\leq i\leq
m-2\atop\scriptstyle 0\leq j\leq n-2}}a_{ij}x^{i}y^{j}.
$$

Le conducteur ${\frak a}$ de $A=k[[t^{n},t^{m}]]$ dans son
normalis\'{e} $\widetilde{A}=k[[t]]$ est \'{e}gal \`{a}
$t^{(m-1)(n-1)}k[[t]]$, soit encore \`{a}
$$
{\frak a}=\{\sum_{i,j\geq 0}a_{ij}x^{i}y^{j}\in k[[x,y]]\mid
a_{ij}=0,~\forall i,j\hbox{ {\rm tels que }}in+jm<(m-1)(n-1)\}/(f).
$$
En particulier, le quotient
$$
V(A)={\frak a}/(\partial_{x}f,\partial_{y}f)\subset
A/(\partial_{x}f,\partial_{y}f)=k[[x,y]]/(x^{m-1},y^{n-1})
$$
admet pour base les classes des mon\^{o}mes $x^{i}y^{j}$ pour lesquels
$0\leq i\leq m-2$, $0\leq j\leq n-2$ et $in+jm\geq (m-1)(n-1)$.

Consid\'{e}rons l'espace affine $S=\mathop{\rm Spec}(k[\ldots,
a_{ij},\ldots ])$, son ferm\'{e} $T\subset S$ d\'{e}fini par les
\'{e}quations $a_{ij}=0,~\forall i,j\hbox{ tels que
}in+jm\geq (m-1)(n-1)$, le plan projectif pond\'{e}r\'{e}
$$
{\Bbb P}_{T}=\mathop{\rm Proj}({\cal O}_{T}[X,Y,Z])
$$
o\`{u} $\mathop{\rm deg}X=n$, $\mathop{\rm deg}Y=m$ et $\mathop{\rm
deg}Z=1$, et la courbe projective et plate relative
$$
\matrix{C_{T}&\kern -2mm\hookrightarrow \kern -2mm&{\Bbb P}_{T}\cr
\noalign{\smallskip}
\downarrow &\swarrow&\cr
\noalign{\smallskip}
T&&\cr}
$$
d\'{e}finie par l'\'{e}quation homog\`{e}ne
$$
\widetilde{F}(X,Y,Z)=X^{m}-Y^{n}+\kern -4mm\sum_{{{\scriptstyle 0\leq i\leq
m-2\atop \scriptstyle 0\leq j\leq n-2}\atop\scriptstyle in+jm<
(m-1)(n-1)}}\kern -4mm a_{ij}X^{i}Y^{j}Z^{mn-in-jm}=0
$$
de degr\'{e} $mn$. Pour chaque $t\in T$, l'intersection de la fibre
$C_{t}$ de $C_{T}\rightarrow T$ en $t$ avec le diviseur \`{a} l'infini
$Z=0$ du plan projectif pond\'{e}r\'{e} ${\Bbb P}_{t}$ est r\'{e}duite au
point $(1;1;0)$ et on voit comme ci-dessus que $C_{t}$ est lisse en ce
point.

On a des actions compatibles du groupe multiplicatif ${\Bbb G}_{{\rm
m},k}$ sur $T$ et $C_{T}$ donn\'{e}es par
$$
\lambda\cdot (\ldots ,a_{ij},\ldots )=(\ldots
,\lambda^{mn-in-jm}a_{ij},\ldots )
$$
et
$$
\lambda\cdot (\ldots ,a_{ij},\ldots ,X;Y;Z))=(\ldots
,\lambda^{mn-in-jm}a_{ij},\ldots ,\lambda^{n}X,\lambda^{m}Y,Z).
$$
Par d\'{e}finition du ferm\'{e} $T$ de $S$, l'action sur $T$ est contractante.

Consid\'{e}rons la composante $\mathop{\overline{\rm
Pic}}\nolimits_{C_{T}/T}^{\,0}$ de degr\'{e} $0$ du sch\'{e}ma de
Picard compactifi\'{e} relatif de $C_{T}$ sur $T$.  L'action de ${\Bbb
G}_{{\rm m},k}$ sur $C_{T}$ induit une action de ${\Bbb G}_{{\rm
m},k}$ sur $\mathop{\overline{\rm Pic}}\nolimits_{C_{T}/T}^{0}$ qui
rel\`{e}ve celle sur $T$.

L'espace tangent au ferm\'{e} $T\subset S$ est un suppl\'{e}mentaire
du sous-espace $V(A)$ de $A/(\partial_{x}f,\partial_{y}f)
=T_{(0,0)}S$.  Il s'en suit, d'apr\`{e}s le th\'{e}or\`{e}me 4.1.5,
que le sch\'{e}ma $\mathop{\overline{\rm Pic}}\nolimits_{C_{T}/T}$ est
lisse sur le corps de base $k$ le long de sa fibre
$\mathop{\overline{\rm Pic}}\nolimits_{C/k}$ \`{a} l'origine $0$ de
$T$.  Compte tenu de l'action de ${\Bbb G}_{{\rm m},k}$,
$\mathop{\overline{\rm Pic}}\nolimits_{C_{T}/T}$ est partout lisse sur
$k$.

Soient $\ell$ un nombre premier distinct de la caract\'{e}ristique de
$k$ et $K\in\mathop{\rm ob}D_{{\rm c}}^{{\rm b}}(T,{\Bbb Q}_{\ell})$
l'image directe du faisceau constant ${\Bbb Q}_{\ell}$ par la
projection $\mathop{\overline{\rm Pic}}\nolimits_{C_{T}/T}^{\,0}
\rightarrow T$. L'action de ${\Bbb G}_{{\rm m},k}$ sur $T$ se
rel\`{e}ve en une {\og}{action}{\fg} sur $K$ et on peut appliquer la
Proposition 1 de [Sp] qui assure que
$$
H^{i}(T,K)=\openup\jot\cases{K_{0} & si $i=0$,\cr
(0) & sinon.\cr}
$$
Or $R\Gamma (T,K)$ n'est autre que la cohomologie $\ell$-adique
$R\Gamma (\mathop{\overline{\rm Pic}}\nolimits_{C_{T}/T}^{\,0},{\Bbb
Q}_{\ell})$ du $k$-sch\'{e}ma lisse et quasi-projectif
$\mathop{\overline{\rm Pic}}\nolimits_{C_{T}/T}^{\,0}$, et le
th\'{e}or\`{e}me de changement de base propre implique que la fibre
$K_{0}$ de $K$ \`{a} l'origine de $T$ n'est autre que la cohomologie
$\ell$-adique $R\Gamma (\mathop{\overline{\rm
Pic}}\nolimits_{C_{k}/k}^{\,0},{\Bbb Q}_{\ell})$ du $k$-sch\'{e}ma
projectif $\mathop{\overline{\rm Pic}}\nolimits_{C_{k}/k}^{\,0}$, soit
encore la cohomologie $\ell$-adique $R\Gamma (X^{0},{\Bbb Q}_{\ell})$
de la composante de degr\'{e} $0$ de la fibre de Springer puisque
cette composante $X^{0}=Z^{0}$ est universellement hom\'{e}omorphe
\`{a} $\mathop{\overline{\rm Pic}}\nolimits_{C/k}^{0}$.

Si $k$ est de caract\'{e}ristique $p>0$, toute la situation ci-dessus
est d\'{e}finie sur ${\Bbb F}_{p}$.  Appliquant alors la forme forme
de la conjecture de Weil prouv\'{e}e par Deligne comme l'a fait
Springer dans [Sp] , on d\'{e}duit de qui pr\'{e}c\`{e}de:

\thm PROPOSITION 4.2.1
\enonce
Chaque groupe de cohomologie $H^{i}(X^{0},{\Bbb Q}_{\ell})$ est pur de
poids $i$.
\hfill\hfill$\square$
\endthm

Bien s\^{u}r, cette proposition r\'{e}sulte aussi de fait que $X^{0}$
peut \^{e}tre pav\'{e} par des espaces affines (cf.  [Lu-Sm]).

\section{4.3}{Construction de d\'{e}formations de fibres de Springer
affine}

Revenons \`{a} la situation g\'{e}n\'{e}rale de la section 4.1.  Nous
avons donc la courbe int\`{e}gre, projective $C_{s}=C_{I,s_{I}}$ sur
$k$, de normalis\'{e}e $\pi_{s}:\widetilde{C}_{s}\rightarrow C_{s}$
une droite projective, avec son unique point singulier $c$ en lequel
la singularit\'{e} est plane et l'ensemble des branches
$\pi_{s}^{-1}(c)=\{\widetilde{c}_{i}\mid i\in I\}$ est index\'{e} par
$I$.  Nous allons utiliser les r\'{e}sultats g\'{e}n\'{e}raux de la
section 2.5 et les r\'{e}sultats de la section 2.6 pour
{\og}{d\'{e}former}{\fg} le rev\^{e}tement
$$
X_{s}\rightarrow Z_{s}
$$
de groupe de Galois
$$
\Lambda_{s}^{0}=H^{0}(\widetilde{C}_{s}\setminus\pi_{s}^{-1}(c),
{\Bbb G}_{{\rm m}})/H^{0}(\widetilde{C}_{s},{\Bbb G}_{{\rm m}})
=\mathop{\rm Ker}({\Bbb Z}^{I}\rightarrow {\Bbb Z}),
$$
de la section $2.2$. Plus exactement, nous d\'{e}formerons le rev\^{e}tement
$$
\overline{P}_{s}^{\,\natural}\rightarrow\overline{P}_{s}
$$
qui lui est hom\'{e}omorphe.
\vskip 2mm

Consid\'{e}rons l'alg\'{e}brisation $(C=C_{I},\infty )\rightarrow
S_{I}=S$ de la d\'{e}formation miniverselle de $(C_{s},\infty_{s})$
introduite en 4.1 et la strate \`{a} $\delta$ constant
$S^{\delta}=S_{I}^{\delta}\subset S$ (cf.  la fin de la section
$3.4$).  Consid\'{e}rons aussi le morphisme de normalisation
$\pi_{S^{\delta}}:\widetilde{S^{\delta}}=\widetilde{S_{I}^{\delta}}
\rightarrow S^{\delta}$ de cette strate et la courbe
$$
C_{\widetilde{S^{\delta}}}=C_{I,\widetilde{S^{\delta}}}\rightarrow
\widetilde{S^{\delta}}
$$
d\'{e}duite de $C\rightarrow S$ par le changement de base
$$
\widetilde{S^{\delta}}\twoheadrightarrow S^{\delta}\hookrightarrow S.
$$
D'apr\`{e}s le Corollaire 3.3.4, cette courbe relative admet une
normalisation en famille $\widetilde{C_{\widetilde{S^{\delta}}}}
\rightarrow C_{\widetilde{S^{\delta}}}$, et on a en fait
$$
\widetilde{C_{\widetilde{S^{\delta}}}}\cong {\Bbb
P}_{\widetilde{S^{\delta}}}^{1}
$$
puisque la normalis\'{e}e de $C_{s}$ est une droite projective sur
$k$.  On peut choisir l'isomorphisme ci-dessus de telle sorte que la
section de $\widetilde{S^{\delta}}\rightarrow
\widetilde{C_{\widetilde{S^{\delta}}}}$ induite par la section $\infty
:S\rightarrow C$ corresponde \`{a} la section \`{a} l'infini de ${\Bbb
P}_{\widetilde{S^{\delta}}}^{1}$ sur $\widetilde{S^{\delta}}$.

On a vu (cf.  Th\'{e}or\`{e}me 3.4.11) que, au moins pour $k$ de
caract\'{e}ristique nulle, le morphisme $C_{\widetilde{S^{\delta}}}
\rightarrow\widetilde{S^{\delta}}$ a les propri\'{e}t\'{e}s
suivantes :

\decale{\rm (1)} $\widetilde{S^{\delta}}$ est un sch\'{e}ma
strictement local r\'{e}gulier (formellement lisse sur $k$) de
dimension $\tau (C_{s},\infty_{s})-\delta_{I}$,

\decale{\rm (3)} toutes les fibres g\'{e}om\'{e}triques de
$\widetilde{C_{\widetilde{S^{\delta}}}}\rightarrow
\widetilde{S^{\delta}}$ sont des courbes int\`{e}gres \`{a}
singularit\'{e}s planes,

\decale{\rm (2)} la fibre g\'{e}n\'{e}rique g\'{e}om\'{e}trique de
$\widetilde{C_{\widetilde{S^{\delta}}}}\rightarrow
\widetilde{S^{\delta}}$ n'a comme seules singularit\'{e}s que
$\delta_{I}$ points doubles ordinaires.

\thm HYPOTH\`{E}SES{\pc ~ET} {\pc NOTATIONS} 4.3.1
\enonce
{\it Dans la suite, nous supposerons que les propri\'{e}t\'{e}s $(1)$
\`{a} $(3)$ sont v\'{e}rifi\'{e}es sur notre corps $k$ {\rm (}c'est le
cas par exemple si la carac\-t\'{e}ristique de $k$ nulle ou {\og}{assez
grande}{\fg}{\rm )}.

De plus, comme la d\'{e}formation totale $C\rightarrow S$
n'interviendra plus, nous noterons simplement $C\rightarrow S$ la
courbe relative $C_{\widetilde{S^{\delta}}}\rightarrow
\widetilde{S^{\delta}}$ pour all\'{e}ger l'exposition.}
\endthm

Nous avons donc un sch\'{e}ma strictement local $S$ formellement lisse
sur $k$, de dimension $\tau (C_{s},\infty_{s})-\delta_{I}$, et une
courbe relative $C$ sur $S$, munie d'une section globale $\infty$ le
long de laquelle $C$ est lisse sur $S$, de fibre sp\'{e}ciale $C_{s}$,
qui admet une normalisation en famille $\pi_{C}:\widetilde{C}
\rightarrow C$ dont l'espace total $\widetilde{C}$ est une droite
projective sur $S$.  Le morphisme structural $f:C\rightarrow S$ est
lisse en dehors du sous-sch\'{e}ma $D\subset C$ fini sur $S$ qui est
d\'{e}fini par l'Id\'{e}al conducteur ${\frak a}$ de $\pi_{\ast}{\cal
O}_{\widetilde{C}}$ dans ${\cal O}_{C}$ et
$$
\widetilde{f}=f\circ\pi_{C}:\widetilde{C}\rightarrow S
$$
est identifi\'{e} \`{a} la droite projective standard ${\Bbb
P}_{S}^{1}\rightarrow S$ de telle sorte que la section \`{a} l'infini
corresponde \`{a} la section $\infty :S\rightarrow C$, et donc ne
rencontre pas $\widetilde{D}=\pi_{C}^{-1}(D)$.

Consid\'{e}rons les $S$-sch\'{e}mas de Picard $P$ et de Picard
compactifi\'{e} $\overline{P}$ relatifs de $C$ sur $S$.  Le
$S$-sch\'{e}ma en groupes $P$ est isomorphe \`{a} $P^{0}\times {\Bbb
Z}$ puisque le lieu de lissit\'{e} de $C\rightarrow S$ admet une
section.  Pour chaque entier $d$, la composante connexe
$\overline{P}^{d}$ contient $P^{d}=P^{0}\times\{d\}$ comme ouvert
dense.

La composante neutre $P^{0}$ du sch\'{e}ma de Picard est par
d\'{e}finition le $S$-sch\'{e}ma affine et lisse qui repr\'{e}sente le
faisceau fppf
$$
\widetilde{f}_{\ast}{\Bbb G}_{{\rm m},\widetilde{C}}/f_{\ast}{\Bbb
G}_{{\rm m},C}=f_{\ast}(\pi_{C,\ast}{\Bbb G}_{{\rm m},\widetilde{C}}/
{\Bbb G}_{{\rm m},C}).
$$
Notons $f_{D}$ {\rm (}resp.  $\widetilde{f}_{D}${\rm )} les
restrictions de $f$ et $\widetilde{f}$ aux ferm\'{e}s $D\subset C$
{\rm (}resp.  $\widetilde{D}\subset \widetilde{C}${\rm )}.  On a le
$S$-sch\'{e}ma en groupes affine et lisse
$$
\mathop{\rm Res}\nolimits_{D/S}{\Bbb G}_{{\rm m}}\hbox{ {\rm (}resp.
}\mathop{\rm Res}\nolimits_{\widetilde{D}/S}{\Bbb G}_{{\rm m}}{\rm )}
$$
restriction \`{a} la Weil du groupe multiplicatif de $D$ \`{a} $S$
{\rm (}resp.  de $\widetilde{D}$ \`{a} $S${\rm )}, qui repr\'{e}sente
le faisceau fppf $f_{D,\ast}{\Bbb G}_{{\rm m},D}$ {\rm (}resp.
$\widetilde{f}_{D,\ast}{\Bbb G}_{{\rm m}, \widetilde{D}}${\rm )}, et
on a la fl\`{e}che d'adjonction
$$
\alpha :\pi_{C,\ast}{\Bbb G}_{{\rm m},\widetilde{C}}/{\Bbb G}_{{\rm
m}, C}\rightarrow i_{\ast}(\pi_{D,\ast}{\Bbb G}_{{\rm m},
\widetilde{D}}/{\Bbb G}_{{\rm m},D})
$$
o\`{u} $i:D\hookrightarrow C$ est l'inclusion et $\pi_{D}:
\widetilde{D}\rightarrow D$ est la restriction de $\pi$.

\thm LEMME 4.3.2
\enonce
La fl\`{e}che d'adjonction $\alpha$ est un isomorphisme et induit un
isomorphisme de $S$-sch\'{e}mas en groupes
$$
P^{0}\buildrel\sim\over\longrightarrow
\mathop{\rm Res}\nolimits_{\widetilde{D}/S}{\Bbb G}_{{\rm m}}/
\mathop{\rm Res}\nolimits_{D/S}{\Bbb G}_{{\rm m}}.
$$
\endthm

\rem Preuve
\endrem
Soit $\mathop{\rm Spec}(A)\rightarrow C$ une carte affine de $C$.  Au
dessus cette carte, $\widetilde{C}$ est \'{e}gal \`{a} $\mathop{\rm
Spec}(B)$ pour une $A$-alg\`{e}bre finie $B$, la trace de $D$ est
d\'{e}finie par l'id\'{e}al conducteur $I$ de $B$ dans $A$ et $\alpha$
est donn\'{e}e par la fl\`{e}che naturelle
$$
B^{\times}/A^{\times}\rightarrow (B/I)^{\times}/(A/I)^{\times}
$$
(on rappelle que $I\subset A\subset B$ est \`{a} la fois un id\'{e}al
de $A$ et un id\'{e}al de $B$). Comme cette derni\`{e}re fl\`{e}che
est trivialement injective, l'injectivit\'{e} de $\alpha$ est
d\'{e}montr\'{e}e.

Maintenant, si $b,b'\in B$ sont tels que $bb'=1+a$ avec $a\in I\subset
A$, quitte \`{a} remplacer $\mathop{\rm Spec}(A)$ par le voisinage
ouvert $\mathop{\rm Spec}(A[(1+a)^{-1}])$ du ferm\'{e} $\mathop{\rm
Spec}(A/I)\subset \mathop{\rm Spec}(A)$, on voit que $b$ est
inversible dans $B$, ce qui d\'{e}montre la surjectivit\'{e} de
$\alpha$.
\hfill\hfill$\square$
\vskip 3mm

La fibre sp\'{e}ciale $\mathop{\rm Res}\nolimits_{D_{s}/s}{\Bbb
G}_{{\rm m}}$ du $S$-sch\'{e}ma en groupes $\mathop{\rm
Res}\nolimits_{D/S}{\Bbb G}_{{\rm m}}$ admet pour tore maximal le tore
${\Bbb G}_{{\rm m},k}$ puisque $(D_{s})_{{\rm red}}$ est r\'{e}duit au
point $c\in D_{s}\subset C_{s}$.  De m\^{e}me, $\mathop{\rm
Res}\nolimits_{\widetilde{D}_{s}/s}{\Bbb G}_{{\rm m}}$ admet pour tore
maximal le tore ${\Bbb G}_{{\rm m},k}^{I}$ puisque
$(\widetilde{D}_{s})_{{\rm red}}=\{\widetilde{c}_{i}\mid i\in I\}$. Le
tore maximal $T_{s}$ de $P_{s}^{0}$ est donc canoniquement isomorphe
\`{a} ${\Bbb G}_{{\rm m},k}^{I}/{\Bbb G}_{{\rm m},k}$ (plongement
diagonal).

Le tore maximal ${\Bbb G}_{{\rm m},k}$ de $\mathop{\rm Res}
\nolimits_{D_{s}/s}{\Bbb G}_{{\rm m}}$ est la fibre sp\'{e}ciale du
tore canonique ${\Bbb G}_{{\rm m},S}\subset \mathop{\rm
Res}\nolimits_{D/S}{\Bbb G}_{{\rm m}}$ d\'{e}fini par la fl\`{e}che
d'adjonction $\mathop{\rm id}\rightarrow f_{D,\ast}f_{D}^{\ast}$.
Comme $S$ est strictement hens\'{e}lien, $\widetilde{D}$ se casse en
autant de composantes connexes qu'il y a de points dans
$(\widetilde{D}_{s})_{{\rm red}}$ et le tore maximal ${\Bbb G}_{{\rm
m},k}^{I}$ de $\mathop{\rm Res} \nolimits_{\widetilde{D}_{s}/s}{\Bbb
G}_{{\rm m}}$ est aussi la fibre sp\'{e}ciale du tore canonique ${\Bbb
G}_{{\rm m},S}^{I}\subset \mathop{\rm Res}
\nolimits_{\widetilde{D}/S}{\Bbb G}_{{\rm m}}$ d\'{e}fini par la
fl\`{e}che d'adjonction $\mathop{\rm id}\rightarrow
\widetilde{f}_{D,\ast}\widetilde{f}_{D}^{\ast}$, composante connexe
par composante connexe de $\widetilde{D}$.

Par suite, le tore maximal $T_{s}$ de $P_{s}^{0}$ est la fibre
sp\'{e}ciale d'un tore canonique
$$
T={\Bbb G}_{{\rm m},S}^{I}/{\Bbb G}_{{\rm m},S}\subset P^{0}.
$$
Plus g\'{e}n\'{e}ralement, pour chaque point g\'{e}om\'{e}trique $t$
de $S$, la fibre $T_{t}$ de $T$ en $t$ est contenue dans le tore
maximal de
$$
P_{t}^{0}\buildrel\sim\over\longrightarrow
\mathop{\rm Res}\nolimits_{\widetilde{D}_{t}/t}{\Bbb G}_{{\rm m}}/
\mathop{\rm Res}\nolimits_{D_{t}/t}{\Bbb G}_{{\rm m}},
$$
tore maximal qui est isomorphe \`{a} $\prod_{j\in J_{t}}({\Bbb
G}_{{\rm m},\kappa (t)}^{I_{t,j}}/{\Bbb G}_{{\rm m},t})$ o\`{u}
$\{c_{j,t}\mid j\in J_{t}\}$ est l'ensemble des points singuliers de
$C_{t}$ et, pour chaque $j\in J_{t}$, $\{\widetilde{c}_{t,i}\mid i\in
I_{t,j}\}$ est l'ensemble des branches du germe formel de $C_{t}$ en
son point singulier $c_{t,j}$.

On identifie le groupe des caract\`{e}res de $T_{s}$, ou ce qui
revient au m\^{e}me celui de $T$, \`{a} $\Lambda_{s}^{0}:=\Lambda^{0}$
comme dans la section $2.6$.  Pour chaque point g\'{e}om\'{e}trique
$t$ de $S$, $\Lambda^{0}$ est donc un quotient du groupe des
caract\`{e}res
$$
\Lambda_{t}^{0}=\prod_{j\in J}\mathop{\rm Ker}({\Bbb Z}^{I_{t,j}}
\rightarrow {\Bbb Z})
$$
du tore maximal de $P_{t}^{0}$, quotient que l'on peut voir
concr\`{e}tement en consid\'{e}rant la mani\`{e}re dont les points
singuliers $c_{t,j}$ et leurs branches $\widetilde{c}_{t,i}$ confluent
vers le point singulier $c$ et ses branches $\widetilde{c}_{i}$.
\vskip 2mm

Nous sommes maintenant en mesure de construire la d\'{e}formation de
$\overline{P}_{s}^{\,\natural}\rightarrow\overline{P}_{s}$.

Pour chaque entier $d$, la restriction de l'homomorphisme $\beta$
d'Esteves, Gagn\'{e} et Kleiman \`{a} $T$ est un homomorphisme
$T\rightarrow\mathop{\rm Pic}\nolimits_{\overline{P}/S}$ qui
d\'{e}finit d'apr\`{e}s la proposition 2.5.1 un $\Lambda^{0}$-torseur
$$
Q^{d}\rightarrow\overline{P}{}^{d}
$$
dont la formation commute \`{a} tout changement de base $S'\rightarrow
S$.  Notons
$$
Q\rightarrow\overline{P}
$$
la somme disjointe de ces torseurs. C'est la d\'{e}formation cherch\'{e}e.

En effet, la fibre sp\'{e}ciale est par d\'{e}finition le
rev\^{e}tement $\overline{P}_{s}^{\,\natural}\rightarrow
\overline{P}_{s}$.  Plus g\'{e}n\'{e}ralement, pour chaque point
g\'{e}om\'{e}trique $t$ de $S$, $Q_{t}\rightarrow \overline{P}_{t}$
est le quotient du rev\^{e}tement $\overline{P}_{t}^{\,\natural}
\rightarrow\overline{P}_{t}$ dans le corollaire $2.2.6$, quotient qui
correspond au quotient $\Lambda_{t}^{0}\twoheadrightarrow \Lambda^{0}$
entre les groupes de Galois.
\vskip 5mm

\centerline{5. BIBLIOGRAPHIE}
\vskip 5mm

\newtoks\ref \newtoks\auteur \newtoks\titre \newtoks\annee
\newtoks\revue \newtoks\tome \newtoks\pages \newtoks\reste

\def\bibitem#1{\parindent=20pt\itemitem{#1}\parindent=24pt}

\def\article{\bibitem{[\the\ref]}%
\the\auteur~-- \the\titre, {\sl\the\revue} {\bf\the\tome},
({\the\annee}), \the\pages.\smallskip\filbreak}

\def\autre{\bibitem{[\the\ref]}%
\the\auteur~-- \the\reste.\smallskip\filbreak}

\ref={A-I-K}
\auteur={A. {\pc ALTMAN}, A. {\pc IARROBINO}, S. {\pc KLEIMAN}}
\reste={Irreducibility of the Compactified Jacobian, dans ``Real and
complex singularities. Proceedings, Oslo 1976, P. Holm (ed.)'',
Sijthoff \& Nordhoff, (1977), 1-12}
\autre

\ref={A-K 1}
\auteur={A. {\pc ALTMAN}, S. {\pc KLEIMAN}}
\reste={Introduction to Grothendieck Duality Theory, Lecture Notes in
Math. {\bf 146}, Springer-Verlag, (1970)}
\autre

\ref={A-K 2}
\auteur={A. {\pc ALTMAN}, S. {\pc KLEIMAN}}
\titre={Compactifying the Jacobian}
\revue={Bull. Am. Math. Soc.}
\tome={82}
\annee={1976}
\pages={947-949}
\article

\ref={A-K 3}
\auteur={A. {\pc ALTMAN}, S. {\pc KLEIMAN}}
\titre={Compactifying the Picard Scheme}
\revue={Advances in Math.}
\tome={35}
\annee={1980}
\pages={50-112}
\article

\ref={A-K 4}
\auteur={A. {\pc ALTMAN}, S. {\pc KLEIMAN}}
\titre={Bertini Theorems for Hypersurface Sections Containing a
Subscheme}
\revue={Communication in Algebra}
\tome={7}
\annee={1979}
\pages={775-790}
\article

\ref={Be}
\auteur={R. {\pc BEZRUKAVNIKOV}}
\titre={The Dimension of the Fixed Point Set on Affine Flag Manifolds}
\revue={Math. Res. Letters}
\tome={3}
\annee={1996}
\pages={185-189}
\article

\ref={B-L-K}
\auteur={D. {\pc BOSCH}, W. {\pc L\"{U}TKEBOHMERT}, M. {\pc RAYNAUD}}
\reste={N\'{e}ron Models, Ergebnisse der Mathematik und ihrer
Grenzgebiete 3.Folge, Band 21, Springer-Verlag, (1990)}
\autre

\ref={De 1}
\auteur={P. {\pc DELIGNE}}
\reste={La classe de cohomologie associ\'{e}e \`{a} un cycle, dans
{\it Cohomologie \'{e}tale, SGA 4${1\over 2}$}, Lecture Notes in Math.
{\bf 569}, Springer-Verlag, (1977), 129-153}
\autre

\ref={De 2}
\auteur={P. {\pc DELIGNE}}
\reste={Intersection sur les surfaces r\'{e}guli\`{e}res, dans ``SGA 7
II, Groupe de Monodromie en G\'{e}om\'{e}trie Alg\'{e}brique,
dirig\'{e} par P. Deligne et N. Katz'', Lecture Notes in Math.  {\bf
340}, Springer-Verlag, (1973), 1-38}
\autre

\ref={D-H}
\auteur={S. {\pc DIAZ}, J. {\pc HARRIS}}
\titre={Ideals Associated to Deformations of Singular Plane Curves}
\revue={Trans. of the Amer. Math. Soc.}
\tome={309}
\annee={1988}
\pages={433-468}
\article

\ref={E-G-K}
\auteur={E. {\pc ESTEVES}, M. {\pc GAGN\'{E}}, S. {\pc KLEIMAN}}
\reste={Autoduality of the compactified Jacobian,
http://arxiv.org/abs/math.AG//9911071, (1999)}
\autre

\ref={Fo}
\auteur={J. {\pc FOGARTY}}
\titre={Algebraic Families on an Algebraic Surface}
\revue={Amer. J. Math.}
\tome={90}
\annee={1968}
\pages={511-521}
\article

\ref={F-G-S}
\auteur={D. {\pc FANTECHI}, L. {\pc G\"{O}TTSCHE}, D. {\pc VAN STRATEN}}
\titre={Euler Number of the Compactified Jacobian and Multiplicity of
Rational Curves}
\revue={J. Algebraic Geometry}
\tome={8}
\annee={1999}
\pages={115-133}
\article

\ref={Gr 1}
\auteur={A. {\pc GROTHENDIECK}}
\reste={Formule de Lefschetz et rationalit\'{e} des fonctions $L$,
S\'{e}minaire Bourbaki 279, 1964/65, dans {\it Dix expos\'{e}s sur la
cohomlologie des sch\'{e}mas}, Masson, Paris, et North-Holland,
Amsterdam, (1968)}
\autre

\ref={Gr 2}
\auteur={A. {\pc GROTHENDIECK}}
\reste={Fondements de la g\'{e}om\'{e}trie alg\'{e}brique,
S\'{e}minaire Bourbaki 232, Benjamin, New York, (1966)}
\autre

\ref={Gr 3}
\auteur={A. {\pc GROTHENDIECK}}
\reste={\'{E}l\'{e}ments de g\'{e}om\'{e}trie alg\'{e}brique, III,
\'{E}tude cohomologique des faisceaux coh\'{e}rents, {\it Publications
Math\'{e}\-ma\-ti\-ques de l'I.H.\'{E}.S.} {\bf 11}, (1961)}
\autre

\ref={Gr 4}
\auteur={A. {\pc GROTHENDIECK}}
\reste={\'{E}l\'{e}ments de g\'{e}om\'{e}trie alg\'{e}brique, IV,
\'{E}tude locale des sch\'{e}mas et des morphismes de sch\'{e}mas
(Troisi\`{e}me partie), {\it Publications Math\'{e}\-ma\-ti\-ques de
l'I.H.\'{E}.S.} {\bf 28}, (1966)}
\autre

\ref={Ia}
\auteur={A. {\pc IARROBINO}}
\reste={Punctual Hilbert Schemes, Mem. Am. Math. Soc. {\bf 188},
(1977)}
\autre

\ref={Il}
\auteur={L. {\pc ILLUSIE}}
\reste={Complexe Cotangent et D\'{e}formations I, Lecture Notes in
Math. {\bf 239}, Springer-Verlag, (1971)}
\autre

\ref={K-L}
\auteur={D {\pc KAZHDAN}, G. {\pc LUSZTIG}}
\titre={Fixed Point Varieties on Affine Flag Manifolds}
\revue={Israel J. of Math.}
\tome={62}
\annee={1988}
\pages={129-168}
\article

\ref={L-R}
\auteur={G. {\pc LAUMON}, J. {\pc RAPOPORT}}
\reste={A geometric approach to the fundamental lemma for unitary
groups, http://arxiv.org/abs/math.AG/9711021, (1997)}
\autre

\ref={La-Sh}
\auteur={R. {\pc LANGLANDS}, D. {\pc SHELSTAD}}
\titre={On the definition of the transfer factors}
\revue={Math. Ann.}
\tome={278}
\annee={1987}
\pages={219-271}
\article

\ref={Lu-Sm}
\auteur={G. {\pc LUSZTIG}, J.M. {\pc SMELT}}
\titre={Fixed point varieties on the space of lattices}
\revue={Bull. London. Math. Soc.}
\tome={23}
\annee={1991}
\pages={213-218}
\article

\ref={M-F}
\auteur={D. {\pc MUMFORD}, J. {\pc FOGARTY}}
\reste={Geometric Invariant Theory, Second Enlarged Edition,
Springer-Verlag, (1982)}
\autre

\ref={Ra}
\auteur={M {\pc RAYNAUD}}
\titre={Sp\'{e}cialisation du foncteur de Picard}
\revue={Publications Math\'{e}\-ma\-ti\-ques de l'I.H.\'{E}.S.}
\tome={38}
\annee={1970}
\pages={27-76}
\article

\ref={Re}
\auteur={C.J. {\pc REGO}}
\titre={The Compactified Jacobian}
\revue={Ann. Scient. \'{E}c. Norm. Sup.}
\tome={13}
\annee={1980}
\pages={211-223}
\article

\ref={Ri}
\auteur={D.S. {\pc RIM}}
\reste={Formal Deformation Theory, dans ``SGA 7 I, Groupe de
Monodromie en G\'{e}om\'{e}trie Alg\'{e}brique, dirig\'{e} par A.
Grothendieck'', Lecture Notes in Math.  {\bf 288}, Springer-Verlag,
(1972), 32-132}
\autre

\ref={Sp}
\auteur={T.A. {\pc SPRINGER}}
\titre={A purity result for fixed point varieties in flag manifolds}
\revue={J. Fac. Sci. Univ. Tokyo}
\tome={31}
\annee={1984}
\pages={271-282}
\article

\ref={Te 1}
\auteur={B. {\pc TEISSIER}}
\reste={Cycles \'{e}vanescents, sections hyperplanes et condition de
Whitney, dans ``Singularit\'{e}s \`{a} Carg\`{e}se'', Ast\'{e}risque
7 et 8, (1973), 285-362}
\autre

\ref={Te 2}
\auteur={B. {\pc TEISSIER}}
\reste={The hunting of invariants in the geometry of discriminants,
dans ``Real and complex singularities.  Proceedings, Oslo 1976, P.
Holm (ed.)'', Sijthoff \& Nordhoff, (1977), 565-677}
\autre

\ref={Te 3}
\auteur={B. {\pc TEISSIER}}
\reste={R\'{e}solution simultan\'{e}e - I, II, dans ``S\'{e}minaire
sur les Singularit\'{e}s des Surfaces, Palaiseau, France 1976-1977,
\'{E}dit\'{e} par M. Demazure, H. Pinkham, B. Teissier'',
Lecture Notes in Math. {\bf 777}, (1980), 71-146}
\autre

\ref={Wa}
\auteur={J.-L. {\pc WALDSPURGER}}
\reste={Int\'{e}grales orbitales nilpotentes et endoscopie pour les
groupes classiques non ramifi\'{e}s, Ast\'{e}risque 269, (2001)}
\autre

\bye